\magnification=\magstep1
\vbadness=10000
\hbadness=10000
\tolerance=10000

\def\B{{\bf B}}   
\def\C{{\bf C}}   
\def\d{\partial } 
\def\e{{\bf e}}   
\def\E{{\bf E}}   
\def\G{{\bf G}}   
\def\P{{\bf P}}   
\def\Q{{\bf Q}}   
\def\R{{\bf R}}   
\def\Z{{\bf Z}}   
\def\Aut{{\rm Aut}}
\def\Hom{{\rm Hom}}
\def\\{\backslash}

Please send comments, errors,  etc. to reb@dpmms.cam.ac.uk

\proclaim
Automorphic forms with singularities on Grassmannians. 
29 September 1996, corrected 29 May 1997.

Richard E. Borcherds,\footnote{$^*$}{ Supported by a Royal Society
professorship and  by NSF grant
DMS-9401186.}

D.P.M.M.S.,

16 Mill Lane,

Cambridge,

CB2 1SB,

England.

e-mail: reb@dpmms.cam.ac.uk

home page: http://www.dpmms.cam.ac.uk/\~{}reb

\bigskip

We construct some families of automorphic forms on Grassmannians which
have singularities along smaller sub Grassmannians, using  Harvey
and Moore's extension of
the Howe (or theta) correspondence to modular forms with poles at
cusps.
Some of the applications are as follows. We construct families
of holomorphic automorphic forms which can be written as infinite
products, which give many new examples of generalized Kac-Moody
superalgebras.  We extend the Shimura and Maass-Gritsenko
correspondences to modular forms with singularities. We prove some
congruences satisfied by the theta functions of positive definite
lattices, and find a sufficient condition for a Lorentzian lattice to
have a reflection group with a finite volume fundamental domain.  We
give some examples suggesting that these automorphic forms with
singularities are related to Donaldson polynomials and to mirror
symmetry for K3 surfaces.

\proclaim Contents.

1. Introduction.

Notation and terminology.

2. Modular forms.

3. Fourier transforms.

4. Siegel theta functions.

5. Reduction to smaller lattices.

6. The singularities of $\Phi$.

7. The Fourier expansion of $\Phi$.

8. Anisotropic lattices.

9. Definite lattices.

10. Lorentzian lattices.

11. Congruences for lattices.

12. Hyperbolic reflection groups.

13. Holomorphic infinite products.

14. The Shimura-Doi-Naganuma-Maass-Gritsenko-... correspondence.

15. Examples related to mirror symmetry and Donaldson polynomials.

16. Open problems.

\proclaim 1.~Introduction.

In this paper we construct automorphic functions $\Phi_M(v,p,F)$ with
known singularities on the Grassmannians $G(b^+,b^-)$ of
$b^+$-dimensional positive definite subspaces of $\R^{b^+,b^-}$.
As a special case, for $b^+=2$
we recover the results of [B95] giving examples of holomorphic
automorphic forms which can be written as infinite products.

The main tool we use is Harvey and Moore's extension of the Howe (or
theta) correspondence to automorphic forms with singularities [H-M].
We briefly recall the Howe correspondence; see the articles by Howe,
Gelbart, and Rallis in [B-C] for more details. If we have a commuting
pair of subgroups in the metaplectic group (a double cover of the
symplectic group) then we get a correspondence between representations
of the two subgroups, by decomposing the metaplectic representation of
the metaplectic group into a sum of tensor products of representations
of the two subgroups. As some representations of groups over the adeles
tend to correspond to automorphic forms, we can get a correspondence
between automorphic forms on these two groups. We will take
the commuting pair of subgroups to be the double cover $Mp_2$ of $SL_2$
and the orthogonal group $O_{b^+,b^-}$. If we unravel all the definitions
we find that the correspondence between automorphic forms
can be described explicitly as follows. To simplify slightly we will
take $b^+=2$ and
take  $M$ to be an even unimodular
lattice in $\R^{2,b^-}$ and $F$ to be a holomorphic modular form
for $SL_2(\Z)$ of weight $1-b^-/2$. (Actually we have simplified
a bit too much because the only such forms $F$ are constant, but we
will ignore this as we are only giving a rough idea of things.)
The Siegel theta function $\Theta_M(\tau;v^+)$ (defined in section 4)
is a function of
$v^+\in G(2,b^-)$ and of $\tau$ in the upper half plane and is invariant
under the action of $\Aut(M)$ on $v^+$. The function
$$\Phi_M(v,F)=\int_{\tau\in SL_2(\Z)\\H}\bar\Theta_M(\tau;v^+)F(\tau)dxdy/y
$$
is then  an automorphic function on $G(2,b^-)$ invariant under
the discrete group $O_M(\Z)$,
and this is (roughly) the Howe correspondence
from automorphic forms $F$ on $SL_2$ to automorphic forms
$\Phi_M$ on $O_{2,b^-}$ (at least for the purposes
of this paper).

Now suppose that we allow $F$ to have poles at cusps but but still
insist that it be holomorphic on the upper half plane $H$.  Then the
integral above diverges wildly. However Harvey and Moore discovered
that it is still possible to make sense of the integral by
regularizing it.  They showed that the results of [B95] could be given
much simpler proofs using this ``singular Howe correspondence'',
because the regularized integral turns out to be more or less the real
part of the logarithm of the infinite products used in [B95] to define
automorphic forms.  For example, the singularities of $\Phi_M$ can be
easily read off from the singularities of $F$, and this immediately
gives the locations of the zeros of the corresponding infinite
product.

In this paper we will generalize this construction
in the following ways.
\item{1.} We  replace $M$ by a lattice of any determinant,
or by a coset of such a lattice, and $F$ by a form of higher level.
This is similar to the level 1 case but with the usual extra
complications for higher level: we have to deal with more than one
cusp, and we often have to replace automorphic forms by finite
dimensional spaces of automorphic forms. We deal with both these
problems by using vector valued modular forms; a useful bonus of using
these is that forms of any level can be considered as vector valued
modular forms of level 1, so we immediately reduce the higher level
case to the level 1 case.
\item{2.} We   allow $M$ to be a lattice
in any space $\R^{b^+,b^-}$, with $b^+$ not necessarily 2, although
the Grassmannian is then usually no longer hermitian.
\item{3.} We  replace the Siegel theta function $\Theta_M$
by a function depending on some homogeneous polynomial $p$ on
$\R^{b^+,b^-}$. (We do not insist that this polynomial $p$ should be
harmonic.) In the case of $\R^{2,1}$ and holomorphic forms $F$ this
gives Niwa's description [Ni] of the Shimura correspondence for a
suitable choice of harmonic function. For $\R^{2,3}$ and holomorphic forms $F$
we recover Maass'
correspondence (see [E-Z]) and for $\R^{2,b^-}$ we recover the higher
dimensional generalization of Maass' correspondence due to Gritsenko
[Gr].  See theorem 14.3.
\item{4.} We allow $F$ to be an ``almost holomorphic'' form;
for example, we allow powers of the modular form $\E_2(\tau)$. (This
case was also done in [H-M].)
\item{5.} Up to section 7 we can allow $F$ to be a real analytic
modular form with singularities at cusps. (However in the later
sections we assume $F$ is almost holomorphic. See problem 16.13.)

In the cases of lattices in $\R^{2,b^-}$ we recover the results of
[B95] constructing holomorphic automorphic forms as infinite products,
but with much simpler proofs.  (We have taken advantage of the
simplifications to do everything in much greater generality, so in
fact the proofs end up looking more complicated. The reader who wants
to extract a simpler proof can take $M$ unimodular and $p=1$ in
sections 2 to 7, when a lot of the complexity vanishes; see [Kon] for
an expository account of this case.)  The main improvements are as
follows.  In [B95] the proof starts with an explicit Fourier expansion
of $\Phi$ (or rather an infinite product expansion of $\exp(\Phi)$).
We then have to analytically continue $\Phi$, find its zeros, and
check that it is an automorphic form, none of which are easy to do.
For example, to prove that $\Phi$ is an automorphic form using this
method it is necessary to find a set of generators for the discrete
groups and check by long calculations that  $\Phi$ transforms
correctly under each generator.  In particular it seems hard to
generalize the method to higher levels. (It is possible to do this in
some particular cases; for example, Gritsenko and Nikulin [G-N] have
recently used the method of [B95] to produce many higher level
examples of automorphic forms, such as a Siegel modular form of genus
2 and weight 5, which can be written as infinite products.) In Harvey
and Moore's approach used in this paper, we start with an expression
for $\Phi$ which is obviously invariant under $\Aut(M)$ and for which
it is trivial to read off the singularities of $\Phi$.  The only
problem is to calculate the Fourier expansion of $\Phi$. To do this we
use a modification of Harvey and Moore's calculation in [H-M appendix
A] in order to get a recursive formula (Theorem 7.1) relating the
Fourier series of $\Phi$ to that of an automorphic form on
$G(b^+-1,b^--1)$.  This can be thought of as a version of the
Rankin-Selberg method: we first write the theta function of a lattice
$M$ containing a norm 0 vector $z$ as a sort of Poincar\'e series
involving theta functions of lattices $z^\perp/z$, and then unravel
the integral over a fundamental domain of $SL_2(\Z)$ to get an integral
over the rectangular region $\Re(\tau)\le 1/2$.

We now describe the sections of this paper in more detail.

Sections 2 to 7 are mainly concerned with proving theorem 7.1, which is
the main tool used in the later sections and gives a complete
description of the Fourier expansion of $\Phi$. In the nonsingular case
this theorem is essentially well known; the main point is to check
that it remains true for functions with singularities at cusps.
Sections 8 and 9 give some minor auxiliary results.

In section 10 we give a complete description of the functions
$\Phi$ in the case $b^+=1$, when the Grassmannian is $b^-$-dimensional
hyperbolic space. In this case the main theorem (10.3) is that
the function $\Phi$ is
a piecewise polynomial function on hyperbolic space.
The Weyl vector $\rho(K,W,F_K)$ in the product formula
below comes from a piecewise linear automorphic form on hyperbolic space.
The function $\Phi$ has properties similar
to those of Donaldson polynomials for 4-manifolds with $b^+=1$;
for example they both have similar ``wall crossing formulas''
depending on the coefficients of modular forms. (See problem 16.7.)

In section 11 we apply the results of section 10 to
find some congruences for the coefficients of theta
functions of all positive definite lattices, generalizing the fact
that the number of roots of a Niemeier lattice is divisible by 24. The
idea of the proof of these congruences is to show that certain ``Weyl
vectors'' related to automorphic forms are in certain lattices, and
then calculate the coefficients of these Weyl vectors explicitly in
terms of coefficients of theta functions of lattices.

In section 12 we give an application of section 10
to hyperbolic reflection groups,
by giving a sufficient condition (in terms of the existence of a
modular form with certain properties) for the reflection group of a
Lorentzian lattice to have finite index in its automorphism group.
This gives many of the known examples of such lattices; for example,
it takes only a couple of lines to recover Vinberg and Kaplinskaja's
result that the reflection group of $I_{1,19}$ has finite index in the
automorphism group.

In section 13 we construct holomorphic automorphic forms
on the hermitian symmetric space $G(2,b^-)$ as infinite products
of the form
$$\e((Z,\rho(K,W,F_K)))
\prod_{\lambda\in K'\atop (\lambda,W)>0}
\prod_{\delta\in M'/M\atop \delta|L=\lambda}
(1-\e((\lambda,Z)+(\delta,z')))^{c_{\delta}(\lambda^2/2)}.
$$
(where $\e(x)=e^{2\pi i x}$, $M$ is a lattice in $\R^{2,b^-}$, and the
$c_\delta$'s are the coefficients of a vector valued modular form of
weight $1-b^-$).  This generalizes the level 1 case of [B95 theorem
10.1] to arbitrary levels and to non unimodular lattices. One new
phenomenon that appears in the higher level case is that a holomorphic
automorphic form can have several apparently quite different infinite
product expansions, one for each orbit of cusps; see example 13.7.
When the automorphic
form has singular weight it usually turns out to be the denominator
formula of a generalized Kac-Moody algebra or superalgebra, which
allows us to construct many new examples of such Lie algebras. In
particular we can recover the examples in [B95] and the higher level
examples worked out by Gritsenko and Nikulin in [G-N].  (The first
suggestion that infinite dimensional Kac-Moody algebras might be
related to automorphic forms seems to be due to Feingold and Frenkel
[F-F]. They suggested that the hyperbolic Kac-Moody algebra associated
to $A_1$ might be related to some Siegel automorphic form of genus
2. To fit into the framework of this paper their Kac-Moody algebras need
to be embedded into larger generalized Kac-Moody algebras, and Siegel
automorphic forms need to be replaced by automorphic forms for
$O_{2,n}$, although of course in the genus 2 case Siegel automorphic
forms are essentially the same as automorphic forms for $O_{2,3}$.)

In section 14 we give a common
generalization of several well known correspondences,
including the Shimura and Maass-Gritsenko correspondences,
to modular forms with poles at cusps.
More precisely we show how to construct automorphic forms
(possibly with poles along rational quadratic divisors)
of weight $m^+>0$ on a Grassmannian of $\R^{2,b^-}$ from
modular forms (possibly with poles at cusps) of weight
$1+m^+-b^-$. For example, if we take $b^-=1$, then we go
from modular forms of weight $m^++1/2$ to automorphic forms
of weight $m^+$ on $O_{2,1}(\R)$, which are essentially the same
as modular forms of weight $2m^+$, and in the special case of holomorphic
modular forms this is the Shimura correspondence.
The case when the modular forms have no poles is essentially
due to Oda [O], Rallis-Schiffmann [R-S], and Gritsenko [G].
See example 14.4 for examples of the Shimura correspondence
for modular forms with singularities.

In section 15 we give a few miscellaneous examples. In particular
find an automorphic form on the moduli space of Ricci flat K3 surfaces
with a $B$-field that is invariant under mirror symmetry,
and show that some Donaldson invariants of 4 manifolds are related
to some of the automorphic forms constructed in section 10.

Finally section 16 lists some
possible topics for further research.

\proclaim Notation and terminology.

\item{$\lambda_V$} is the orthogonal
 projection of a vector $\lambda$ onto a subspace $V$.
\item{$G^+$} If $G$ is a subgroup of a real orthogonal group then
$G^+$ means the elements of $G$ whose spinor norm
has the same sign as the determinant.
\item {$M'$} If $M$ is a lattice then $M'$ means the dual of $M$.
\item {$v^\perp$} Orthogonal complement of a vector (or sublattice)
of a lattice $v$.
\item{$\bar{f}$} Complex conjugate of a function $f$.
\item{$\hat{f}$} The  Fourier transform of $f$.
\item {$\sqrt{}$} We always use the principal value with positive real part,
or zero real part and nonnegative imaginary part.
\item{${A\choose B}$} is 0 if $B$ is a negative integer, 1 if $B=0$,
and $A(A-1)\cdots(A-B+1)$ if $B$ is a positive integer.
\item{$\alpha$} A vector of $M\otimes \R$.
\item{$a$} An entry of a matrix ${ab\choose cd}$ in $SL_2(\Z)$.
\item{$\beta$} A vector of $M\otimes \R$.
\item{$b^\pm$} The lattice $M$ has signature $(b^+,b^-)$.
\item{$b$} An entry of a matrix ${ab\choose cd}$ in $SL_2(\Z)$.
\item{$B_n$} A Bernoulli number.
\item{$\B_n(x)$} A Bernoulli piecewise
polynomial $-n!\sum_{j\ne 0}\e(jx)/(2\pi i j)^n$.
\item{$\gamma$} An element of $M'/M$
\item{ $\Gamma_0(N)$} $ \{{ab\choose cd}\in SL_2(\Z)|c\equiv 0 \bmod N\}$
\item{$\Gamma(z)$} Euler's gamma function.
\item{$c$} An integer, often an entry of a matrix ${ab\choose cd}$ in
$SL_2(\Z)$.
\item{$c_\gamma(m,k)$} A coefficient of $F$ when $F$ is almost holomorphic.
\item{$c_{\gamma,m}(y)$} A coefficient of $F$.
\item {\C} The complex numbers.
\item {$C$} The positive open cone in a Lorentzian lattice.
\item{$\delta$} An element of $M'/M$ or $\Z/N\Z$.
\item{$\delta^m_n$} 1 if $m=n$, 0 otherwise.
\item {$\Delta$} The  delta function,
$\Delta(\tau)=q\prod_{n>0}(1-q^n)^{24}$, or a Laplacian operator.
\item{$d$} An integer, often an entry of a
matrix ${ab\choose cd}$ in $SL_2(\Z)$.
\item{$\e(x)$} $\exp({2\pi i x})$.
\item{$e_\gamma$} An element of a basis of $\C[M'/M]$.
\item{$E_k$} An Eisenstein series of weight $k$,
equal to $1-(2k/B_{k})\sum_{n>0}\sigma_{k-1}(n)q^n$ if $k\ge 2$.
\item{$\E_2$} The non-holomorphic modular form  $E_2(\tau)-3/\pi \Im(\tau)$
of weight 2.
\item {$\zeta$} The Riemann zeta function.
\item {$f$} A function.
\item {$f_\gamma$} A component of $F$.
\item {$F$} A vector valued modular form with components $f_\gamma$.
See theorem 5.3 for $F_M$, $F_K$.
\item{$F_w$} The set of complex numbers $\tau$ with $|\Re(\tau)|\le 1/2$,
$|\tau|\ge 1$, and $0<\Im(\tau)\le w$.
\item{$\gamma$} An element of $M'/M$.
\item {$G(M)$} A Grassmannian, equal to  the set of maximal positive definite
subspaces of some real vector space with a symmetric bilinear form.
\item{$G_N,\G_N$} Variations of Zagier's function; see section 9.
\item {$\eta(\tau)$} $=q^{1/24}\prod_{n>0}(1-q^n)$.
\item {$h^\pm$} Integers; see section 5.
\item{$H$} The upper half plane or a Hurwitz class number.
\item{$\theta$, $\theta_K$} Theta functions of lattices or cosets of lattices.
\item{$\Theta$} A vector valued theta function.
\item{$i$} $\sqrt{-1}$.
\item {$I_{m,n}$} The odd unimodular  lattice of dimension
$m+n$ and signature $m-n$.
\item {$II_{m,n}$} The even unimodular  lattice of dimension
$m+n$ and signature $m-n$.
\item{$\Im$} The imaginary part of a complex number.
\item {$j$} The elliptic modular function $j(\tau)=q^{-1}+744+196884q+\cdots$,
or an integer.
\item {$k$} An integer, often an exponent of $1/y$.
\item {$K$} An even lattice of signature  $(b^+-1,b^--1)$ equal to $L/\Z z$.
\item {$K_\mu$} A modified Bessel function.
\item {$\lambda$} An element of $M$.
\item {$\Lambda$} The Leech lattice. See [C-S].
\item {$\log$} We always use the principal
value with $-\pi<\Im(\log(*))\le \pi$.
\item {$L$} An even  singular lattice equal to $M\cap z^\perp$.
\item{$\mu$} A vector of $K\otimes \R$ defined in section 5.
\item {$M$} An even  lattice of signature  $(b^+,b^-)$.
\item{$m$} An integer.
\item {$m^\pm$} The degree of $p$ is $(m^+,m^-)$.
\item {$Mp_2(\Z)$} The metaplectic group, a double cover of $SL_2(\Z)$.
\item{$n$} An integer, often indexing the coefficients of a modular form.
\item{$N$} The largest integer such that $z/N\in M'$.
\item {$O$} An orthogonal group.
\item{$O(q^n)$} A sum of terms  of order at most $q^n$.
\item{$\pi$} $3.14159\ldots$
\item {$p$} A homogeneous polynomial on $\R^{b^+,b^-}$ of degree $(m^+,m^-)$.
\item {$p_{w,h^+,h^-}$} See section 5.
\item {$P$} A principal $\C^*$ bundle over $H$.
\item {$q$} $e^{2\pi i \tau}$
\item {$\Q$} The rational numbers.
\item {$\rho(M,W,F)$} A Weyl vector (see section 10).
\item {$\rho_M$} A representation of $Mp_2(\Z)$ (see section 4).
\item {$\Re$} The real part of a complex number.
\item {$\R$} The real numbers.
\item{$\sigma_{k-1}(n)$} $=\sum_{d|n}d^{k-1}$ if $n>0$, $-B_k/2k$ if $n=0$.
\item{$S$} The element
$S=({0-1\choose
1\phantom{-}0},\sqrt\tau)$ of $Mp_2(\Z)$.
\item {$SL$} A special linear group.
\item {$\tau$} A complex number $x+iy$ with positive imaginary part $y$.
\item{$T$} The element  $T=({11\choose 01},1)$ of $Mp_2(\Z)$.
\item{$\Phi_M$} An automorphic form with singularities on $G(M\otimes \R)$
defined in section 6.
\item{$\Psi_M$} A meromorphic  automorphic form of weight $k$ on $P$.
See sections 13 and 14.
\item{$\Psi_z$} A restriction of $\Psi_M$ to the hermitian symmetric
space $K\otimes\R+iC$.
\item{$v$} An isometry from $M\otimes\R$ to $\R^{b^+,b^-}$.
\item {$v^\pm$} The inverse image of $\R^{b^\pm}$ under $v$.
\item {$V$} A vector space.
\item {$W$} A Weyl chamber. See section 10.
\item{$x$} A real number, often equal to
$\Re(\tau)$.
\item{$X$} The real part of $Z$,
which is in the Lorentzian space $K\otimes \R$.
\item{$X_M$} The real part of $Z_M$.
\item{$y$} A real number, often equal to $\Im(\tau)$.
\item{$Y$} The imaginary part of $Z$, which is in $C$.
\item{$Y_M$} The imaginary  part of $Z_M$.
\item{$z$} A primitive norm 0 vector of $M$.
\item{$z'$} A vector of $M'$ such that $(z,z')=1$.
\item{$Z$} The element  $Z=({-1\phantom{-}0\choose \phantom{-}0-1},i)$
generating the center of order 4 of $Mp_2(\Z)$, or the element
$X+iY\in M\otimes\C$.
\item{$Z_M$} $=(Z,1,-Z^2/2-z'^2/z)=X_M+iY_M$
\item {$\Z$} The integers.

\proclaim Terminology.

\item{}{\bf Automorphic form.} See section 13.
\item{}{\bf Koecher  principle.} An  automorphic form holomorphic
everywhere except possibly at the cusps on a simple group of
rank greater than 1 is automatically holomorphic at the cusps.
\item{}{\bf Primitive.} A sublattice $K$ of $M$ is primitive if
$M/K$ is torsion free. A vector of $M$ is primitive if it generates
a primitive sublattice.
\item{}{\bf Rational quadratic divisor.} The zero set of $a(y,y)+(b,y)+c$
where $a\in \Z, b\in K,c\in \Z$.
\item{}{\bf Singular weight.} Weight $b^--1$ or 0 (for automorphic
forms on $G(\R^{2,b^-})$.
\item{}{\bf Spinor norm.} A homomorphism from a real orthogonal group to
$\R^*/\R^{*2}$
taking reflections of vectors of positive or negative norm to $1$
or $-1 $ respectively.
\item{}{\bf Theta function.} A modular form or Jacobi form
depending on a lattice.
\item{}{\bf Weyl chamber.} A generalization of the Weyl chamber of
a root system. See section 6.
\item{}{\bf Weyl vector.} A vector $\rho(M,W,F)$ such that the inner
product with $\rho(M,W,F)$ is a multiple of $\Phi$; see section 10.

\proclaim 2.~Modular forms.

In this section we summarize some slightly nonstandard facts about
modular forms that we will use later. The main differences to the
common treatments of modular forms are that we replace the concept of
a modular form of high level by the more precise and more general
concept of a vector valued modular form associated to some projective
representation $\rho$ of $SL_2(\Z)$, and we also allow modular forms
to be non holomorphic and to have singularities at the cusps.

Recall that the group $SL_2(\Z)$ has a double cover $Mp_2(\Z)$
called the metaplectic group whose
elements can be written in the form
$$\left({ab\choose cd},\pm\sqrt{c\tau +d}\right)
$$
where $\left({ab\choose cd})\right)\in SL_2(\Z)$
and $\sqrt{c\tau+d}$ is considered as a holomorphic function
of $\tau$ in the upper half plane whose square is $c\tau+d$.
The multiplication is defined so that the usual formulas for
the transformation of modular forms work for half integer weights,
which means that
$$
(A,f(\cdot))(B,g(\cdot))= (AB,f(B(\cdot))g(\cdot))
$$
for $A,B\in SL_2(\Z)$ and $f,g$ suitable functions on $H$.

Suppose that $\rho$ is a representation of $Mp_2(\Z)$ on a vector
space $V$, and suppose that $m^+$, $m^-$ are integers
or half integers.  We define a modular form of weight $(m^+,m^-)$ and type
$\rho$ to be a real analytic function $F$ on the upper half plane $H$
with values in $V$ such that
$$\eqalign{
&F((a\tau+b)/(c\tau+d))\cr
=& (c\tau+d)^{m^+}(c\bar\tau+d)^{m^-}
\rho\left({ab\choose cd},\sqrt{c\tau+d} \right)
F(\tau)\cr
}
$$
for  elements $\left({ab\choose cd},\sqrt{c\tau +d}\right)$
of the metaplectic group. Note that the factor
$(c\tau+d)^{m^+}$ means $\sqrt{c\tau+d}^{2m^+}$ when
$2m^+$ is odd, and similarly for $(c\bar\tau+d)^{m^-}$.
 If $\tau\in H$ we write
$x$ and $y$ for the real and imaginary parts of $\tau$.  We will say
that $F$ is almost holomorphic of weight $(m^+,m^-)$ if all components
of $F$ can be written as
$$\sum_{m\in \Q} \sum_{k\in \Z}c(m,k)\e(m\tau)y^{-k}
$$
(where $\e(x)$ means $e^{2\pi i x}$)
and if the coefficients $c(m,k)$ vanish whenever $m<<0$ or $k<0$
or $k>>0$ (note that we allow $F$ to have ``poles''
of finite order at cusps). We say that $F$ is holomorphic on $H$ if
it has weight $(m^+,0)$ for some $m^+$ and
the coefficients $c(m,k)$ vanish whenever  $k\ne 0$, and we say that
$F$ is holomorphic if in addition the coefficients vanish
whenever $m<0$.

{\bf Example 2.1.} The function $F(\tau)=y$ is an almost holomorphic
modular form of weight $(-1,-1)$. In particular any modular
form of weight $(m^+,m^-)$ can be turned into one of
weight $(m^+-m^-,0)$ in a canonical way by multiplying it by $y^{m^-}$.
(So we would lose little generality by only considering forms
of weights $(m^+,0)$, but this seems a little unnatural;
for example, Siegel theta functions of lattices of signature
$(b^+,b^-)$ have weights $(b^+/2,b^-/2)$.)

{\bf Example 2.2.} If $f$ is a (classical) holomorphic modular form of
level $N$ corresponding to some character $\chi$ of some subgroup
$\Gamma$ of finite index in $SL_2(\Z)$, then $f$ induces a holomorphic
modular form $F$ of type $V$ where $V$ is the induced representation
$Ind_\Gamma^{SL_2(\Z)}(\chi)$. The components of $F$ are (more or
less) the Fourier expansions of $f$ at the cusps of $\Gamma$. In
particular we do not lose any generality by only considering ``level
1'' vector valued modular forms. The induced representation is often
reducible, so we can specify level $n$ forms more precisely by
specifying some sub representation of the induced representation that
their image has to lie in; see lemma 2.6 below.

{\bf Example 2.3.} Jacobi forms as in [E-Z] can all be considered as
vector valued modular forms. More precisely theorem 5.1 of [E-Z]
implies that their space $J_{k,m}$ of Jacobi forms of weight $k$ and
index $m$ is naturally isomorphic to the space of holomorphic modular
forms of weight $k-1/2$ and representation $\bar\rho_M$ dual to $\rho_M$,
where $M$ is a 1-dimensional lattice generated by a vector of norm
$2m$.

{\bf Example 2.4.} The Kohnen ``plus space'' [Ko] has a natural interpretation
in terms of vector valued modular forms as in the proof of [E-Z, theorem 5.4].
In particular modular forms of level 4 and half integer weight
satisfying the plus space condition are essentially the same
as certain level 1 vector valued modular forms.

{\bf Example 2.5.} The real analytic function
$\E_2(\tau)=E_2(\tau)-3/\pi \Im(\tau)$ is an almost holomorphic
modular form of weight $(2,0)$.

\proclaim Lemma 2.6. Suppose that $f$ is a complex valued modular form of
weight $(m^+,m^-)$ for the group $\Gamma_1(N)=\{{ab\choose cd}\in SL_2(\Z)
|a\equiv d\equiv 1\bmod N, c\equiv 0 \bmod N\}$, and write
$f({*\tau+*\over c\tau+d}) $ for the number $f({a\tau+b\over c\tau+d})$
for any ${ab\choose cd}\in SL_2(\Z)$ (which is well defined
as $f(\tau+1)=f(\tau)$).
Let $M$ be the 2 dimensional lattice generated by norm 0 vectors
$z$, $z'$ with $(z,z')=N$. Define a vector valued function $F(\tau)$
by
$$
f_{cz'/N+\gamma z/N}(\tau)=
\sum_{d\in \Z/N\Z\atop (d,c,N)=1}
\e(-\gamma d/N)
f({*\tau+*\over c\tau+d}).
$$
Then $F$ is a modular form of type $\rho_M$ (defined in section 4).

Proof. We have to check that $F$ transforms correctly under the
elements $S$ and $T$.

For $T$ we see that
$$\eqalign{
&f_{cz'/N+\gamma z/N}(\tau+1)\cr
=&\sum_{d\in \Z/N\Z\atop (d,c,N)=1}
\e(-\gamma d/N)f({*\tau+*\over c\tau+c+d})\cr
=&\sum_{d\in \Z/N\Z\atop (d,c,N)=1}
\e(-\gamma d/N)\e(\gamma c/N)f({*\tau+*\over c\tau+d})\cr
=& \e((cz'/N+\gamma z/N)^2/2)f_{cz'/N+\gamma z/N}(\tau).\cr
}
$$

For the generator $S$ we see that
$$\eqalign{
&f_{cz'/N+\gamma z/N}(-1/\tau)\cr
=&\tau^{m^+}\bar\tau^{m^-}
\sum_{d\in \Z/N\Z\atop (d,c,N)=1} \e(-\gamma d/N)f({*\tau+*\over -c+d\tau})\cr
=&{\tau^{m^+}\bar\tau^{m^-}\over N}\sum_{d,\delta\in \Z/N\Z}
\e(-c\delta N-\gamma d/N)
\sum_{\epsilon\in \Z/N\Z\atop (\epsilon,d,N)=1}
\e(-(\epsilon,\delta))f({*\tau+*\over d\tau+\epsilon})\cr
=&{\tau^{m^+}\bar\tau^{m^-}\over \sqrt{|M'/M|}}\sum_{d,\delta\in \Z/N\Z}
\e(-(cz'/N+\gamma z/N,zd/N+\delta z/N))
f_{z'd/N+\delta z/N}(\tau).\cr
=&{\tau^{m^+}\bar\tau^{m^-}\over \sqrt{|M'/M|}}\sum_{\delta\in M'/M}
\e(-(cz'/N+\gamma z/N,\delta ))
f_{\delta }(\tau).\cr
}
$$
This proves lemma 2.6.

\proclaim 3.~Fourier transforms.

In this section we evaluate some well known Fourier transforms that we
will need later.

We define $\e(x)$ to be $\exp(2\pi i x)$. If $V$ is a real vector
space with a positive definite quadratic form given by
$(x,x)=\sum_jx_j^2$ in some orthonormal basis, then the Laplacian
operator $\Delta $ is defined to be
$$\Delta = \sum_j{d^2\over dx_j^2}.
$$
On $\R^{b^+,b^-}$ we define $\Delta$ to be the Laplacian of
$\R^{b^++b^-}$. (Note that this is not the Laplacian of
$\R^{b^+,b^-}$ and is not invariant under rotations of $\R^{b^+,b^-}$.)

We recall some standard properties of the Fourier transform $\hat f$
of a function $f$ on a vector space $V$ with a nonsingular symmetric
bilinear form of signature $(b^+,b^-)$, defined as $\hat
f(y)=\int_{x\in V}f(x)\e((x,y))dx$.
\proclaim Lemma 3.1.
\item{1.} The Fourier transform of $f(x-a)$ is $\e(ax)\hat f(x)$.
\item{2.} The Fourier transform of $f(x)\e(ax)$ is $\hat f(x+a)$.
\item{3.} The Fourier transform of $xf(x)$ is ${d\over dx}\hat f (x)/2\pi i$.
\item{4.} The Fourier transform of ${d\over dx}f(x)$ is $-2\pi ix\hat f (x)$.
\item{5.} If $a>0$ then the Fourier transform
of $f(ax)$ is $a^{-b^+-b^-}\hat f(x/a)$.
\item{6.} If $b^-=0$ then the
Fourier transform of $e^{-\pi x^2}$ is $e^{-\pi x^2}$.

The proofs of these are all standard (and easy) and will be omitted.

\proclaim Lemma 3.2. Suppose $p$ is a  polynomial  on $b^+$-dimensional
Euclidean space
and $\Im(\tau)>0$.
Write $\Delta$ for the Laplacian operator.
Then the Fourier transform of
$p(x)\e(x^2\tau/2) $ is
$$
(\tau/i)^{-b^+/2}\exp(i\Delta/4\pi\tau)(p)(-x/\tau)\e(- x^2/2\tau).
$$

Proof. This result obviously follows from the 1 dimensional case.
We prove it for $p(x)=x^m$ by induction on $m$, in which case
it is equivalent to showing that the Fourier transform of
$x^m\e(x^2\tau/2)$ is $(\tau/i)^{-1/2}(-\tau)^{-m}
\exp(i\tau\Delta/4\pi)(p)(x)\e(-x^2/2\tau)$.
A short calculation shows that
$$\exp(i\tau\Delta/4\pi)(xp)=
x\exp(i\tau\Delta/4\pi)(p)+i\tau\exp(i\tau\Delta/4\pi)(p')/2\pi
$$
for any polynomial $p$, and in particular for $p(x)=x^m$.
Using this and lemma 3.1 and induction on $m$ we see that
the Fourier transform of
$x\times x^me^{-\pi x^2}$ is
$$\eqalign{
&(\tau/i)^{-1/2}{1\over 2\pi i}{d\over dx}
\left((-\tau)^{-m}\exp(i\tau\Delta/4\pi)(x^m)\e(-x^2/2\tau)\right)\cr
=& (\tau/i)^{-1/2}{1\over 2\pi i}(-\tau)^{-m}
(\exp(i\tau\Delta/4\pi)(mx^{m-1})\e(-x^2/2\tau)+\cr
&\qquad\qquad\qquad\qquad\qquad
-2\pi ix\exp(i\tau\Delta/4\pi)(x^m)\e(-x^2/2\tau)/\tau)\cr
=&(\tau/i)^{-1/2} {1\over 2\pi i}(-2\pi i/\tau)(-\tau)^{-m}
\exp(i\tau\Delta/4\pi)(x^{1+m})\e(-x^2/2\tau)\cr
=&(\tau/i)^{-1/2} (-\tau)^{-m-1}
\exp(i\tau\Delta/4\pi)(x^{1+m})\e(-x^2/2\tau)\cr
}
$$
which proves lemma 3.2 by induction on $m$.

\proclaim Corollary 3.3. The Fourier transform of
$p(x)\e(Ax^2+Bx+C)$ is
$$(2A/i)^{-1/2}\exp(i\Delta/8\pi A)(p)((-x-B)/2A)
\e(-x^2/4A-xB/2A+C-B^2/4A)$$ for
$A,B,C$ complex, $\Im(A)>0$, $x \in V=\R$, $(x,y)=xy$.

Proof. This follows from lemma 3.2 (with $\tau=2A$) by applying
lemma 3.1 part 2 (with $a=B$).

\proclaim Corollary 3.4. If $p$ is a polynomial
on $b^+$-dimensional Euclidean space and $\Im(\tau)>0$ then
the Fourier transform of
$$\exp(-\Delta/8\pi\Im(\tau))(p)(x)\e(\tau x^2/2)$$
is
$$(\tau/i)^{-b^+/2}\exp(-\Delta/8\pi \tau^2\Im(-1/\tau))(p)(-x/\tau)
\e(-x^2/2\tau)
$$
which is equal to
$$(\tau/i)^{-b^+/2} (-\tau)^{-m}\exp(-\Delta/8\pi\Im(-1/\tau))(p)(x)
\e(-x^2/2\tau)$$
if $p$ is homogeneous of degree $m$.

Proof. Applying lemma 3.2 shows that the
Fourier transform of
$$\exp(-\Delta/8\pi \Im(\tau))(p)(x)\e(x^2\tau/2)
$$
is $$\eqalign{
&(\tau/i)^{-b^+/2}\exp(-\Delta/8\pi \Im(\tau)+i\Delta/4\pi\tau)(p)(-x/\tau)
\e(-x^2/2\tau)\cr
=&(\tau/i)^{-b^+/2}\exp(-\Delta\bar\tau/8\pi \tau\Im(\tau))(p)(-x/\tau)
\e(-x^2/2\tau)\cr
=&(\tau/i)^{-b^+/2}\exp(-\Delta/8\pi \tau^2\Im(-1/\tau))(p)(-x/\tau)
\e(-x^2/2\tau).\cr
}
$$
This proves corollary 3.4.

\proclaim Corollary 3.5. Suppose that $p$ is a homogeneous
polynomial of degree $(m^+,m^-)$ on the sum of positive and negative
definite spaces $v^+$ and $v^-$ of dimensions $b^+$ and $b^-$ (which
means that $p$ has degree $m^+$ in the variables of $v^+$ and degree
$m^-$ in the variables of $v^-$).  Then the Fourier transform of
$$\exp(-\Delta/8\pi\Im(\tau))(p)(x)\e(\tau x_{v^+}^2/2+\bar\tau x_{v^-}^2/2)$$
is
$$(\tau/i)^{-b^+/2} (-\tau)^{-m^+}(i\bar\tau)^{-b^-/2} (-\bar\tau)^{-m^-}
\exp(-\Delta/8\pi\Im(-1/\tau))(p)(x)\e(-x_{v^+}^2/2\tau-x_{v^-}^2/2\bar\tau).
$$

Proof. We can assume that $p$ is the product of homogeneous polynomials
of degrees $m^+$ and $m^-$ on $v^+$ and $v^-$.
Corollary 3.5 follows by applying corollary 3.4 to $v^+$ and $v^-$.

\proclaim 4.~Siegel theta functions.

In this section we summarize some standard results about Siegel theta
functions of indefinite lattices. In most cases the proofs are easy
generalizations of the proofs for positive definite lattices, which
can be found in any standard reference about theta functions,
see for example Shintani [S]. We will make some minor modifications
to the usual treatment of theta functions to make later applications
easier; for example, we use vector valued forms of level 1 rather than
forms of higher level.

We let $M$ be an even lattice of signature $(b^+,b^-)$, with dual
$M'$.  Recall that the mod 1 reduction of $(\lambda,\lambda)/2$ is a
$\Q/\Z$-valued quadratic form on $M'/M$, whose associated
$\Q/\Z$-valued bilinear form is the mod 1 reduction of the bilinear
form on $M'$.  We use $v$ to denote an isometry from $M\otimes \R$ to
$\R^{b^+,b^-}$.  We write $v^+$ and $v^-=v^{+\perp}$ for the inverse
images of $\R^{b^+,0}$, $\R^{0,b^-}$ under $v$, so that $M\otimes \R$
is the orthogonal direct sum of the positive definite subspace $v^+$
and the negative definite subspace $v^-$. The Grassmannian $G(M) $ is
the set of positive definite $m^+$-dimensional subspaces $v^+$ of
$M\otimes \R$ and the projection of $\lambda\in M\otimes \R$ into a
subspace $v^\pm$ is denoted by $\lambda_{v^\pm}$, so that
$\lambda=\lambda_{v^+}+\lambda_{v^-}$. The Siegel theta function
$\theta_M$ of $M$ is defined by
$$\theta_M(\tau;v^+)= \sum_{\lambda\in M}\e(\tau\lambda_{v^+}^2/2
+\bar\tau\lambda_{v^-}^2/2)
$$
for $\tau\in H$, $v^+\in G(M)$. It will be useful later to have a more
general theta function defined by
$$\eqalign{
&\theta_{M+\gamma}(\tau,\alpha,\beta;v,p)\cr
=& \sum_{\lambda\in M+\gamma}
\exp(-\Delta/8\pi y)(p)(v(\lambda+\beta))\e(\tau(\lambda+\beta)_{v^+}^2/2
+\bar\tau(\lambda+\beta)_{v^-}^2/2 - (\lambda+\beta/2,\alpha))\cr
}
$$
for $\alpha, \beta\in M\otimes \R$, $\gamma\in M'/M$, $v$ an isometry from
$\R^{b^+,b^-}$ to $M\otimes \R$,  and $p$
a polynomial on $\R^{b^+,b^-}$, homogeneous
of degree  $m^+$ in the first $b^+$ variables, and of degree
$m^-$ on the last $b^-$ variables.
We will sometimes omit some of the arguments: if $\alpha$ and $\beta$
are both 0 we miss them out, if $p=1$ we miss it out, and if
$G(b^+,b^-)$ is a point then we miss out $v$.

It is common to restrict $p$ to be a harmonic homogeneous polynomial,
but there seems to be no good reason for this restriction
and we will not make it.
There are also good reasons
for not restricting $p$ to be harmonic. For example, later on we need to
write $p$ as a linear combination of products of polynomials
on subspaces, and there is no reason for the polynomials on
subspaces to be homogeneous and harmonic even if $p$ is homogeneous
and harmonic.

We let the elements $e_\gamma$ for $\gamma\in M'/M$ be the standard
basis of the group ring $\C[M'/M]$, so that $e_\gamma
e_\delta=e_{\gamma+\delta}$.  Recall that there is a unitary
representation $\rho_M$ of the double cover $Mp_2(\Z)$ of $SL_2(\Z)$
on $\C[M'/M]$ defined by
$$\rho_M(T)(e_\gamma) = \e((\gamma,\gamma)/2)e_\gamma
$$
$$\rho_M(S)(e_\gamma) = {\sqrt{i}^{b^--b^+}\over \sqrt{|M'/M|}}
\sum_{\delta\in M'/M} \e(-(\gamma,\delta))e_\delta
$$
where $T=({11\choose 01},1)$ and $S=({0-1\choose
1\phantom{-}0},\sqrt\tau)$ are the standard generators of $Mp_2(\Z)$,
with $S^2=(ST)^3=Z$, $Z=({-1\phantom{-}0\choose \phantom{-}0-1},i)$,
$Z(e_\gamma)=i^{b^--b^+}e_{-\gamma}$, $Z^4=1$.  The representation
$\rho_M$ is essentially the Weil representation of the self dual
abelian group $M'/M$ with the quadratic character $\gamma^2/2$, and
there is an explicit formula for it in some cases in [W]
and in all cases in [S]; for example
formula 16 of [W] describes it when $c$ is coprime to $|M'/M|$. We
will not use these explicit formulas.  The representation $\rho_M$
factors through the finite group $SL_2(\Z/N\Z)$ if $M$ has even
dimension, and through a double cover of $SL_2(\Z/N\Z)$ if $M$ has odd
dimension, where $N$ is the smallest integer such that
$N(\gamma,\delta)$ and $N\gamma^2/2$ are integers for all
$\gamma,\delta\in M'$.

We will write
$\Theta_M$ for the $\C[M'/M]$-valued function
$$\Theta_M(\tau,\alpha,\beta;v,p)=\sum_{\gamma\in
M'/M}e_\gamma\theta_{M+\gamma}(\tau,\alpha,\beta;v,p)
{}.
$$
Then the transformation formula for
this vector valued function is

\proclaim Theorem 4.1. If $p$ is a homogeneous polynomial
of degree $(m^+,m^-)$
and $({ab\choose cd},\sqrt{c\tau+d})\in Mp_2(\Z)$
then
$$\eqalign{
&\Theta_M((a\tau+b)/(c\tau+d),a\alpha+b\beta,c\alpha+d\beta;v,p)\cr
=& (c\tau+d)^{b^+/2+m^+}(c\bar\tau+d)^{b^-/2+m^-}
\rho_M\left({ab\choose cd},\sqrt{c\tau+d} \right)
\Theta_M(\tau,\alpha,\beta;v,p).\cr
}
$$

Proof. It is easy to check that if this is true for two elements of
$Mp_2(\Z)$ it is true for their product, so it is sufficient to check
it for the standard generators $T=({11\choose 01},1)$ and
$S=({0-1\choose 1\phantom{-}0},i)$. For the first generator $T$ we
have to show that
$$\theta_{M+\gamma}(\tau+1,\alpha+\beta,\beta;v,p)
= \e((\gamma,\gamma)/2)\theta_{M+\gamma}(\tau,\alpha,\beta;v,p)$$
which is trivial to check.
For the second generator $S$ we
have to show
$$\eqalign{
&\sqrt{|M'/M|}\theta_{M+\gamma}(-1/\tau,-\beta,\alpha;v,p)\cr
=&
\sqrt{\tau/i}^{b^+}\tau^{m^+}\sqrt{i\bar\tau}^{b^-}\bar\tau^{m^-}
\sum_{\delta\in M'/M}\e(-(\gamma,\delta))
\theta_{M+\delta}(\tau,\alpha,\beta;v,p).\cr
}
$$
By corollary 3.5 the  Fourier transform of
$$
\sqrt{\tau/i}^{-b^+}(-\tau)^{-m^+}{\sqrt{\bar\tau/i}}^{-b^-}(-\bar\tau)^{-m^-}
\exp({-\Delta\over8\pi\Im(-1/\tau)})(p)(v(x))
\e(-x_{v^+}^2/2\tau-x_{v^-}^2/2\bar\tau )
$$
is
$$
(-1)^{m^++m^-}
 \exp(-\Delta/8\pi \Im(\tau))(p)(v(x))\e(\tau x_{v^+}^2/2
+\bar\tau x_{v^-}^2/2).
$$

If we change $x$ to $x+\alpha+\gamma$ and then multiply by
$\e((x+\gamma+\alpha/2,\beta))$ and use lemma 3.1 we find that the
function $f$ given by
$$\eqalign{
&\sqrt{\tau/i}^{-b^+}\tau^{-m^+}{\sqrt{\bar\tau/i}}^{-b^-}
\bar\tau^{-m^-}\exp(-\Delta/8\pi\Im(-1/\tau))(p)(v(x+\alpha+\gamma))
\times\cr
\times&\e((-1/\tau)(x+\alpha+\gamma)_{v^+}^2/2
+(-1/\bar\tau)(x+\alpha+\gamma)_{v^-}^2/2 - (x+\gamma+\alpha/2,-\beta))\cr
}
$$
has Fourier transform $\hat f(x)$ equal to
$$\exp(-\Delta/8\pi \Im(\tau))(p)(v(x+\beta))\e(\tau(x+\beta)_{v^+}^2/2
+\bar\tau(x+\beta)_{v^-}^2/2 - (x+\beta/2,\alpha)-(x,\gamma)).
$$

We apply the Poisson summation formula
$\sqrt{|M'/M|}\sum_M   f=\sum_{M'}{\hat f}$ to the function
above and find
$$\eqalign{
&
\sqrt{\tau/i}^{-b^+}\tau^{-m^+}\sqrt{\bar\tau/i}^{-b^-}\bar\tau^{-m^-}
\sqrt{|M'/M|}\theta_{M+\gamma}(-1/\tau,-\beta,\alpha;v,p)\cr
=&
\sqrt{|M'/M|}\sum_{\lambda\in M}
 f(\lambda)\cr
=&
\sum_{\delta\in M'/M}\sum_{\lambda\in M}\hat f(\lambda+\delta)\cr
=&\sum_{\delta\in M'/M}
\e(-(\delta,\gamma))
\theta_{M+\delta}(\tau,\alpha,\beta;v,p).\cr
}
$$

This verifies the transformation formula for $\Theta$ under $S$ and
completes the proof of theorem 4.1.

The relation $(ST)^3=Z$ gives the following well known generalization
of the law of quadratic reciprocity.

\proclaim Corollary 4.2. (Milgram)
$$\sum_{\gamma\in M'/M} \e(\gamma^2/2) = \sqrt{|M'/M|}\e((b^+-b^-)/8).
$$

Proof. An explicit calculation shows that
$$\eqalign{
&\left(\sqrt{|M'/M|}\sqrt i ^{b^+-b^-}ST\right)^3(e_\gamma)\cr
 &=
\sum_{\delta,\epsilon,\zeta\in M'/M}
\e(-(\gamma,\delta))\e(\delta^2/2)
\e(-(\delta,\epsilon))\e(\epsilon^2/2)
\e(-(\epsilon,\zeta))\e(\zeta^2/2)e_\zeta\cr
&=
\sum_{\delta,\epsilon,\zeta\in M'/M}\e((\epsilon-\delta-\zeta)^2/2)
\e(-(\delta,\zeta+\gamma))e_\zeta\cr
&=
\sum_{\epsilon\in M'/M}\e(\epsilon^2/2)|M'/M|e_{-\gamma}.\cr
}
$$
Comparing this with $(ST)^3(e_\gamma)=Z(e_\gamma)=i^{b^--b^+}e_{-\gamma}$
proves corollary 4.2.

\proclaim 5.~Reduction to smaller lattices.

In this section and section 7 we will work out the Fourier expansion
of the function $\Phi$. The calculations look rather complicated but
the essential idea is easy (and well known) and is as follows. Suppose
that $\Theta(\tau)$ and $F(\tau)$ are modular forms of level 1 and
weights $k$ and $-k$ and we wish to work out the integral
$$\int_{SL_2(\Z)\\H}\Theta(\tau)F(\tau)dxdy/y^2
$$
(ignoring convergence problems for the moment). If we can find an expression
for $\Theta(\tau)$ of the form
$$\Theta(\tau)=\sum_{(c,d)=1}(c\tau+d)^kg({a\tau+b\over  c\tau+d})
$$
then the integral is formally equal to
$$\int_{y>0}\int_{x\in \R/\Z}g(\tau)F(\tau)dxdy/y^2
$$
which is much easier to evaluate as we are integrating over a rectangle.
(The same idea appears in several places in the theory of modular forms
and is known as the Rankin-Selberg method;
for example, we could take
$\Theta$ to be a real analytic Eisenstein series
and take $g(\tau)$ to be a power of $\Im(\tau)$, to see that
the Peterson inner product of $F$ with a real analytic Eisenstein series
is essentially the Mellin transform of the constant term of $F$.)
The rest of this section is mainly concerned with finding
an expression (theorem 5.2) analogous to the one above for $\Theta$ the
Siegel theta function of $M$, when $g$ turns out to be related to
the theta function of a smaller lattice $K$.
The rest of this section consists mainly of computations,
and the reader may skip everything except the statements of
theorems 5.2 and 5.3 and the definition of $F_K$ without any great loss.

We will do this when $\Theta$ is a Siegel theta function by
taking a partial Fourier transform in one variable. (In terms of the
Weil representation this Fourier transform is essentially given
by an element of the Weyl group of $Sp_4$ exchanging two copies of
$SL_2$ corresponding to positive roots.)

Suppose that $z$ is a primitive norm 0 vector of $M$.
In this section we will find a
certain expression for the theta function of $M$ in terms of theta
functions of the lattice $K=(M\cap z^\perp)/\Z z$.

Recall that $z_{v^\pm}$ is the projection of $z$ onto $v^{\pm}$.
We let $w^+$ be the orthogonal complement of $z_{v^+}$ in $v^+$, and we let
$w^-$ be the orthogonal complement of $z_{v^-}$ in $v^-$.
We define the linear map $w$ from $M\otimes \R$ to $\R^{b^+,b^-}$ by
$w(\lambda)=v(\lambda_{w^+}+\lambda_{w^-})$, so that
$w$ is an isomorphism from $w^+$ and $w^-$ to their images,
and $w$ vanishes on $z_{v^+}$ and $z_{v^-}$.

Given a homogeneous polynomial $p$ of degree $(m^+,m^-)$ we define
homogeneous
polynomials $p_{w,h^+,h^-}$
of degrees $(m^+-h^+,m^--h^-)$
on the vector spaces $w(M\otimes \R)$ by
$$p(v(\lambda))=\sum_{h^+,h^-}(\lambda,z_{v^+})^{h^+}
(\lambda,z_{v^-})^{h^-}p_{w,h^+,h^-}(w(\lambda)).
$$

\proclaim Lemma 5.1.
$$\eqalign{
&\theta_{M+\gamma}(\tau;v,p) = \cr
=&{1\over \sqrt{2yz_{v^+}^2}}
\sum_{\lambda\in M/z+\gamma}\sum_{n\in \Z}\sum_{h^+,h^-}
\exp(-\Delta/8\pi y)(p_{w,h^+,h^-})(w(\lambda))\times\cr
&\times\sum_h {h!(-yz_{v^+}^2/\pi)^h\over (-2iy)^{h^++h^-}}
{h^+\choose h}{h^-\choose h}
((\lambda,z)\bar\tau+n)^{h^+-h}((\lambda,z)\tau+n)^{h^--h}\times\cr
&\times
\e\left(\tau\lambda_{w^+}^2/2+\bar\tau\lambda_{w^-}^2/2
-n(\lambda,(z_{v^+}-z_{v^-})/2z_{v^+}^2)-{|(\lambda,z)\tau+n|^2\over
4iyz_{v^+}^2}\right).\cr
}
$$

Proof.
Consider the function
$$\eqalign{
g(\lambda,n)=&\exp(-\Delta/8\pi y)(p)(v(\lambda+nz))
\e(\tau(\lambda+nz)_{v^+}^2/2+\bar\tau(\lambda+nz)_{v^-}^2/2)\cr
=&\exp(-\Delta/8\pi y)(p)(v(\lambda+nz))\times\cr
&\times\e\left((\tau-\bar\tau)z_{v^+}^2n^2/2
+ (\tau(\lambda,z_{v^+})+\bar\tau(\lambda,z_{v^-}))n
+\tau\lambda_{v^+}^2/2+\bar\tau\lambda_{v^-}^2/2\right).
}
$$
The theta function $\theta_{M+\gamma}(\tau;v,p) $ is equal to
$$\sum_{\lambda\in \gamma+M/z}\left(\sum_{n\in \Z}g(\lambda,n)\right)
 =\sum_{\lambda\in \gamma+M/z}\left(\sum_{n\in \Z}\hat g(\lambda,n)\right)
$$
by the Poisson summation formula,
where $\hat g$ is the Fourier transform with respect to the variable
$n$.

We prove lemma 5.1 by working out the Fourier transform $\hat g$
explicitly and substituting it in.
We can work out the Fourier transform $\hat g$ using corollary 3.3
with $A=(\tau-\bar\tau)z_{v^+}^2/2=iyz_{v^+}^2$,
$B=\tau(\lambda,z_{v^+})+\bar\tau(\lambda,z_{v^-})$,
and $C=\tau\lambda_{v^+}^2/2+\bar\tau\lambda_{v^-}^2/2$.
Using the fact that
$$\eqalign{
&\exp(-\Delta/8\pi y)(p)(v(\lambda+nz))\cr
=&
\exp(-{d^2\over dn^2}/8\pi y z_{v^+}^2)
\left((\lambda+nz,z_{v^+})^{h^+}\right)\times\cr
&\times\exp(-{d^2\over dn^2}/8\pi y z_{v^+}^2)
\left((\lambda+nz,z_{v^-})^{h^-}\right)\times\cr
&\times\exp(-\Delta/8\pi y)(p_{w,h^+,h^-})(w(\lambda))
\cr
}
$$
we find $\hat g(n)$ is equal to
$$\eqalign{
&{1\over \sqrt{2yz_{v^+}^2}}\sum_{h^+,h^-}
\exp(-\Delta/8\pi y)(p_{w,h^+,h^-})(w(\lambda))\times\cr
&\times\exp({1\over 8\pi yz_{v^+}^2}{d^2\over dn_2^2})\cr
&\left(
\exp({-1\over8\pi y z_{v^+}^2}{d^2\over dn_2^2}) (\lambda+n_2z,z_{v^+})^{h^+}
\exp({-1\over8\pi y z_{v^+}^2}{d^2\over dn_2^2}) (\lambda+n_2z,z_{v^-})^{h^-}
\right)
\times\cr
&\times\e\left(
{-n_1^2/2
\over
(\tau-\bar\tau)z_{v^+}^2}
+{\tau\lambda_{v^+}^2\over2}+{\bar\tau\lambda_{v^-}^2\over 2}
\right)\cr
}$$
where $n_1=n+\tau(\lambda,z_{v^+})+\bar\tau(\lambda,z_{v^-})$
and $n_2=-n_1/2iyz_{v^+}^2$.
We want to show this is equal to the expression appearing in the lemma.

We evaluate the last factor of this
using the equalities
$\lambda_{v^+}^2=\lambda_{w^+}^2+(\lambda,z_{v^+})^2/z_{v^+}^2$,
$\lambda_{v^-}^2=
\lambda_{w^-}^2+(\lambda,z_{v^-})^2/z_{v^-}^2$,
and $z_{v^+}^2+z_{v^-}^2=0$ to see that
$$\eqalign{
&\e\left(
{-n_1^2/2
\over
(\tau-\bar\tau)z_{v^+}^2}
+{\tau\lambda_{v^+}^2\over2}+{\bar\tau\lambda_{v^-}^2\over 2}
\right)\cr
=&\e\left(
{-n^2/2 -n(\tau(\lambda,z_{v^+})+\bar\tau(\lambda,z_{v^-}))
-(\tau(\lambda,z_{v^+})+\bar\tau(\lambda,z_{v^-}))^2/2
\over
(\tau-\bar\tau)z_{v^+}^2}
+{\tau\lambda_{v^+}^2\over2}+{\bar\tau\lambda_{v^-}^2\over 2}
\right)\cr
&=
\e\left(\tau\lambda_{w^+}^2/2+\bar\tau\lambda_{w^-}^2/2
-n(\lambda,(z_{v^+}-z_{v^-})/2z_{v^+}^2)-{|(\lambda,z)\tau+n|^2\over
2(\tau-\bar\tau)z_{v^+}^2}\right).\cr
}
$$

Next we note that
$$\eqalign{
&\exp(A({d\over dn_3}+{d\over dn_4})^2)
\exp(-A{d^2\over dn_3^2})\exp(-A{d^2\over dn_4^2})\cr
&= \exp(2A{d\over dn_3}{d\over dn_4})\cr
&=\sum_h{(2A)^h\over h!}{d^h\over dn_3^h}{d^h\over dn_4^h}.\cr
}
$$

If we apply this with $n_3=(\lambda+n_2z,z_{v^+})$,
$n_4=(\lambda+n_2z,z_{v^-})$, $A=z_{v^+}^2/8\pi y$, we find that
$$\eqalign{
&\exp({1\over 8\pi yz_{v^+}^2}{d^2\over dn_2^2}
)
\Bigg(
\exp({-1\over8\pi yz_{v^+}^2 }{d^2\over dn_2^2})
(\lambda+n_2z,z_{v^+})^{h^+}\cr
&\qquad\qquad\qquad\qquad
\exp({-1\over8\pi yz_{v^+}^2 }{d^2\over dn_2^2}) (\lambda+n_2z,z_{v^-})^{h^-}
\Bigg)\cr
&=
\sum_h {h!(z_{v^+}^2)^h\over (4\pi y)^h}{h^+\choose h}{h^-\choose h}
(\lambda+n_2z,z_{v^+})^{h^+-h}(\lambda+n_2z,z_{v^-})^{h^--h}.\cr
}
$$
Substituting in
$n_2=-(n+\tau(\lambda,z_{v^+})+\bar\tau(\lambda,z_{v^-}))/2iyz_{v^+}^2$
this becomes
$$\eqalign{ &\sum_h {(-1)^hh!y^{h-h^+-h^-}(z_{v^+}^2/\pi)^h\over
(-2i)^{h^++h^-}}{h^+\choose h}{h^-\choose h}
((\lambda,z)\bar\tau+n)^{h^+-h}((\lambda,z)\tau+n)^{h^--h}\cr }
$$
because
$$\eqalign{ (\lambda+n_2z,z_{v^+})=
(\lambda-{n+\tau(\lambda,z_{v^+})+\bar\tau(\lambda,z_{v^-})\over
2iyz_{v^+}^2}z,z_v^+)&= -{(\lambda,z)\bar\tau+n\over 2iy}\cr
(\lambda+n_2z,z_{v^-})=
(\lambda-{n+\tau(\lambda,z_{v^+})+\bar\tau(\lambda,z_{v^-})\over
2iyz_{v^-}^2}z,z_v^-)&= -{(\lambda,z)\tau+n\over 2iy}.\cr }
$$

If we substitute these expressions into the formula
for $\hat g(n)$ we
find that $\hat g(n)$ is equal to
$$\eqalign{
&{1\over \sqrt{2yz_{v^+}^2}}
\sum_{h^+,h^-}
\exp(-\Delta/8\pi y)(p_{w,h^+,h^-})(w(\lambda))\times\cr
&\times\sum_h {h!(-yz_{v^+}^2/\pi)^h\over (-2iy)^{h^++h^-}}{h^+\choose
h}{h^-\choose h}
((\lambda,z)\bar\tau+n)^{h^+-h}((\lambda,z)\tau+n)^{h^--h}\times\cr
&\times
\e\left(\tau\lambda_{w^+}^2/2+\bar\tau\lambda_{w^-}^2/2
-n(\lambda,(z_{v^+}-z_{v^-})/2z_{v^+}^2)-{|(\lambda,z)\tau+n|^2\over
4yiz_{v^+}^2}\right).\cr
}
$$
Inserting this into the formula giving $\theta_{M+\gamma}$
in terms of $\hat g$ proves lemma 5.1.

We can use this to express the theta function of $M$ in terms of
that of $K=L/\Z z$, where $L=M\cap z^\perp$.

\proclaim Theorem 5.2.
Suppose that $z$ is a primitive norm 0 vector of $M$ and
choose a vector $z'\in M'$ with $(z,z')=1$. We write $N$ for the smallest
positive value of the inner product of $z$ with something in $M$,
so that $|M'/M|=N^2|K'/K|$. If $c\equiv (\gamma,z)\bmod N$
then by abuse of notation
we write $K+\gamma-cz'$ for the   coset of $K$ in $K'$
given by $K'\cap(M+\gamma-cz')/\Z z$.
Let $\mu $ be the vector
$$\mu = -z'+z_{v^+}/2z_{v^+}^2+z_{v^-}/2z_{v^-}^2
$$
of
$L\otimes \R/z = K\otimes \R$.
Then
$$\eqalign{
&\theta_{M+\gamma}(\tau;v,p)\cr
=&{1\over \sqrt{2yz_{v^+}^2}}
\sum_{h\ge 0}
\sum_{h^+,h^-}
{h!(-yz_{v^+}^2/\pi)^h\over (-2iy)^{h^++h^-}}{h^+\choose h}{h^-\choose h}
(c\bar\tau+d)^{h^+-h}(c \tau+d)^{h^--h}\times\cr &\times\sum_{c\equiv
(\gamma,z)\bmod N\atop d\in \Z}
\e\left({-|c\tau+d|^2\over 4iyz_{v^+}^2}
-{(\gamma,z')d}+{(z',z')cd\over 2}\right)
\theta_{K+(\gamma-cz')}(\tau,\mu d,-c\mu,w,p_{w,h^+,h^-}).\cr
}
$$

Proof. We use lemma 5.1, and rewrite the sum over $M/z+\gamma$ using
the fact that every element $\lambda$ of $M/z+\gamma$ can be uniquely
written in the form $\lambda=\lambda_K+cz'$ with $\lambda_K\in
K+\gamma-cz'=(M/z+\gamma-cz')\cap z^\perp$ and $c\equiv
(\gamma,z)\bmod N$.  We find that
$$\eqalign{
&\theta_{M+\gamma}(\tau;v,p) = \cr
&{1\over \sqrt{2yz_{v^+}^2}}
\sum_{c\equiv (\gamma,z)\bmod N\atop d\in \Z}
\sum_{h^+,h^-}
\sum_{\lambda_K\in K+\gamma-cz'}
\e\left(-{|c\tau+d|^2 \over 4iyz_{v^+}^2}\right)\times\cr
&\times\sum_h {h!(-yz_{v^+}^2/\pi)^h\over (-2iy)^{h^++h^-}}
{h^+\choose h}{h^-\choose h}
\times\cr&\times
((\lambda_K+cz',z)\bar\tau+d)^{h^+-h}((\lambda_K+cz',z)\tau+d)^{h^--h}
\times\cr
&\times
\exp(-\Delta/8\pi y)(p_{w,h^+,h^-})(w(\lambda_K))\times\cr
&\times
\e\left(\tau(\lambda_K+cz')_{w^+}^2/2 +\bar\tau(\lambda_K+cz')_{w^-}^2/2
-(\lambda_K+cz',(z_{v^+}-z_{v^-})/2z_{v^+}^2)d\right)
\cr}
$$
so to prove theorem 5.2 we have to check that
$$\eqalign{
&\tau(\lambda_K+cz')_{w^+}^2/2 +\bar\tau(\lambda_K+cz')_{w^-}^2/2
-(\lambda_K+cz',(z_{v^+}-z_{v^-})/2z_{v^+}^2)d\cr
=&\tau(\lambda_K-c\mu)_{w^+}^2/2+\bar\tau(\lambda_K-c\mu)_{w^-}^2/2
-(\lambda_K-c\mu/2,\mu d)-(\lambda_K,z')d-cd(z',z')/2\cr
}
$$
(because $(\lambda_K,z')d=(\lambda-cz',z')d\equiv
(\gamma-cz',z')d\bmod 1$).

But $z'$ differs from $-\mu$ by multiples of $z_{v^+}$ and $z_{v^-}$
which have zero projections into $w^+$ and $w^-$, so we only have to
check that
$$
-(\lambda_K+cz',(z_{v^+}-z_{v^-})/2z_{v^+}^2)d
=
-(\lambda_K-c\mu/2,\mu d)-cd(z',z')/2-(\lambda_K,z')d.
$$
But  $\mu =-z'+z_{v^+}/2z_{v^+}^2+z_{v^-}/2z_{v^-}^2$,  so
$$\eqalign{
&-(\lambda_K+cz',(z_{v^+}-z_{v^-})/2z_{v^+}^2)d\cr
=&-(\lambda_K+cz',z_{v^+}/2z_{v^+}^2+z_{v^-}/2z_{v^-}^2)d\cr
=&-(\lambda_K+c(z'-z_{v^+}/2z_{v^+}^2-z_{v^-}/2z_{v^-}^2)/2,
z_{v^+}/2z_{v^+}^2+z_{v^-}/2z_{v^-}^2-z')d \cr
&-cd(z',z')/2-(\lambda_K,z')d\cr
=& -(\lambda_K-c\mu/2,\mu d)-cd(z',z')/2-(\lambda_K,z')d.\cr
}
$$
This proves theorem 5.2.

Suppose that $F_M=\sum_\gamma e_\gamma f_{M+\gamma}$ is a modular form
of type $\rho_M$ and weight $(-b^-/2-m^-,-b^+/2-m^+)$.
Define a
$\C[K'/K]$-valued function
$$F_K(\tau,\alpha,\beta)=
\sum_{\gamma\in K'/K}e_\gamma f_{K+\gamma}(\tau,\alpha,\beta)
$$
by putting
$$f_{K+\gamma}(\tau,\alpha,\beta) =
\sum_{\lambda\in M'/M\atop \lambda|L=\gamma}
\e(-(\lambda,\alpha z')-\alpha\beta (z',z')/2))
f_{M+\lambda+\beta z'}(\tau)
$$
for $\alpha,\beta\in \Z$, $\gamma\in K'/K$.
The notation $\lambda|L$ means the restriction of $\lambda\in \Hom(M,\Z)$
to $L$, and $\gamma\in \Hom(K,\Z)$ is considered an element of
$\Hom(L,\Z)$ using the quotient map from $L$ to $K$. The
elements of $M'$ whose restriction to $L$ is 0 are exactly the integer
multiples of $z/N$.
Therefore if $\lambda$ is one of the elements in the
sum above, then the remaining elements in the sum
are the elements $\lambda+nz/N$ for $n\in \Z/N\Z$,
and $\lambda^2/2\equiv \gamma^2/2\bmod 1$.
\proclaim Theorem 5.3. With notation as above, the function $F_K$
satisfies the transformation formula
$$\eqalign{
&F_K((a\tau+b)/(c\tau+d),a\alpha+b\beta,c\alpha+d\beta)\cr
=& (c\tau+d)^{-b^-/2-m^-}(c\bar\tau+d)^{-b^+/2-m^+}
\rho_K\left({ab\choose cd},\sqrt{c\tau+d} \right)
F_K(\tau,\alpha,\beta)\cr
}
$$
for all $({ab\choose cd},\sqrt{c\tau+d})\in Mp_2(\Z)$.

Proof. As in the proof of theorem 4.1,
it is sufficient to
check it for the standard generators $T=({11\choose 01},1)$ and
$S=({0-1\choose 1\phantom{-}0},\sqrt\tau)$.
 For the  generator $T$ we have to show that
$$f_{K+\gamma}(\tau+1,\alpha+\beta,\beta)
= \e((\gamma,\gamma)/2)f_{K+\gamma}(\tau,\alpha,\beta).
$$
We prove this as follows.
$$\eqalign{
&f_{K+\gamma}(\tau+1,\alpha+\beta,\beta)\cr
=&\sum_{\lambda\in M'/M\atop \lambda|L=\gamma}
f_{M+\lambda+\beta z'}(\tau)\e((\lambda+\beta z')^2/2)
\e(-(\lambda,(\alpha+\beta) z')-(\alpha+\beta)\beta (z',z')/2))\cr
=&\sum_{\lambda\in M'/M\atop \lambda|L=\gamma}
f_{M+\lambda+\beta z'}(\tau)\e({\lambda}^2/2)
\e(-(\lambda,\alpha z')-\alpha\beta (z',z')/2))\cr
=& \e((\gamma,\gamma)/2)f_{K+\gamma}(\tau,\alpha,\beta)\cr
}
$$

For the  generator $S$ we
have to show
$$\eqalign{
&\sqrt{i\bar\tau}^{b^+}\bar\tau^{m^+}\sqrt{\tau/i}^{b^-}\tau^{m^-}
\sqrt{|K'/K|}f_{K+\gamma}(-1/\tau,-\beta,\alpha)\cr
=&
\sum_{\delta\in K'/K}\e(-(\gamma,\delta))f_{K+\delta}(\tau,\alpha,\beta).\cr
}
$$

We prove this as follows.
$$\eqalign{
&\sqrt{i\bar\tau}^{b^+}\bar\tau^{m^+}\sqrt{\tau/i}^{b^-}\tau^{m^-}
\sqrt{|K'/K|}f_{K+\gamma}(-1/\tau,-\beta,\alpha)\cr
=&\sqrt{i\bar\tau}^{b^+}\bar\tau^{m^+}\sqrt{\tau/i}^{b^-}\tau^{m^-}
{\sqrt{|M'/M|}\over N}
\sum_{\lambda\in M'/M\atop \lambda|L=\gamma} f_{M+\lambda+\alpha z'}(-1/\tau)
\e((\lambda,\beta z')+\beta\alpha(z',z')/2)\cr
=&{1\over N}\sum_{\lambda\in M'/M\atop \lambda|L=\gamma}\sum_{\delta\in M'/M}
\e((-\lambda-\alpha z',\delta))
f_{M+\delta}(\tau)\e((\lambda,\beta z')+\beta\alpha(z',z')/2)\cr
=&\sum_{\delta\in M'/M\atop (\delta,z)=\beta }
\e((-\alpha z',\delta))f_{M+\delta}(\tau)\e(\beta\alpha(z',z')/2)
\e((\gamma,\beta z'-\delta))\cr
=&\sum_{\delta\in M'/M\atop (\delta,z)=0 }
\e((-\alpha z',\delta ))f_{M+\delta+\beta z'}(\tau)\e(-\beta\alpha(z',z')/2)
\e((-\gamma,\delta))\cr
=&\sum_{\delta\in K'/K\atop (\delta,z)=0}\e(-(\gamma,\delta))
\sum_{\lambda\in M'/M\atop \lambda|L=\delta}
\e(-(\lambda,\alpha z')-\alpha\beta(z',z')/2)f_{M+\beta z'+\lambda}(\tau)\cr
=&\sum_{\delta\in K'/K}\e(-(\gamma,\delta))f_{K+\delta}(\tau,\alpha,\beta)\cr
}
$$
This proves theorem 5.3.

\proclaim 6.~The singularities of $\Phi$.

We set up some notation for the rest of this paper.  We let $M$ be an
even lattice of signature $(b^+,b^-)$.  We write $z$ for a primitive
norm 0 vector of $M$ (if one exists) and write $z'$ for a vector of
$M'$ with $(z,z')=1$. We let $L$ be the singular lattice $M\cap
z^\perp$ and let $K$ be the nonsingular lattice $L/\Z z$. We can
identify $K\otimes \R$ with the orthogonal complement of $z$ and $z'$
in $M\otimes \R$, and hence can identify $K$ with a subset of
$M\otimes \R$ (but note that $K$ is not necessarily a subset of $M$ in
this identification, though it is if $z'\in M$).  We recall that $v$
is an isometry from $M\otimes \R$ to $\R^{b^+,b^-}$, so that $v^+$ is
an element of the Grassmannian $G(M\otimes\R)$.

We suppose that $F_M(\tau)=y^{b^+/2+m^+}F(\tau)$ is some
$\C[M'/M]$-valued function on the upper half plane $H$ transforming
under $SL_2(\Z)$ with weight $(-b^-/2-m^-,-b^+/2-m^+)$ and
representation $\rho_M$.  We write $f_{M+\gamma}$ for the component of
$F$ corresponding to $\gamma\in M'/M$, and we will usually assume that
$F$ can be written in the form
$$F(\tau)=\sum_{\gamma}e_\gamma f_\gamma(\tau)
=\sum_{\gamma}e_\gamma \sum_{n\in \Q}\sum_{k\ge 0}c_\gamma(n,k)\e(n\tau)y^{-k}
$$
for complex numbers $c_\gamma(n,k)$ which are zero for all but a
finite number of values of $k$ and for all sufficiently small values
of $n$. The functional equation for $F$ under $Z\in Mp_2(\Z)$ implies that
$f_{M-\gamma}=(-1)^{m^++m^-} f_{M+\gamma}$, so that
$c_{-\gamma}(n,k)=(-1)^{m^++m^-}c_\gamma(n,k)$.

We define $\Phi_M(v,p,F_M)$ by
$$
\Phi_M(v,p,F_M)=
\int_{SL_2(\Z)\\H} \bar \Theta_M(\tau;v,p)F_M(\tau)dxdy/y^2
$$
and define $\Phi_M(v,p,F)$ by
$$
\Phi_M(v,p,F)=\Phi_M(v,p,F_M)=
\int_{SL_2(\Z)\\H} \bar \Theta_M(\tau;v,p)F(\tau)y^{b^+/2+m^+}dxdy/y^2.
$$
(where we define complex conjugation in $\C[M'/M]$ by
putting $\bar e_\gamma=e_{-\gamma}$, and the product
of $\bar \Theta_M$ and $F_M$ means we take their inner product
using $(e_\gamma,e_\delta)=1$ if $\gamma+\delta=0$ and 0 otherwise.)
The power
of $y$ and the weight of $F$ are chosen so that the integrand has
weight $(0,0)$; recall that $\bar\Theta_M$ has weight
$(b^-/2+m^-,b^+/2+m^+)$, $y$ has weight $(-1,-1)$, and $dxdy$ has
weight $(-2,-2)$.

The integral is often divergent and has to be regularized as follows.
We integrate over the region $F_w$, where $F_\infty=\{\tau||\tau|\ge
1, |\Re(\tau)|\le 1/2\}$ is the usual fundamental domain of $SL_2(\Z)$
and $F_w$ is the subset of $F_\infty$ of points $\tau$ with
$\Im(\tau)\le w$.  Suppose that the $\lim_{w\rightarrow
\infty}\int_{F_w}F(\tau)y^{-s}dxdy/y^2$ exists for $\Re(s)>>0$ and can
be continued to a meromorphic function defined for all complex
$s$. Then we define $\int_{SL_2(\Z)\\H} F(\tau)dxdy/y^2$ to be the
constant term of the Laurent expansion of this function at $s=0$.
This regularized integral exists for much more general functions $F$;
it is sufficient that the function $\bar\Theta F$ should have a
Fourier series expansion in $x$, whose constant coefficient has an
asymptotic expansion whose terms are constants times complex powers of
$y$ times nonnegative integral powers of $\log(y)$.  (When there is a
pole at 0 this definition is a bit clumsy and it might be better to
define $\Phi$ as the residue at $s=0$ of $\Phi_M(v,p,-2FE^*(*,s))ds$
where $E^*(\tau,s)$ is a real analytic Eisenstein series with
a pole of residue $-1/2$ at $s=0$.)

If we have a function invariant under a subgroup $\Gamma$ of
$SL_2(\Z)$ of finite index then we can define its regularized
integral by first averaging it over $SL_2(\Z)/\Gamma$ to get a
function invariant under $SL_2(\Z)$, and then taking the regularized
integral of this average, but we will not need this in this paper.

The function $\Phi_M$ is invariant under $\sigma\in \Aut(M)$ in the
sense that $\Phi_M(\sigma(v),\sigma(p),\sigma(F)) = \Phi_M(v,p,F)$,
where the action on $F$ is given by the action on $M'/M$.  We define
$\Aut(M,F)$ to be the subgroup of $\Aut(M)$ fixing $F$.  If $p$ is the
constant function 1 then $\Phi_M(v,F)=\Phi_M(v,p,F)$ is a function on
the Grassmannian $G(M)$ that is invariant under $\Aut(M,F)$. For more
general functions $p$ we can interpret $\Phi_M$ as an
$\Aut(M,F)$-invariant section of an $O_M(\R)$-equivariant vector
bundle over $G(M)$ as follows.  Suppose that $V$ is some subspace of
the polynomials on $\R^{b^+,b^-}$ that is invariant under the action
of $O_{b^+,b^-}(\R)$. We define a vector bundle over $G(M)$ to be the
set of pairs $(v,p)\in Iso(M\otimes \R,\R^{b^+,b^-})\times V$ modulo
the action of the group $O_{b^+,b^-}(\R)$ (which acts on both
factors). This gives us a $O_M(\R)$-invariant vector bundle over
$G(M)$, which is of type $V$ in the sense that the action of the
stabilizer of a point of $G(M)$ on the fiber is a representation of
type $V$.  Now we can see that if $p\in V$ then the function
$\Phi_M(v,p,F)$ is just an invariant section of the dual vector bundle
of type $V^*$ (which is isomorphic to the bundle of type $V$ as $V$ is
self dual as a representation of $O_{b^+,b^-}(\R)$).

We will say that a function $f$ has singularities of type $g$ at a
point if $f-g$ can be redefined on a set of codimension at least 1 so
that it becomes real analytic near the point. In the rest of this
section we will find all singular points of $\Phi_M$ and find what
type of singularities $\Phi_M$ has at its singular points.

\proclaim Lemma 6.1. For real $r$ the function
$$f(r)=\int _1^\infty e^{-r^2y}y^{s-1} dy
= |r|^{-2s}\Gamma(s,r^2)
$$
has a singularity at $r=0$ of type $|r|^{-2s}\Gamma(s)$ unless $s$ is
a non-positive integer, in which case $f$ has a singularity of type
$(-1)^{s+1}r^{-2s}\log(r^2)/(-s)!$.

Proof. If $s>0$ then $\int _0^1 e^{-r^2y}y^{s-1} dy$ is nonsingular at $r=0$
(as can be seen by expanding $e^{-r^2 y}$ as a power series in $y$)
so $f$ has singularities of type
$$\int _0^\infty e^{-r^2y}y^{s-1} dy
 = (r^2)^{-s}\Gamma(s)=|r|^{-2s}\Gamma(s).$$
(Warning: note that $(r^2)^{-s}$ is not the same as $r^{-2s}$ for $r<0$
if $s$ is not an integer.)

If $s=0$ we integrate by parts to see that $f$ has singularities of type
$$r^2\int _1^\infty e^{-r^2y}\log(y) dy.$$ As $\log(y)$ is integrable
near 0 we can again change the range of integration to $[0,\infty]$
without affecting the type of the singularity. Changing $y$ to $y/r^2$
we see that the singularity has type
$$r^2\int_0^\infty e^{-y}\log(1/r^2)dy/r^2 = -\log(r^2).
$$

If $s<0$ we integrate by parts to see that $f$ has a singularity of type
$$r^2\int _1^\infty e^{-r^2y}y^{s}s^{-1} dy.$$
If $s$ is a negative integer this shows that $f$ has a singularity
of type
$$-{r^2\over s}{r^2\over s+1}\cdots {r^2\over -1}\log(r^2)
= (-1)^{s+1}r^{-2s}\log(r^2)/(-s)!.
$$
If $s$ is not a negative integer then we see by reducing to the case when
$s>0$ that $f$ has a singularity of type
$$|r|^{-2s}\Gamma(s).
$$

This completes the proof of lemma 6.1.

The singularities of $\Phi_M$ can be worked out using the method
of Harvey and Moore [H-M] as follows.
\proclaim Theorem 6.2.
Near the point $v_0\in G(M)$, the function
$\Phi_M(v,p,F)$ has a singularity of type
$$\eqalign{ &\sum_{\lambda\in M'\cap v_0^-\atop \lambda\ne 0}
\sum_{j,k} c_\lambda(\lambda^2/2,k)(1/j!)(-\Delta/8\pi)^j(\bar
p)(v(\lambda))
\times\cr&\times
(2\pi\lambda_{v^+}^2)^{1+j+k-b^+/2-m^+}\Gamma(-1-j-k+b^+/2+m^+)
\cr
}
$$
except that whenever $1+j+k-b^+/2-m^+$ is a non-negative integer
the corresponding term in the sum has to be replaced by
$$\eqalign{
&-c_\lambda(\lambda^2/2,k)(1/j!)(-\Delta/8\pi)^j(\bar p)(v(\lambda))
\times\cr&\times
(-2\pi\lambda_{v^+}^2)^{1+j+k-b^+/2-m^+}
\log(\lambda_{v^+}^2)/(1+j+k-b^+/2-m^+)!.
\cr
}
$$
In particular $\Phi_M$, considered as a section of a vector bundle
over $G(M)$, is nonsingular except along a locally finite set of
codimension $b^+$ sub Grassmannians (isomorphic to $G(b^+,b^--1)$) of
$G(M)$ of the form $\lambda^\perp = \{v_+|v_+\perp \lambda\}$ for some
negative norm vectors $\lambda\in M$.

Proof. The function $\Phi_M(v,p,F)$ is defined by an integral
$$\int_{y>0} \int_{|x|\le 1/2 \atop x^2+y^2\ge 1}
\bar\Theta(\tau;v,p)F(\tau)y^{b^+/2+m^+}dxdy/y^2.
$$
The integral over any compact region is a real analytic function of
$v$, so we may assume that the integral is taken over the region
$|x|\le 1/2$, $y\ge 1$ as this does not change the types of
singularities.  If we substitute in the definitions
$$\theta_{M+\gamma}(\tau;v,p) =\sum_{\lambda\in
M+\gamma}\exp(-\Delta/8\pi y)
p(v(\lambda))\e(\tau\lambda_{v^+}^2/2+\bar\tau\lambda_{v^-}^2/2)
$$
and $f_\gamma(\tau)=\sum_{n,k} c_\gamma(n,k)\e(n\tau)y^{-k}$
and carry out the integral over $x$
we get a sum of terms of the form
$$\sum_{\lambda\in M'\atop \lambda\ne 0}\sum_{j,k}{1\over j!}
(-\Delta/8\pi)^j(\bar p)(v(\lambda))c_\lambda(\lambda^2/2,k)\int_{y\ge 1}
 \exp(-2\pi y\lambda_{v^+}^2)y^{-2-j-k+b^+/2+m^+}dy
$$
plus a term for $\lambda=0$ which does not depend on $v$ and therefore
does not contribute to the singularity.  The singularities of these
terms occur only when $\lambda_{v^+}^2$ becomes 0, or in other words
when $\lambda\in v_0^-$. In this case the singularity of the integral
can be read off from lemma 6.1 with $r^2=2\pi \lambda_{v^+}^2$ and
$s=-1-j-k+b^+/2+m^+$.  This proves theorem 6.2.

Suppose that $b^+ = 1$, so that $G(M)$ is real hyperbolic space of
dimension $b^-$ and the singularities of $\Phi_M$ lie on hyperplanes
of codimension 1. Then the set of points where $\Phi_M$ is real
analytic is not connected, so we would like a wall crossing formula
telling us how the function changes as we pass through the singular
set.  The set of norm 1 vectors of $M\otimes \R$ has two components,
each isomorphic to $b^-$-dimensional hyperbolic space.  The components
of the points where $\Phi_M$ is real analytic are called the Weyl
chambers of $\Phi_M$. We will also call the positive cones generated
by these sets Weyl chambers. If $W$ is a Weyl chamber and $\lambda\in
M$ then $(\lambda,W)>0$ means that $\lambda$ has positive inner
product with all elements in the interior of $W$.

\proclaim Corollary 6.3 (The wall crossing formula). We use notation
as above.  Suppose that $\Phi_1 $ and $\Phi_2$ are the real analytic
restrictions of $\Phi_M$ to two adjacent Weyl chambers $W_1$ and
$W_2$, separated by a wall $W_{12}$. Then $\Phi_1$ and $\Phi_2$ can
both be extended to real analytic functions on the closure of the
union $W_1\cup W_2$, and their difference $\Phi_1(v)-\Phi_2(v)$ is
given by
$$ \eqalign{
\sum_{{\lambda\in M', \atop \lambda\perp W_{12}}
\atop (\lambda,W_1)>0}\sum_{j,k}
&{4\over j!}(-\Delta/8\pi)^j(\bar
p)(v(\lambda))c_\lambda(\lambda^2/2,k)\times\cr
&\times(\sqrt{2\pi}\times(\lambda,v_1))^{1+2j+2k-2m^+}
\Gamma(-1/2-k-j+m^+).\cr }
$$
where $v_1$ is a norm 1 vector in $W_1$ or $W_2$ generating $v^+$.
Moreover this expression is a polynomial in $v_1$
of degree at most $m^--m^++1+2k_{max}$, where $k_{max}$ is
the largest value of $k$ with some $c_\gamma(m,k)$ nonzero.

Proof. The formula for $\Phi_1-\Phi_2$ follows immediately from theorem 6.2,
because the function $\Phi_M$ has a singularity of type
$$ \eqalign{
\sum_{{\lambda\in M', \atop \lambda\perp W_{12}}}\sum_{j,k}
&{1\over j!}(-\Delta/8\pi)^j(\bar
p)(v(\lambda))c_\lambda(\lambda^2/2,k)\times\cr
&\times(\sqrt{2\pi}\times|\lambda_{v^+}|)^{1+2j+2k-2m^+}
\Gamma(-1/2-k-j+m^+)\cr }
$$
along the wall
$W_1\cap W_2$.
The terms of the sum for $\lambda$ and $-\lambda$ are the same because
$c_{-\lambda}=(-1)^{m^++m^-}c_\lambda$ and $p$ is homogeneous
of degree $m^++m^-$. So we may as well only sum over
the elements $\lambda$ with $(\lambda,W_1)>0$, in which case
$|\lambda_{v^+}| = (\lambda,v_1)$.
The factor of 4 appears because  we pick up a factor
of 2 by summing over only half the vectors of $M$, and another factor of
2 because the difference of $|x|$ and $-|-x|$ is $2x$ for $x>0$.

Now we have to prove that the expression above is a polynomial
in $(\lambda,v_1)$, or in other words that the terms with negative powers
of $(\lambda,v_1)$ all cancel out.
The function $p(v(\lambda))$ can be written as $(\lambda, v_1)^{m^+}
p^-(v(\lambda))$ for some polynomial $p^-$ of degree $m^-$
on $\R^{1,b^-}/\R^1=\R^{b^-}$.
Then
$$\eqalign{
&{1\over j!}(-\Delta/8\pi)^j(p)(v(\lambda))
\cr
=&\sum_{j^+,j^-\atop j^++j^-=j}
{1\over j^-!}
(-\Delta/8\pi)^{j^-}(p^-)(v(\lambda))
{m^+!\over (m^+-2j^+)!}(-1/8\pi)^{j^+}(1/j^+!)
(\lambda,v_1)^{m^+-2j^+}.
\cr
}
$$
The power of $(\lambda,v_1)$ in the terms with some fixed $\lambda$,
$k$ and $j^-$ in the wall crossing formula is $(1+2j-2m^++2k)+
(m^+-2j^+)= 2j^--m^++1+2k$, so we assume that this is negative and we
want to prove that the sum over $j^+$ vanishes.  In particular we then
have $-1/2-k-j^--j^++m^+>0$ because $2j^+\le m^+$.

We use the duplication formula $\Gamma(z+1/2)=
\sqrt\pi2^{1-2z}\Gamma(2z)/\Gamma(z)$ for the $\Gamma$ function to see
that if $z$ is a nonnegative integer then $\Gamma(z+1/2)=\sqrt\pi
(2z)!/2^{2z}z!$.  From this it follows that if $-1/2-k-j^--j^++m^+>0$
then
$$
\eqalign{
&\sum_{j^+}
{m^+!\over (m^+-2j^+)!}(-1/8\pi)^{j^+}(1/j^+!)
(\lambda,v_1)^{m^+-2j^+}
\times\cr&\times
(\sqrt{2\pi}(\lambda,v_1))^{2j^+}
\Gamma(-1/2-k-j^--j^++m^+)
\cr
&=
\sum_{j^+}
{
(\lambda,v_1)^{m^+}
(-1)^{j^+}m^+!
(-2-2k-2j^--2j^++2m^+)!\sqrt \pi
\over
j^+!(m^+-2j^+)!4^{j^+}
(-k-j^--j^++m^+-1)!2^{-2-2k-2j^--2j^++2m^+}
}
\cr
&= {\sqrt \pi (\lambda,v_1)^{m^+}(-2-2k-2j^-+m^+)!\over
(-1-k-j^-+m^+)!2^{-2-2k-2j^-+2m^+}}
\times\cr&\times
\sum_{j^+}
(-1)^{j^+}{-1-k-j^-+m^+\choose j^+}{ -2-2k-2j^--2j^++2m^+\choose m^+-2j^+}
}
$$
If we put
$A=m^+$, $C=-1-k-j^-+m^+$, $B=k+j^-$.
then $0\le B<C$ and $B+C<A$, so by lemma 14.1
the sum above vanishes.
Hence all the non-polynomial terms in the wall crossing formula
cancel out,
which  proves corollary 6.3.

One consequence of corollary 6.3 is that the difference
$\Phi_1-\Phi_2$ is a polynomial.  We will later use this to prove the
stronger result that $\Phi_1$ is itself a polynomial.  (Warning:
sometimes the function $\Phi_1$ is given by an odd polynomial in
$v_1$, but in spite of this $\Phi_M$ has the same value on $v_1$ and
$-v_1$.)

\proclaim Corollary 6.4. Suppose that in corollary 6.3 we take $p=1$ and
take $F$ to be holomorphic on $H$.
Then the difference $\Phi_1(v)-\Phi_2(v)$
is given by
$$
8\pi \sqrt 2\sum_{{\lambda\in M', \atop \lambda\perp W_{12}}\atop
(\lambda,W_1)>0} c_\lambda(\lambda^2/2)(\lambda,v_1).
$$

Proof. This is just a special case of corollary 6.3.

\proclaim 7.~The Fourier expansion of $\Phi$.

We calculate the Fourier expansion of the function $\Phi_M(v,p,F)$
recursively in terms of a similar function $\Phi_K$, where $K$ is a
lattice of signature $(b^+-1,b^--1)$.
There is a similar result in the nonsingular case in [R-S].

\proclaim Theorem 7.1.
Let $M$, $b^\pm$, $K$, $z$, $z'$, $p$, $m^\pm$ be defined as in
section 6.  Suppose
$$F_M(\tau)=\sum_{\gamma\in M'/M}e_\gamma\sum_{m\in \Q}
c_{\gamma,m}(y)\e(mx)
$$
is a modular form of weight $(-b^-/2-m^-,-b^+/2-m^+)$ and type
$\rho_M$ with at most exponential growth as $y\rightarrow +\infty$.
Assume that each function $c_{\gamma,m}(y)\exp(-2\pi |m|y)$ has an
asymptotic expansion as $y\rightarrow +\infty$ whose terms are
constants times products of complex powers of $y$ and nonnegative
integral powers of $\log(y)$.  If $z_{v^+}^2$ is sufficiently small
then the Fourier expansion of $\Phi_M(v,p,F_M)$ is given by the
constant term of the Laurent expansion at $s=0$ of the analytic
continuation of
$$\eqalign{
& {1\over \sqrt 2 |z_{v^+}|}
\sum_{h\ge 0}h!(z_{v^+}^2/4\pi)^h\Phi_K(w,p_{w,h,h},F_K)
+\cr+&
{ \sqrt 2\over |z_{v^+}|}
\sum_{h\ge 0}
\sum_{h^+,h^-}
{h!(-z_{v^+}^2/\pi)^h\over (2i)^{h^++h^-}}{h^+\choose h}{h^-\choose h}
\sum_{j}
\sum_{\lambda\in K'}
{(-\Delta)^j(\bar p_{w,h^+,h^-})(w(\lambda))\over (8\pi)^j j!}
\times\cr
&\quad\times
\sum_{n> 0}
\e((n\lambda,\mu))
n^{h^++h^--2h}
\sum_{\delta\in M'/M\atop \delta|L=\lambda}
\e(n(\delta,z'))
\times\cr
&\quad\times
\int_{y>0}c_{\delta,\lambda^2/2}(y)
\exp(-\pi n^2/ 2yz_{v^+}^2-\pi y(\lambda_{w^+}^2-\lambda_{w^-}^2))
y^{h-h^+-h^--s-j-5/2} dy
}
$$
(which converges for $\Re(s)>>0$ to a holomorphic functions of $s$
which can be analytically continued to a meromorphic function of all
complex $s$).

Remark. The conditions of this theorem hold for a wide class of forms
$F$. For example they hold for forms which are holomorphic on $H$ and
meromorphic at cusps, Maass wave forms, real analytic Eisenstein
series, Siegel theta functions, $\E_2(\tau)$, Zagier's function
$G(\tau)$, and any products of these functions.  The reason for the
conditions on $F_M$ is that the condition about the asymptotic
expansion implies that $\int_1^\infty c_{\gamma,m}(y)\exp(-2\pi
|m|y)y^{-s-1}dy$ converges for $\Re(s)>>0$ to a function which can be
analytically continued to a meromorphic function for all complex $s$,
so that certain integrals have well defined regularizations (as in
section 6). The condition about exponential growth is needed to
exchange the order of a sum and an integral later in the proof.  For
functions $F$ which are holomorphic on $H$ the condition about
exponential growth is equivalent to saying that $F$ is meromorphic at
the cusps.

Proof of theorem 7.1. We expand $\Theta_M(\tau;v,p)$
in the formula
$$
\Phi_M(v,p,F)=
\int_{SL_2(\Z)\\H} \bar \Theta_M(\tau;v,p)F(\tau)dxdy/y^{2+s}
$$
into a sum over $c,d\in \Z$ using
theorem 5.2. (In the formula above, we implicitly assume that we
analytically continue the function in $s$ and then take the constant term
at $s=0$.) The terms with $c=d=0$ vanish unless $h^+=h=h^-$
in which case we get the terms in theorem 7.1 involving $\Phi_K$.
We rewrite the remaining terms as follows.
Inserting the complex conjugate of the
series for $\Theta_M$ from theorem 5.2 gives
$$
\eqalign{
&
\int_{\tau\in SL_2(\Z)\\H}
{1\over \sqrt {2y}|z_{v^+}|}
\sum_{(c,d)\ne (0,0)}
\sum_{\gamma\in M'/M\atop (\gamma,z)\equiv c\bmod N}\cr
&
\sum_{h\ge 0}
\sum_{h^+,h^-}
{h!(-z_{v^+}^2/\pi)^h\over (2i)^{h^++h^-}}{h^+\choose h}{h^-\choose h}
\times\cr&\times
(c\tau+d)^{h^+-h}(c \bar\tau+d)^{h^--h}
\e\left({-|c\tau+d|^2\over 4iy z_{v^+}^2}\right)
\e\left({(\gamma,z')d}-{(z',z')cd\over 2}\right)\times\cr
&\times\bar\theta_{K+\gamma-cz'}(\tau,\mu d,-c\mu,w,p_{w,h^+,h^-})
f_{M+\gamma}(\tau) y^{h-h^+-h^-}dxdy/y^{2+s}
.\cr}
$$
 From now on we will fix $h^+$, $h^-$, and $h$, and drop the factor
of
$${h!(-z_{v^+}^2/\pi)^h\over \sqrt 2
|z_{v^+}|(2i)^{h^++h^-}}{h^+\choose h}{h^-\choose h}$$ (and remember
to put it back in at the end of the calculation).  We also change
$\gamma$ to $\gamma+cz'$ and substitute in the definition
$$f_{K+\gamma}(\tau,-d,c) =
\sum_{\lambda\in M'/M\atop \lambda|L=\gamma}
\e((\lambda,z'd)+cdz'^2/2)f_{M+\lambda+cz'}(\tau)
$$
of $f_{K+\gamma}$ and get
$$
\eqalign{
&
\int_{\tau\in SL_2(\Z)\\H}
\sum_{(c,d)\ne (0,0)}
\sum_{\gamma\in K'/K}
(c\tau+d)^{h^+-h}(c \bar\tau+d)^{h^--h}
\e\left({-|c\tau+d|^2\over 4iy z_{v^+}^2}\right)
\times\cr
&\times\bar\theta_{K+\gamma}(\tau,\mu d,-c\mu,w,p_{w,h^+,h^-})
f_{K+\gamma}(\tau,-d,c) y^{h-h^+-h^--1/2}dxdy/y^{2+s}
.\cr}
$$

We replace the sum over all $(c,d)\ne (0,0)$ by a sum over
$(nc,nd)$ with $c,d$ coprime and $n>0$ and get
$$
\eqalign{
&
\int_{\tau\in SL_2(\Z)\\H}
\sum_{c,d\atop (c,d)=1}\sum_{n>0}
(c\tau+d)^{h^+-h}(c \bar\tau+d)^{h^--h}n^{h^++h^--2h}
\e\left({-|c\tau+d|^2n^2\over 4iy z_{v^+}^2}\right)
\times\cr
&\times\bar\Theta_{K}(\tau,n\mu d,-nc\mu,w,p_{w,h^+,h^-})
F_{K}(\tau,-nd,nc) y^{h-h^+-h^--1/2}dxdy/y^{2+s}
.\cr}
$$

We use the transformation of $\Theta_K$ of weight
$(b^+/2-1/2+m^+-h^+,b^-/2-1/2+m^--h^-)$ (theorem 4.1) and $F_K$
of weight
$(-b^-/2-m^-,-b^+/2-m^+)$ (theorem 5.3)
and $y$ of weight $(-1,-1)$ under $Mp_2(\Z)$ to get
$$\eqalign{
\int_{\tau\in SL_2(\Z)\\H}
&
\sum_{n>0\atop{ab\choose cd}\in SL_2(\Z)/\Z}\!\!\!\!\!\!\!\!
\exp(-\pi n^2/ 2\Im({a\tau+b\over c\tau+d})z_{v^+}^2)
\bar\Theta_{K}({a\tau+b\over c\tau+d},n\mu,0,w,p_{w,h^+,h^-})
\times\cr&\times
n^{h^++h^--2h}
F_K({a\tau+b\over c\tau+d},-n,0)
\Im({a\tau+b\over c\tau+d})^{-1 / 2-h^+-h^-+h }dxdy/y^{2+s}
.\cr
}
$$

We want to replace the integral over a fundamental domain of
$SL_2(\Z)$ by an integral over a fundamental domain of $\Z$. This
would be trivial to justify if the final integral below were
absolutely convergent (in the region $|x|\le 1/2$). It is in general
exponentially divergent as $y$ increases or as $y$ tends to a rational
cusp because of the singularities of $F$ at these points, so we need
to justify this exchange of sum and integral.  We first note the it is
enough to show the final integral is absolutely convergent in the
region $|x|\le 1/2, y\le 1$, because although both integrals are
divergent for $y\ge 1$, they have the same divergences and are
regularized in the same way.  As $F_M$ has at most exponential growth
$\exp(Ay)$ for some constant $A$ as $y\rightarrow\infty$ and $F_M$ is
an automorphic form we see that for small $y$, $F_M$ is bounded by
$\exp(-(A+\epsilon)/y)$ for any positive $\epsilon$. The main point is
that if $z_{v^+}^2$ is sufficiently small (more precisely, less than
$2/\pi A$) then the term $\exp(-\pi / 2yz_{v^+}^2)$ is sufficiently
small near the cusps on the real line to kill off the rapid growth of
$F_M$ near these cusps and ensure that the integral over $y\le 1$ is
absolutely convergent.  Hence the exchange of order of the sum and
integral is valid for sufficiently small $z_{v^+}^2$. (It is not
always valid for larger values of $z_{v^+}^2$ and sometimes gives the
wrong answer in this case; see footnote 22 of [HM] for a comment on
this. This is also the reason why the functions in section 10 are only
piecewise polynomials and not polynomials.)  We can keep the term
$y^{-s}$ because the value of $\Im((a\tau+b)/(c\tau+d))^s$ at $s=0$ is
the same as that of $y^s$ at $s=0$, so the value of the integral over
$y\le 1$ at $s=0$ is not affected if we make this replacement.  Doing
this gives
$$\eqalign{
&2\int_{y>0}\int_{x\in \R/\Z}
\sum_{n>0}
\exp(-\pi n^2/ 2yz_{v^+}^2)
\bar\Theta_{K}(\tau,n\mu,0,w,p_{w,h^+,h^-})
\times\cr&\times
n^{h^++h^--2h}
F_K(\tau,-n,0) y^{h-h^+-h^--5/ 2-s}dxdy
.\cr
}
$$
(The extra factor of 2 in front comes from the fact that $SL_2(\Z)$
has a center of order 2 acting trivially on $H$.)  We can replace
$\Theta_{K}(\tau,n\mu,0,w,p_{w,h^+,h^-})$ by the series defining it
and interchange the summation and integral because if we remove a
finite number of terms from the sum then the sum of the absolute
values of the remaining terms has a convergent integral. So the
expression above is equal to
$$\eqalign{
&
2\sum_{\gamma\in K'/K}
\sum_{\lambda\in K+\gamma}
\int_{y>0}\int_{x\in \R/\Z}
\sum_{n>0}
n^{h^++h^--2h}
\sum_{j}{(-\Delta)^j(\bar p_{w,h^+,h^-})(w(\lambda))\over (8\pi)^j j!}
\times\cr&\times
\exp(-\pi n^2 / 2yz_{v^+}^2)
\overline{\e(\tau\lambda_{w^+}^2/2+\bar\tau\lambda_{w^-}^2/2 -(\lambda,n\mu))}
\times\cr&\times
e_{-\gamma}F_K(\tau,-n,0) y^{h-h^+-h^--j-5/2-s}dxdy
.\cr}
$$

Next we expand $f_{K+\gamma}(\tau,-n,0)$ as a series
$$\eqalign{
 &\sum_{\delta\in M'/M\atop \delta|L=\gamma}
f_{M+\delta}(\tau)\e((\delta,n z'))
=\sum_{\delta\in M'/M\atop \delta|L=\gamma}
\sum_{m}c_{\delta,m}(y)\e(mx)\e((\delta,n z'))\cr
}
$$
to get
$$
\eqalign{
&
2\sum_{\gamma\in K'/K}
\sum_{\lambda\in K+\gamma}
\sum_{\delta\in M'/M\atop \delta|L=\gamma}
\sum_{j}{(-\Delta)^j(\bar p_{w,h^+,h^-})(w(\lambda))\over (8\pi)^j j!}
\times\cr&\times
\int_{y>0}
\sum_{n>0}
\sum_{m}c_{\delta,m}(y)
\exp(-\pi n^2/ 2yz_{v^+}^2)
\exp(-\pi y\lambda_{w^+}^2+\pi y \lambda_{w^-}^2)
\times\cr
&\times
\e((n\lambda,\mu))\e(n(\delta,z'))
n^{h^++h^--2h}
\times\cr&\times
\int_{x\in \R/\Z}\e(-x\lambda_{w^+}^2/2-x\lambda_{w^-}^2/2+xm)dx\quad
 y^{h-h^+-h^--j-5/2-s}dy.\cr
}
$$

We carry out the integral over $x$, which is 0 unless
$m=\lambda_{w^+}^2/2+\lambda_{w^-}^2/2=\lambda^2/2$ to get
$$\eqalign{&
2\sum_{\lambda\in K'}
\sum_{\delta\in M'/M\atop \delta|L=\lambda}
\sum_{j}{(-\Delta)^j(\bar p_{w,h^+,h^-})(w(\lambda))\over (8\pi)^j j!}
\times\cr&\times
\int_{y>0}
\sum_{n>0}
n^{h^++h^--2h}\e((n\lambda,\mu))
\exp(-\pi n^2/ 2yz_{v^+}^2)
\exp(-\pi y\lambda_{w^+}^2+\pi y \lambda_{w^-}^2)
\times\cr&\times
\e(n(\delta,z'))
c_{\delta,\lambda^2/2}(y) y^{h-h^+-h^--j-5/2-s}dy.
\cr}
$$
This proves theorem 7.1.

We now calculate the integral over $y$ in theorem 7.1 in several cases.
When $F$ is almost holomorphic we have to distinguish
between the cases $\lambda_{w^+}=0$ (which is always
the case if $b^+=1$ or $\lambda=0$),
and the case $\lambda_{w^+}\ne 0$ (which is true for generic $v$
provided that $b^+>1$ and $\lambda\ne 0$).

\proclaim Lemma 7.2.
Suppose
$$F(\tau)=y^{-b^+/2-m^+}F_M(\tau)=\sum_{\gamma\in M'/M}\sum_{m\in
\Q}\sum_{k\ge 0} c_\gamma(m,k)\e(m\tau)y^{-k}e_{\gamma}
$$
is an almost holomorphic modular form of weight
$(b^+/2+m^+-b^-/2-m^-,0)$.
If $\lambda_{w^+}$ is nonzero
then the integral over $y$ at $s=0$ in theorem 7.1  is equal to
$$\eqalign{
&\sum_k 2c_\delta(\lambda^2/2,k)
\left({n\over 2|z_{v^+}||\lambda_{w^+}|}
\right)^{h-h^+-h^--j-k+b^+/2+m^+-{3/2}}
\times\cr&\times
K_{h-h^+-h^--j-k+b^+/2+m^+-3/2}(2\pi n|\lambda_{w^+}|/|z_{v^+}|).\cr }
$$
If $h-h^+-h^--j-k+b^+/2+m^+=1$ this is equal to
$$
\sqrt 2 c_\delta(\lambda^2/2)
{|z_{v^+}|\over n}
\exp(-2\pi n|\lambda_{w^+}|/|z_{v^+}|).
$$

Proof.
We know that
$$c_{\delta,\lambda^2/2}(y)
=\sum_kc_\delta(\lambda^2/2,k)y^{b^+/2+m^+-k}\exp(-2\pi
\lambda^2 y/2).
$$
and therefore when $s=0$ the integral over $y$ in theorem 7.1 is
$$\int_{y>0}\sum_kc_\delta(\lambda^2/2,k)
\exp(-\pi n^2/ 2yz_{v^+}^2-2\pi y\lambda_{w^+}^2)
y^{h-h^+-h^--j+b^+/2+m^+-k-5/2} dy.
$$
We can evaluate this integral using the formula
$$
\int_{y>0}\exp(-\beta/y-\alpha y)y^{\nu-1}dy
=2(\beta/\alpha)^{\nu/2}K_\nu(2\sqrt{\alpha\beta})
$$
([E, IT vol 1, p. 313, 6.3.17])
which is valid if $\alpha$ and $\beta$ are both positive.

So we put
$$
\eqalign{
\alpha&= -2\pi \lambda_{w^+}^2 \cr
\beta&=  -\pi n^2/2z_{v^+}^2  \cr
\nu&=    h-h^+-h^--j-k+b^+/2+m^+-3/2      \cr
}
$$
and find that the integral over $y$ has the value stated in the lemma.
The special case $h-h^+-h^--j-k+b^+/2+m^+=1$
 follows from the
fact the $K_{-1/2}(z)=K_{1/2}(z)=\sqrt{\pi/2z}\exp(-z)$.
This proves lemma 7.2.

\proclaim Lemma 7.3.
Suppose that $\lambda_{w^+}=0$ and
$$F(\tau)=y^{-b^+/2-m^+}F_M(\tau)=
\sum_{\gamma\in M'/M}\sum_{m\in \Q}\sum_{k\ge 0}
c_\gamma(m,k)\e(m\tau)y^{-k}
$$
is an almost holomorphic modular form of weight
$(b^+/2+m^+-b^-/2-m^-,0)$.
Then the integral over $y$ in theorem 7.1  is equal to
$$\eqalign{
&\sum_k
c_{\delta}(\lambda^2/2,k)
\left({\pi n^2\over 2z_{v^+}^2}\right)^{h-h^+-h^--s-j-k+b^+/2+m^+-3/2}
\times\cr&\times
\Gamma(-h+h^++h^-+j+k-b^+/2-m^++s+3/2).\cr
}
$$

Proof. We are given that
$$
c_{\delta,\lambda^2/2}(y) = \sum_k
c_\delta(\lambda^2/2,k)y^{b^+/2+m^+-k}\exp(-2\pi y \lambda^2/2).
$$
and $\lambda_{w^+}=0$ (so $\lambda=\lambda_{w^-}$)
and therefore
the integral over $y$ in 7.1 is equal to
$$\eqalign{
&\int_{y>0}\sum_k
c_{\delta}(\lambda^2/2,k)
\exp(-\pi n^2/2yz_{v^+}^2)y^{h-h^+-h^--j-k+b^+/2+m^+-s-3/2}dy/y\cr
=&\sum_k c_{\delta}(\lambda^2/2,k)
\left({\pi n^2\over 2z_{v^+}^2}\right)^{h-h^+-h^--s-j-k+b^+/2+m^+-3/2}
\times\cr&\times
\int_{y>0}
\exp(-y)y^{-h+h^++h^-+j+k-b^+/2-m^++s+3/2}dy/y\cr
=&\sum_k
c_{\delta}(\lambda^2/2,k)
\left({\pi n^2\over 2z_{v^+}^2}\right)^{h-h^+-h^--s-j-k+b^+/2+m^+-3/2}
\times\cr&\times
\Gamma(-h+h^++h^-+j+k-b^+/2-m^++s+3/2).\cr
}
$$
This proves lemma 7.3.

\proclaim 8.~Anisotropic lattices.

The calculation of the Fourier expansion in section 7 depends on the
existence of a primitive norm 0 vector in $M$, which automatically
exists whenever $M$ is indefinite and of dimension at least 5.  We
will later use this Fourier expansion to show that $\Phi_M$ has
various local properties (such as being locally a polynomial, or
holomorphic, or the real part of a holomorphic function), and we would
also like to show that $\Phi_M$ has these properties in dimensions at
most 4 when the lattice $M$ can be anisotropic. We show how to do this
in this section. The main idea is to embed the lattice $M$ in larger
lattices which have primitive norm 0 vectors, and then show that the
function $\Phi_M$ is a linear combination of restrictions of functions
associated to these larger lattices.

\proclaim Lemma 8.1 (the embedding trick).
Given $M$, $F$, and $p$ as in section 6 with $M$ not negative definite
we can write $\Phi_M(v,p,F)$
as a linear combination of functions
each of which is the restriction  to $G(M)$
of a function of the form
$$\Phi_{M\oplus M_j}(v,p,F_j)- singularities$$ where $M_j$ is a
nonzero negative definite even unimodular lattice, $F_j$ is a modular
form of the same type as $F$ and weight
$b^+/2+m^+-b^-/2-m^--\dim(M_j)/2$, and $p$ is extended by projecting
$\R^{b^+,b^-+\dim(M_j)}$ to $\R^{b^+,b^-}$.  The lattices $M\oplus
M_j$ contain primitive norm 0 vectors.

Proof. We will
write $\Phi_M$ as a difference of two functions as in the lemma.
We take $M_1$ and $M_2$ to be the Niemeier lattices with
root systems $A_2^{12}$ and $A_3^{8}$ with norms multiplied by $-1$,
and we take $F_1(\tau)=F_2(\tau)=
F(\tau)/24\Delta(\tau)$. Then $M_1$ and $M_2$ are even and self dual, the
functions $F_1$ and $F_2$ are modular forms of the same type as
$F$ because $\Delta$ is a modular form of level 1,
and
$$\bar\Theta_{M+M_2}(\tau)F_2(\tau)-\bar\Theta_{M+M_1}(\tau)F_1(\tau)
=\bar\Theta_M(\tau)F(\tau)
{\bar\Theta_{M_2}(\tau)-\bar\Theta_{M_1}(\tau)\over24\Delta(\tau)}
=\bar\Theta_M(\tau)F(\tau).
$$
If the functions $\Phi_{M\oplus M_j}(v,p,F_j)$ has no singularities
along $G(M)$ this would show that it was equal to the difference
$\Phi_{M\oplus M_2}(v,p,F_2)-\Phi_{M\oplus M_1}(v,p,F_1)$ and lemma
8.1 would now follow from this. In general they do have singularities
corresponding to the negative norm vectors of $M_j$, but we can get
round this by first subtract the singularities corresponding to nonzero
vectors of $M_j$ (given in theorem 6.2) before restricting.

We also need to check that that lattices $M\oplus M_j$ contain
primitive norm 0 vectors. But this follows immediately from the fact
that these lattices are indefinite lattices (as $M$ is not negative
definite) of dimension at least 5 (as $M_j$ has dimension at least 8)
and therefore have nonzero norm 0 vectors. (It is not hard to prove
the slightly stronger statement that the lattices $M\oplus M_j$
contain $II_{1,1}$ as a direct summand.)  This proves lemma 8.1.

\proclaim 9.~Definite lattices.

Theorem 7.1 reduces the calculation of the Fourier
expansion of $\Phi_M$ to the case of positive or negative definite lattices
$M$.
In this section we show how to do these calculations in some cases.
We get two different cases depending on whether the lattice $M$
is positive  or negative definite.

In the case of a negative definite lattice
we have to evaluate $\int \bar\Theta_M(\tau,p)F(\tau)dxdy/y^2$
where $\bar\Theta_M$ is a holomorphic modular form of weight $(m^-+b^-/2,0)$,
so that $F=F_M$ is a modular form of weight $(-b^-/2-m^-,0)$.
In particular $\bar\Theta_MF$ is an almost holomorphic modular
function. So it is sufficient to evaluate the integral
of any almost holomorphic modular function.

The main idea for evaluating $\int F(\tau)dxdy/y^2$ for a modular
function $F$ is to write $F(\tau)dxdy/y^2$ as $d(\omega d\tau)=
2i{\d\omega\over \d\bar\tau}dxdy$ for some modular form $\omega$ of
weight $(2,0)$.  Then we can convert the integral of $F(\tau)dxdy/y^2$
over the subset $F_w$ of the fundamental domain $F_\infty$ of
$SL_2(\Z)$ into an integral of $2i{\d\omega\over \d\bar\tau}dxdy$ over
the line $y=w$, $|x|\le 1/2$, which can usually be evaluated
explicitly.

We will show how to evaluate
$$\int_{SL_2(\Z)\\H} F(\tau)dxdy/y^2
$$
where $F$ is an almost holomorphic function invariant under $SL_2(\Z)$
(possibly with singularities at cusps) and where the divergent integral
is regularized as in section 6.
Recall that $\E_2(\tau)= E_2(\tau) - 3/\pi y$ is a non holomorphic
modular form of weight 2.

\proclaim Lemma 9.1. Any
almost holomorphic $SL_2(\Z)$-invariant function $g$
is a linear combination of functions of the form
$F(\tau)\E_2(\tau)^n$ where $F$ is a holomorphic function
(possibly singular at cusps)
transforming like a modular form of weight $-2n$.

Proof. The coefficient $F(\tau)$ of the highest power $y^{-n}$ of $1/y$
transforms like a modular form of weight  equal to weight$(g)-2n$,
so we can subtract a multiple of $F(\tau)\E_2(\tau)^n$ to reduce
this highest power. Lemma 9.1 now follows by induction on this highest
power.

\proclaim Theorem 9.2. The regularized divergent integral
$$\int_{SL_2(\Z)\\H} F(\tau)\E_2(\tau)^ndxdy/y^2
$$
for $F$ a modular form  of weight $-2n$ which is holomorphic on $H$
is  the constant term
of $$E_2(\tau)^{n+1}F(\tau)\pi/3(n+1).$$

Proof. We reproduce the proof of this
given in [L-S-W].
The main point is that $F(\tau)\E_2(\tau)^ndxdy/y^2$
is an exact differential, equal to $-d(F(\tau)
\E_2(\tau)^{n+1}\pi d\tau/3(n+1))$ because $${\d \E_2(\tau)\over \d \bar\tau}
= (1/2)({\d\over \d x}+i{d\over \d y})(-3/\pi y)=3i/2\pi y^2$$
and $d\bar\tau d\tau = 2idxdy$.
This implies that the integral over $F_w$ is equal (by Stokes'
theorem) to
$$
\int_{ x= 1/2+iw}^{-1/2+iw}-F(\tau)
\E_2(\tau)^{n+1}\pi d\tau/3(n+1).
$$
  This integral is the
constant term of
$$F(\tau)\E_2(\tau)^{n+1}\pi/3(n+1).
$$
The constant terms involving negative powers of $y$ tend to zero as
$w$ tends to $+\infty$ so we can drop them, and find that the
regularized integral is the constant term of
$E_2(\tau)^{n+1}F(\tau)\pi/3(n+1)$.  This proves theorem 9.2.

\proclaim Corollary 9.3. If $K$ is a negative definite even lattice
of dimension $b^-$ and $F$ is a modular form of weight $(-b^-/2,0)$ of
type $\rho_K$ which is holomorphic on $H$ and meromorphic at the cusps
then
$$\Phi_K(\cdot,1,F)= {\pi\over 3}\times\hbox{constant term of~}
\bar\Theta_{K}FE_2.$$

Proof. This follows immediately from theorem 9.2 and the fact
that $\Phi_K(\cdot,1,F)$ is by definition equal to the regularized
integral of $\bar\Theta_KFdxdy/y^2$ over a fundamental domain.
This proves corollary 9.3.

We now discuss the case of a positive definite lattice $M$, where we
have to evaluate $\int
\bar\Theta_M(\tau,p)F(\tau)y^{b^+/2+m^+}dxdy/y^2$ where
$\Theta_M(\tau,p)$ is now an almost holomorphic modular form of weight
$(b^+/2+m^+,0)$ and $F$ is an almost holomorphic modular form of
weight $(-b^+/2-m^+,0)$.  If $F$ is an almost holomorphic cusp form
this is just the usual Peterson inner product of $\Theta_M$ and $F$,
for which there is in general no known explicit finite elementary
formula.  If the lattice $M$ is one dimensional and generated by a
vector of norm $2N>0$ and $m^+=0$ we can usually evaluate the integral
$$\int  \bar\Theta_M(\tau) F(\tau) y^{-3/2}dxdy$$
explicitly using Zagier's non holomorphic modular form $\G(\tau)$
of weight $3/2$.
We recall the basic properties of $G$ from [Z].
We write $H(n)$ for the Hurwitz class number of $n$, so that
$$G(\tau)=\sum_nH(n)q^n= -1/12
+q^3/3+q^4/2+q^7+q^8+q^{11}+(4/3)q^{12}+O(q^{15}).$$ The function $\G$
is defined by
$$\eqalign{
\G(\tau)&=\sum_{n}H(n)q^n +
\sum_{n}q^{-n^2}
{1\over 16\pi\sqrt y}\int_{1\le u\le \infty}\exp(-4\pi un^2y)du/u^{3/2}\cr
&=\sum_{n}H(n)q^n +
{1\over 16\pi}
\sum_nq^{-n^2}
\int_{y\le u\le \infty}\exp(-4\pi un^2)du/u^{3/2}\cr
}
$$
so that
$$\eqalign{
{\d \G\over \d\bar\tau}
&= {1\over 16\pi}
\sum_nq^{-n^2}
{\d\over \d\bar\tau}\int_{y\le u\le \infty}\exp(-4\pi un^2)du/u^{3/2}\cr
&= {-i\over 32\pi}
\sum_n\exp(2\pi i (-n^2(x+iy)))\exp(-4\pi n^2y)/y^{3/2}\cr
&= {-i\over 32\pi}
y^{-3/2}\sum_n\bar\e(n^2\tau).\cr
}
$$
Zagier showed that  $\G$ is a modular form for the group $\Gamma_0(4)$
of weight $3/2$.

We now convert $\G$ into a level 1 modular form $\G_1$ of type
$\rho_M$ where $M$ is generated by a vector of norm $2$.
\proclaim Lemma 9.4. There is a function $\G_1$ with the following properties.
\item {1.} $\G_1(\tau)$ is a (non holomorphic) modular form
of type $\bar\rho_M$ and weight $(3/2,0)$.
\item {2.} ${\d \G_1(\tau)\over \d\bar\tau}
= {-i\over 16\pi}y^{-3/2}\bar\Theta_M(\tau)$.

Proof.
If $\sum_nc(n,y)\e(nx)$ is any modular form of weight $3/2\bmod 2$
satisfying the ``plus space'' condition $c(n,y)=0$ for $y\ne 0,-1\bmod 4$
then $e_0(\sum_nc(4n,y)\e(nx))+e_1\sum_nc(4n-1,y)\e((n-1/4)x)$
is a vector valued modular form of type $\bar\rho_M$, as follows
from the proof of theorem 5.4 of [E-Z]. (Formula (16) on page 64 of
[E-Z] is equivalent to the definition of a modular form of
type $\bar\rho_M$.)
If we apply this construction to $\G$ we get a modular
form $\G_1$ of weight $3/2$ and type $\bar\rho_M$. The statement
about the derivative of $\G_1$ follows from the calculation of the derivative
of $\G$ above and the fact that if we apply this construction
to $\sum_n\bar \e(n^2\tau)$ we get $\bar\Theta_M(\tau)$
This proves lemma 9.4.

\proclaim Lemma 9.5.
Let $N$ be a positive integer and let $M$ be the lattice generated by
a vector of norm $2N$.  There is a function $\G_N$ with the following
properties.
\item {1.} $\G_N(\tau)$ is a (non holomorphic) modular form
of type $\bar\rho_M$ and weight $(3/2,0)$.
\item {2.} ${\d \G(\tau)\over \d\bar\tau}= {-i\sqrt N\over 16\pi}
y^{-3/2}\bar\Theta_M(\tau)$.

Proof. Roughly speaking, we construct $\G_N$ from $\G_1$ in the same
way that $\Theta_{2N}$ is constructed from
$\Theta_2=\theta_0e_0+\theta_1e_1$ (where we write $\Theta_{2N}$ for
the theta function of a lattice generated by a vector of norm $2N$).
Define $U$ to be the representation of dimension
$2|SL_2(\Z)/\Gamma_0(N)|=2\sigma_1(N)$ obtained by restricting
$\rho_2$ to the double cover of $\Gamma_0(N)$ and then inducing back
up to $Mp_2(\Z)$ and define $V$ to be the space spanned by all the
functions of the form $\theta_j((a\tau+b)/d)$ for $ad=N$, $j=0,1$ so
there is a natural map from $U$ onto $V$.  There is also a map from
the space of $\rho_M$ to $V$ given by taking $\Theta_M$ to its
components (usually neither injective nor surjective). As all these
representations are completely reducible we can find a morphism from
$U$ to $\rho_M$ such that the image of the functions
$\theta_j((a\tau+b)/d)$ is $\Theta_M$. We define $\G_N$ to be the
image of the components of $\G_1((a\tau+b)/d)$ under the complex
conjugate of this map. This proves lemma 9.5.

\proclaim Corollary 9.6.
Let $M$ be a one dimensional lattice generated by a vector of norm $2N>0$.
Suppose $F$ is a modular form of type
$\rho_M$ and weight $(1/2,0)$ which is holomorphic on $H$. Then the
regularized integral
$$\int_{SL_2(\Z)\\H} \bar\Theta_M(\tau) F(\tau)y^{-3/2}dxdy$$
is equal to the the constant term of
$${-8\pi\over \sqrt N}F(\tau)G_N(\tau)$$
where $G_N$ is the holomorphic part of $\G_N$.

Proof. As $F$ is holomorphic  on $H$ we see from lemma 9.5 that
$${-i\sqrt N\over 16\pi}\bar\Theta_M(\tau) F(\tau)y^{-3/2}2idxdy
= d\left(\G_N(\tau)F(\tau)d\tau\right)
$$
and we can now complete the proof as in theorem 9.2. This proves
corollary 9.6.

\proclaim 10.~Lorentzian lattices.

We work out the functions $\Phi_M$ in the case when $M$ is Lorentzian.
There are two obvious differences between this case and other lattices.
Firstly the projection $\lambda_{v^+}$ is always 0,
while for most other lattices it is generically nonzero, so we
need to use the alternative formula of lemma 7.3 rather than that of
lemma 7.2.
Secondly the singularities of $\Phi_M$ occur on sets of codimension
$b^+=1$ so the set of nonsingular points is disconnected. This means
the formula for $\Phi_M$ we get is only valid in one component
of the nonsingular points, and
we get  ``wall crossing formulas'' telling us how
the formula for $\Phi_M$ changes as we cross a singular hypersurface.

The main result (theorem 10.3) of this section is that if $m^+=0$ then
$\Phi_M$ is the restriction of a piecewise polynomial on the
Lorentzian space.  We can divide Lorentzian (or hyperbolic) space up
into a system of ``Weyl chambers'' (which are sometimes, but not
usually, Weyl chambers for a reflection group) and on each Weyl
chamber $\Phi_M$ is given by a polynomial.  The polynomials of
adjacent Weyl chambers are related by a ``wall crossing formula''
similar to the one appearing in Donaldson theory. When the polynomial
is linear it is given by taking inner products with some vector,
called a Weyl vector, which is sometimes the Weyl vector for a
generalized Kac-Moody algebra.  We prove the main result by
calculating $\Phi_M$ in a Weyl chamber explicitly and finding that it
is given by a rational function, possibly with a pole along a certain
divisor. We show this pole does not exist by calculating $\Phi_M$ via
different Weyl chamber and finding that if we add a polynomial to it
then its only possible singularity lies on a different divisor.

We choose notation as in section 6, and we take $M$ to be a
Lorentzian lattice so that $K$ is negative definite. We let $C$ be one
of the two cones of positive norm vectors of $K\otimes R$, and call it
the positive cone. We assume that $z$ (if it exists) is in the closure
of $C$. We can identify the Grassmannian $G(M)$ with the norm 1
vectors in $C$. If a vector $v\in G(M)$ is represented by
a  norm 1 vector $v_1=(m\mu,m,n)\in M\otimes \R$ then
a short calculation shows that $\mu$ is the same as the vector $\mu$
in section 5, and the other things in 7.1 are given by
$$\eqalign{
z_{v^+}=(z,v_1)v_1/v_1^2&=m(m\mu,m,n)\cr
z_{v^+}^2&= m^2\cr
w^+& =0.\cr
}
$$
We put $p(x)= x_1^{m^+}p^-(x)$ where $p^-$ has degree $(0,m^-)$.
So $p_{w,h^+,h^-}= 0$ if $h^+\ne m^+$, and we put $p^-_{w,h^-}=p_{w,m^+,h^-}$.

We need to use Bernoulli polynomials whose properties we now recall
from [E, vol. 1, 1.13].  We define $\B_m(x)$ by $\B_m(x)=0$ if $m<0$,
$\B_0(x)=1$ and
$$\B_m(x)=-m!\sum_{n\ne 0}\e(nx)/(2\pi i n)^m
$$
for $m>0$
and they have the properties
\item{1.} $\B_m(x+1)=\B_m(x)$
\item{2.} $\B_m'(x) = m\B_{m-1}(x) $ for $x\notin \Z$ or $m\ne 1,2$.
\item{3.} If $0\le x<1$ then $\B_m(x)$
is the Bernoulli polynomial $B_m(x)$ of degree $m$ unless $x=0$ and
$m=1$ in which case $\B_1(0)=0$, and in particular $B_0(x)=1$,
$B_1(x)=x-1/2$ for $x\ne 0$, $B_2(x)= x^2-x+1/6$, and
$B_3(x)=x^3-3x^2/2+x/2$.
\item{4.} If $m\ge 0$ then
$B_m(x+1)-B_m(x)=mx^{m-1}$; in other words $\B_m(x)$ ``jumps down''
by $mx^{m-1}$ as $x$ crosses the origin.

\proclaim Lemma 10.1. If $m\in \Z$ and $x\in \R$
then the analytic continuation of the function
$$\sum_{n\ne 0} {\e(nx)\over n^m|n|^s}$$
is meromorphic for all $s\in \C$ and its value at $s=0$
is 0 if $m<0$ and $-(2\pi i)^m\B_m(x)/m!$ if $m\ge 0$.

Proof. This follows  from [E vol 1, 1.11, formulas 14,17,18]
because the expression above is equal to $F(\e(x),m+s)+(-1)^mF(\e(-x),m+s)$
in the notation of [E].
(There is a misprint in formula (18): the factor $(2\pi i)$
should be $(2\pi i )^m$.)

We can now find a finite formula for $\Phi_M(v,p,F_M)$.

\proclaim Theorem 10.2. Take notation as above.
If $v$ is represented by a norm 1 vector $(\mu m,m,n)$ in a Weyl
chamber $W$, then $\sqrt 2 |z_{v^+}|\Phi_M(v,p,F_M)$ is given by
$$
\eqalign{
& {m^+!(z_{v^+}^2)^{m^+}\over (4\pi)^{m^+}}\Phi_K(w,p_{w,m^+,m^+},F_K)+\cr
+&
\sum_{h\ge 0}
\sum_{h^-} {h!(-z_{v^+}^2/\pi)^h\over (2i)^{m^++h^-}}{m^+\choose h}
{h^-\choose h}
\sum_{j}
\sum_{\lambda\in K'}
{(-\Delta)^j(\bar p_{w,m^+,h^-})(w(\lambda))\over (8\pi)^j j!}
\times\cr
&\quad\times
\sum_{\delta\in M'/M\atop \delta|L=\lambda}
\sum_k
c_{\delta}(\lambda^2/2,k)
\times\cr
&\quad\times
\left({\pi \over 2z_{v^+}^2}\right)^{h-h^--j-k-1}
\Gamma(-h+h^-+j+k+1)
\times\cr&\times
{-(2\pi i)^{-m^++h^-+2j+2k+2}\B_{-m^++h^-+2j+2k+2}
((\lambda,\mu)+(\delta,z'))
\over (-m^++h^-+2j+2k+2)!}.
}
$$

Proof. In the case of Lorentzian lattices the projection
$\lambda_{v^+}$ is always 0 and $h^+=m^+$. This means that the formula
for $\Phi_M$ in theorem 7.1 and lemma 7.3 can be simplified as
follows.
$$
\eqalign{
&\sqrt 2 |z_{v^+}|\Phi_M(v,p,F_M)\cr
=& {m^+!(z_{v^+}^2)^{m^+}\over (4\pi)^{m^+}}\Phi_K(w,p_{w,m^+,m^+},F_K)+\cr
+&
2\sum_{h\ge 0}
\sum_{h^-}
{h!(-z_{v^+}^2/\pi)^{h}\over (2i)^{m^++h^-}}{m^+\choose h}{h^-\choose h}
\sum_{j}
\sum_{\lambda\in K'}
{(-\Delta)^j(\bar p_{w,m^+,h^-})(w(\lambda))\over (8\pi)^j j!}
\times\cr
&\quad\times
\sum_{\delta\in M'/M\atop \delta|L=\lambda}
\sum_k
c_{\delta}(\lambda^2/2,k)
\times\cr
&\quad\times
\left({\pi \over 2z_{v^+}^2}\right)^{h-h^--s-j-k-1}
\Gamma(-h+h^-+j+k+s+1)
\times\cr&\times
\sum_{n> 0}
\e((n\lambda,\mu))
\e(n(\delta,z'))
n^{m^++h^--2h}
(n^2)^{h-h^--s-j-k-1}.\cr
}
$$
If we change the sign of $n$, $\lambda$, and $\delta$ in the summand
above it is unaltered because $c_{-\delta}=(-1)^{m^-+m^+}c_\delta$,
$p_{w,m^+,h^-}(w(-\lambda))=(-1)^{m^-+h^-}p_{w,m^+,h^-}(w(\lambda))$,
and $(-n)^{h^-+m^+-2h}=(-1)^{h^-+m^+}n^{h^-+m^+-2h}$. Therefore we can
replace the sum over $n>0$ above by half the sum over $n\ne 0$ without
affecting the value of the expression.

We can evaluate the sum over $n\ne 0$  in finite terms using  Bernoulli
polynomials and lemma 10.1 as follows:
$$\eqalign{
&{1\over 2}\sum_{n\ne 0}
\e((n\lambda,\mu))
\e(n(\delta,z'))
n^{m^++h^--2h}
(n^2)^{h-h^--s-j-k-1}\cr
=&{-(2\pi i)^{-m^++h^-+2j+2k+2}\B_{-m^++h^-+2j+2k+2}
((\lambda,\mu)+(\delta,z'))
\over 2(-m^++h^-+2j+2k+2)!}.
}
$$
Substituting this into the formula for $\Phi_M$
proves theorem 10.2.

We now show that roughly speaking there is a lot of unexpected
cancelation in the sum for $\Phi_M$.

\proclaim Theorem 10.3.  On the interior of each
Weyl chamber $\Phi_M$ is the restriction of
a polynomial  on $M\otimes \R$
of degree at  most $m^--m^++2k_{max}+1$
where $k_{max}$ is the largest value of
$k$ with a nonvanishing coefficient $c_\gamma(n,k)$,
(and is 0 if this degree is negative).

Proof. By the embedding trick of section 8 we can embed $M$ in larger
lattices which have a norm 0 vector $z$, and also have the property
that the only polynomials fixed by $\Aut(M,F,C)$ are powers of
$(\lambda,\lambda)$ and are therefore constant on the space of norm 1
vectors, and can write $\Phi_M$ as a linear combination of similar
functions on these larger lattices. Hence we may assume that $M$ has
these properties.

We already know that $\Phi_M(v_1)$ restricted to a Weyl
chamber extends to a  rational function
which is a quotient of a polynomial by a power of $z_{v^+}^2=(v_1,z)^2$
by theorem 10.2, because  $(\lambda,\mu)=(\lambda,v_1)/|z_{v^+}|$.
So we first
have to show that this rational function does not have a pole along
$|z_{v^+}|=(z,v_1)=0$. Suppose that we choose two Weyl chambers
in the same positive cone
containing 2 different primitive norm 0 vectors $z_1$ and $z_2$, and
let $\Phi_1$ and $\Phi_2$ be the rational functions that restrict to
$\Phi_M$ on the 2 Weyl chambers. (If we have one primitive
norm 0 vector $z_1$ we can always find another linearly independent one $z_2$
as follows:
take any lattice vector which has nonzero inner product with $z_1$,
add a rational multiple of $z_1$ to make its norm zero, then multiply
it by a rational number to make it primitive.) We know that $\Phi_1(v)$ and
$\Phi_2(v)$ are rational functions whose only poles
lie on the two different irreducible divisors. On
the other hand by
the wall crossing formula 6.3 we know that $\Phi_1-\Phi_2$ is
a polynomial, so that $\Phi_1$ and $\Phi_2$ have the same singularities.
As they have no singularities in common they must both be polynomials.
Hence $\Phi_M$ is the restriction of a polynomial on Lorentzian space.

Now we need to bound the degree of this polynomial.  We note that by
the wall crossing formula 6.3 the functions $\Phi_1$ and $\Phi_2$ on 2
Weyl chambers differ by a polynomial of degree $m^--m^++2k_{max}+1$,
so it is sufficient to show that this is a bound on the degree in some
Weyl chamber.  Firstly, if $m^--m^++2k_{max}+1<-1$ then the Bernoulli
polynomials vanish in theorem 10.2 because $h^-+2j\le m^-$, so
$\Phi_M$ is zero.  Secondly, if $m^--m^++2k_{max}+1=-1$ then the
Bernoulli polynomials in theorem 10.2 are constant, so the expression
in theorem 10.2 is a polynomial vanishing to order at least 2 along
$z_{v^+}=0$. On the other hand $\Phi_M$ is constant because it is a
polynomial invariant under $\Aut(M,F,C)$ up to addition of constants
by the wall crossing formula, so if it is nonzero then $\sqrt
2|z_{v^+}|\Phi_M(v,p,F_M)$ vanishes to order at most 1 along
$z_{v^+}=0$. This is a contradiction, so $\Phi_M$ must be identically
zero. Thirdly, suppose that $m^--m^++2k_{max}+1\ge 0$. Then the image
of $\Phi_M$ is the space of polynomial functions on hyperbolic space
modulo the polynomial functions of degree at most $m^--m^++2k_{max}+1$
is fixed under $\Aut(M,F,C)$ by the wall crossing formula, and hence
must be a constant.  So $\Phi_M$ is the sum of a constant and a
polynomial of degree at most $m^--m^++2k_{max}+1\ge 0$ and hence is a
polynomial of degree at most $m^--m^++2k_{max}+1$.

This completes the proof of
theorem 10.3.

In particular if $p=1$ and $F$ is holomorphic on $H$, so that $k=0$
and $m^+=0$, then by theorem 10.3 $\Phi_M(v,1,F)$
extends to a homogeneous piecewise linear function.
In this case we define the Weyl vector $\rho(M,W,F)$ by
$$8\sqrt 2 \pi(\rho(M,W,F),v)=|v|\Phi_M(v/|v|,1,F)$$
where $v$ is a  vector in the Weyl chamber $W$.
The factor of $8\sqrt 2 \pi$ is put in to give the Weyl vector
good integrality properties, and is the constant in corollary 6.4.
We now derive an explicit formula for the Weyl vector of
any Weyl chamber whose closure contains a primitive norm 0 vector.
Suppose that $z$ is a norm 0 vector of $M$ and $z'$ is a vector
of $M'$ with $(z,z')=1$. We write vectors
of $M\otimes \R$ in the form $(v,m,n)$ with $v\in K\otimes \R$,
$m,n\in \R$, so that $(v,m,n)$ has norm $v^2+2mn +m^2z'^2$.
We have $z=(0,0,1)$, $z'=(0,1,0)$. (Warning: if $v\in K$
this does not imply that $(v,0,0)\in M$.)
Choose a Weyl chamber $W$ of $M$ whose closure contains $z$,
and write $(v,W)>0$ if $v\in M\otimes \R$ has positive
inner product with all elements in the interior of $W$.

\proclaim Theorem 10.4. The Weyl vector $\rho(M,W,F)$ is equal to
$(\rho,\rho_{z'},\rho_z)=\rho+\rho_{z'}z'+\rho_zz$ with
$$\eqalign{
\rho&=-{1\over 2} \sum_{\lambda\in K',(\lambda,0,0)\in M'\atop (\lambda,W)>0}
c_{\delta}
(\lambda^2/2)\lambda\cr
\rho_{z'}&
=\hbox{constant term of~}\bar\Theta_K(\tau)F_K(\tau)E_2(\tau)/24\cr
\rho_z&=-\rho_{z'}z'^2/2 +{1\over 2}\sum_{\lambda\in K'\atop (\lambda,W)>0}
\sum_{\delta\in M'/M\atop \delta|L=\lambda}
c_{\delta}(\lambda^2/2)
\B_2((\delta,z')).\cr
}
$$

Proof.
We have $h=h^-=m^+=m^-=j=k=0$
and $\B_2(x)=x^2-x+1/6$ for $0\le x\le 1$ so theorem 10.2 simplifies
to
$$
\eqalign{
&\sqrt 2 |z_{v^+}|\Phi_M(v,1,F_M)\cr
=&\Phi_K(\cdot,1,F_K)+\cr
&+
4\sum_{\lambda\in K'\atop (\lambda,W)>0}
\sum_{\delta\in M'/M\atop \delta|L=\lambda}
\left({\pi \over 2z_{v^+}^2}\right)^{-1}
c_{\delta}(\lambda^2/2)
\times\cr&\times
\pi^{2}
\Big(((\lambda,\mu)+(\delta,z'))^2-((\lambda,\mu)+(\delta,z'))+1/6\Big)
}
$$
where we take $\mu$ so that $(\lambda,\mu)$ is small and positive,
and we always take the value of $(\delta, z')$ with $0\le (\delta,z')<1$.

Hence for all $m>0$ and all $\mu\in K$ we have (using $|z_{v^+}|=m$)
$$
\eqalign{
&m(\mu,\rho)+m\rho_z+\rho_{z'}/2m+m\rho_{z'}z'^2/2-m\rho_{z'}\mu^2/2\cr
=&((m\mu,m,1/2m-mz'^2/2-m\mu^2/2),(\rho,\rho_{z'},\rho_z))\cr
=&{1\over 8\pi \sqrt 2} \Phi_M((m\mu,m,1/2m-mz'^2/2-m\mu^2/2),1,F_M)\cr
=&{1\over 16 \pi m}\Phi_K(\cdot,1,F_K)+\cr
&+
{1\over 16\pi}\sum_{\lambda\in K'\atop(\lambda,W)>0}
\sum_{\delta\in M'/M\atop \delta|L=\lambda}
c_{\delta}(\lambda^2/2)
\times\cr&\times
8\pi
\Big(((\lambda,m\mu)+m(\delta,z'))^2/m-
((\lambda,m\mu)+m(\delta,z'))+m/6\Big).\cr
}
$$
Now we compare coefficients of both sides, considered as functions
of $m$ and $\mu$. By comparing the coefficients of $1/m$
and using 9.3 we see that
$$\rho_{z'}= {1\over 8\pi}\Phi_K(\cdot,1,F_K)
=\hbox{constant term of~}\bar\Theta_K(\tau)F_K(\tau)E_2(\tau)/24.$$
By comparing the terms linear in $m\mu$ we see that
$$\eqalign{
\rho=
&{1\over 2}\sum_{\lambda\in K'\atop (\lambda,W)>0}
\sum_{\delta\in M'/M\atop \delta|L=\lambda}
c_{\delta}
(\lambda^2/2)(2(\delta,z')-1)\lambda\cr
=&{-1\over 2}\sum_{\lambda\in K', (\lambda,0,0)\in M'\atop(\lambda,W)>0}
c_{\delta}
(\lambda^2/2)\lambda\cr
}$$
because the terms with nonzero values of $(\delta, z')$ cancel out
in pairs.
By comparing coefficients of $m$ we see that
$$\eqalign{
&\rho_z+\rho_{z'}z'^2/2
={1\over 2}\sum_{\lambda\in K'\atop(\lambda,W)>0}
\sum_{\delta\in M'/M\atop \delta|L=\lambda}
c_{\delta}(\lambda^2/2)
\Big((\delta,z')^2-(\delta,z')+1/6\Big).\cr
}$$

This proves theorem 10.4.

We can get an extra
identity as follows  by comparing the terms quadratic in $\mu$.

\proclaim Theorem  10.5. Suppose $K$ is a negative definite even lattice
of dimension $b^-$
and $F_K$  is an automorphic form of weight $(-b^-/2,0)$
and type $\rho_K$
which is holomorphic on the upper half plane and meromorphic at
the cusps with Fourier expansion
$$
F_M(\tau)=
\sum_{\lambda\in K'}c_\lambda(\lambda^2/2)\e(\tau\lambda^2/2)e_\lambda
$$
with integral Fourier coefficients.
Then the numbers $c_\lambda(\lambda^2/2)$ are the coefficients of
a vector system (see below) of index
the constant term of $\bar\Theta_K(\tau)F_K(\tau)E_2(\tau)/24$.

Proof. We recall from [B95] that a vector system
is a set of numbers $c(\lambda)$ for $\lambda\in K'$ which
are zero for all but a finite number of $\lambda$ such that
\item{1.} $c(\lambda)=c(-\lambda)$.
\item {2.} $\sum_{\lambda\in K'}c(\lambda)(\lambda,\mu)^2
= -2m\mu^2$ for some constant $m$ (called the index of the vector system).

(We drop the condition in [B95] that the $c(\lambda)$'s should be nonnegative
from the definition of a vector system, and we get an extra sign
because $K$ is negative definite rather than positive definite.)
The first condition follows from the fact that
$c_{-\lambda}(n)=c_{\lambda}(n)$. The second condition follows because
if we compare the terms in the identity above that are quadratic in
$\mu$ and linear in $m$ we find that
$$\rho_{z'}\mu^2=
-\sum_{\lambda\in K'\atop (\lambda,W)>0}
\sum_{\delta\in M'/M\atop \delta|L=\lambda}
c_{\delta}(\lambda^2/2)
(\lambda,\mu)^2.
$$
(This extra identity is in some sense equivalent to the part of
theorem 10.3 that says that the function $\Phi_M(v,1,F_M)$ is a
polynomial rather than a rational function.)  This proves theorem
10.5.

By theorem 6.5 of [B95] the coefficients $c_\lambda(\lambda^2/2)$
can be used to define an almost holomorphic Jacobi form
of index $m$
by an infinite product.

It is also useful to know the inner product $(\rho(M,W,F),\lambda)$
for positive norm vectors $\lambda\in W$. These inner products can be
evaluated in terms of theta functions as follows.
\proclaim Theorem 10.6. Suppose $\lambda$ is a primitive
norm $2N>0$ vector in the closure of a Weyl chamber of the Lorentzian
lattice $M$, and suppose that $F$ is a modular form of type $\rho_M$
which is holomorphic on $H$. Then the inner product
$$(\lambda, \rho(M,W,F))$$
is the constant term of
$$-F(\tau)G_N(\tau)\Theta_{M,\lambda}(\tau)$$
where $\Theta_{M,\lambda}$ is the modular form of type
$\rho_M\otimes\bar\rho_{2N}$ given by the theta functions
of sublattices of $M$ whose vectors have given inner product
with $\lambda$.

Proof. We have $|\lambda|=2N$, so by definition of the Weyl vector,
the inner product $(\lambda, \rho(M,W,F))$ is equal to the
regularized integral
$${\sqrt {2N}\over 8\sqrt 2 \pi}\int
\bar\Theta_M(\tau,\lambda/|\lambda|,1)F(\tau)y^{1/2}dxdy/y^2.$$ We
substitute in
$$\Theta_M(\tau,\lambda/|\lambda|,1)= \Theta_{M,\lambda}
(\tau)\Theta_{2N}(\tau)$$
and use corollary 9.6 to evaluate the integral
as a constant term. This proves theorem 10.6.

{\bf Example 10.7.} Take $M=II_{1,25}$, $F(\tau)
=1/\Delta(\tau)=q^{-1}+24+O(q)$,
$N=1$ so that $\lambda$ has norm 2.
Recall that  $G_1(\tau)= e_0(-1/12 + q/2+O(q^2))+ e_1(q^{3/4}/3+O(q^{7/4}))$.
Let $c(n)$ be the number of vectors of norm $-n$ in
the dual of $\lambda^\perp\cap II_{1,25}$.
Then we see that
$(\rho(M,W,F), \lambda)$ is the constant term of
$$-(q^{-1}+24+\cdots)(-1/12 +(1/3)q^{3/4}+(1/2)q+\cdots)
(1+c(1/2)q^{1/4}+c(2)q+\cdots)$$
which is $-c(1/2)/3 +3/2 +c(2)/12$.

\proclaim 11.~Congruences for positive definite lattices.

Theorem 12.1 of [B95] states that if $M$ is a nonzero positive
definite even unimodular lattice then the constant term of
$\Theta_M(\tau)/\eta(\tau)^{\dim(M)}$ is divisible by 24,
which is a generalization of the well known fact that the number of roots
of a Niemeier lattice is divisible by 24. In this
section we find analogues of this congruence for lattices of higher
determinant, and also find a few other congruences. The main idea for
finding these congruences is as follows.  By the results of the
previous section we can often write the inner product $(\rho(M,W,F),\lambda)$
as a linear combination of coefficients of some theta function with
rational coefficients. On the other hand the Weyl vector $\rho(M,W,F)$ tends
to have integrality properties; for example, it often belongs to
$M'$. Therefore its inner product with $\lambda$ is often integral,
and this gives us congruences involving the coefficients of theta
functions.

\proclaim Lemma 11.1. Suppose that $K$ is a negative  definite even
lattice and let $(N)$ be the ideal of $\Z$
generated by the inner products of vectors of $K$.
Then the vector $(0,0,N)$ of $M=K\oplus II_{1,1}$ is in the
lattice generated by the norm $-2$ vectors of $M$ that
are not in $2M'$.

Proof. If $v,w\in K$ then
$$
(0,0,(v,w))
=-(0,1,-1)+(v,1,-1-v^2/2)+(w,1,-1-w^2/2)-(v+w,1,-1-(v+w)^2/2)
$$
is in the lattice generated by norm $-2$ vectors of $M$  that
are not in $2M'$.  This proves lemma 11.1 because by assumption $N$ is a
linear combination of the numbers $(v,w)$.

\proclaim Theorem 11.2. Suppose that $K$ is a negative  definite even
lattice of dimension $b^-$ and let $(N)$ be the ideal of $\Z$
generated by the inner products of vectors of $K$.
Suppose that $F$ is a modular form of weight $(-b^-/2,0)$ and type
$\rho_K$ which is holomorphic on $H$ and meromorphic at the cusps
and all of whose Fourier coefficients $c_{\gamma}(m,0)$ for $m< 0$
are integral.
Then the constant term of $NF(\tau)\bar\Theta_K(\tau)$ is divisible by 24.
(If $K$ is positive definite then we get the same result by changing
the sign of $K$, except that now $F$ has type $\bar\rho_K$ and
the constant term of $NF(\tau)\Theta_K(\tau)$ is divisible by 24.)

Proof.
We let $M$ be the Lorentzian lattice $K\oplus II_{1,1}$,
with a norm 0 vector $z=(0,0,1)$. By the wall crossing formula 6.4
any two Weyl vectors differ by an element of $M'$,
and in particular if $\sigma\in \Aut(M,F,C)$ then
$\sigma(\rho)-\rho$ lies in $M'$ because $\sigma(\rho)$ is also a Weyl vector.
The Weyl vector  $\rho(M,W,F)$ has integral inner product
with every norm $-2$ vector $r$ of $M$ not in $2M'$
because reflection $\sigma$ in the hyperplane $r^\perp$
lies in $\Aut(M,F,C)$ and hence takes $\rho$ to $\sigma(\rho)=
\rho+(\rho,r)r$ so $(\rho,r)r\in M'$. As $r$ is a primitive
vector of $ M'$ this
implies that $(\rho,r)\in \Z$.  By lemma 11.1 this implies that $\rho$
has integral inner product with $(0,0,N)$. By theorem 10.4 this inner
product is the constant term of
$NF(\tau)\bar\Theta_K(\tau)E_2(\tau)/24$.  As $E_2(\tau)\equiv 1\bmod
24$ the constant term of $N\bar\Theta_K(\tau)F(\tau)$ is divisible by
24. This proves theorem 11.2.

{\bf Example 11.3.} If we let $K$ be a nonzero even positive definite
unimodular lattice of
dimension $24n$ and take $F$ to be $\Delta(\tau)^{-n}$ we recover the
result that the constant term of $\Theta_K(\tau)/\Delta(\tau)^n$ is
divisible by 24.

{\bf Example 11.4.} The Weyl vector of the first
infinite product of example 13.7  lies in
$M'$ but not in $M$.

{\bf Example 11.5.} If $M=K\oplus II_{1,1}$ then
we see from the explicit formula for $\rho(M,W,F)$
that it is always true that
$24\rho(M,W,F)\in M'$. If we take $M=II_{1,1}$ and $F=1$ then the
smallest multiple of the Weyl vector
 in $M=M'$ is $24\rho(M,W,F)$, so in this case there is no better
congruence.

{\bf Example 11.6.} Take $K$ to be a one dimensional lattice generated by
an element of norm $-2$.
If $r$ is any norm 2 vector
in an even unimodular positive definite lattice then the theta function
of $r^\perp$ is a modular form of type $\rho_K$, whose
coefficients are all even except for the constant term. If
we apply this to the Niemeier lattices with root systems
$A_1^{24}$ and $A_2^{12}$ we get modular forms
$e_0(1+46q+O(q^2))+e_1(O(q^{7/4}))$ and
$e_0(1+66q+O(q^2))+e_1(2q^{3/4}+O(q^{7/4}))$
whose difference $e_0(20q+O(q^2))+e_1(2q^{3/4}+O(q^{7/4}))$
has even coefficients. Dividing by $2\Delta(\tau)$ we get
a modular form $F(\tau)=e_0(10+O(q))+e_1(q^{-1/4}+O(q^{3/4}))$
of weight $-1/2$, type $\rho_K$, and integer coefficients.
Then  the constant term of $\bar\Theta_KF(\tau)$ is 12,
so in this case we cannot omit the factor of $N=2$ in the congruence
of theorem 11.2 for the  lattice $K$.
This also gives an example where the Weyl vector does not lie
in $M'$ even though $F$ has integral coefficients.

The congruences above depend on looking at the inner product
of  a Weyl vector with a norm 0 vector. We can also get
congruences by looking at the inner product of a Weyl vector
with a positive norm vector and using theorem 10.6.

{\bf Example 11.7.} Let $\lambda$ be a norm $2$ vector in $II_{1,25}$ as in
example 10.7. We saw there that $(\rho(M,W,F),\lambda)=-c(1/2)/3+3/2+c(2)/12$.
As $\rho(M,W,F)\in II_{1,25}$ this number must be an integer.
Any even 25 dimensional unimodular positive definite lattice is of the form
$\lambda^\perp\cap II_{25,1}$ for some norm $-2$
vector $\lambda\in II_{25,1}$, and its dual has
norm $1/2$ vectors if and only if
it is not the sum of a Niemeier lattice and a one dimensional  lattice.
Hence we see that if $K$ is any even positive definite lattice which is not
the sum of a Niemeier lattice and a one dimensional lattice then the
number $c(2)$ of norm 2 vectors is congruent to 6 mod 12.

{\bf Example 11.8.} More generally, suppose that $K$
of dimension greater than 1 is the orthogonal
complement in $II_{1,1+8n}$ of a primitive norm $2N>0$ vector
$\lambda\in II_{1,1+8n}$, and let $F$ be any (complex valued) modular
form of level 1 and weight $-4n<0$ which is holomorphic on $H$,
meromorphic at the cusps, and has integer coefficients. Then the
constant term of
$$\bar\Theta_K G_N F$$
is an integer. This follows easily from theorem 10.6, with
$M=K\oplus II_{1,1}$, because the Weyl vector $\rho(M,W,F)$ lies
in $M=M'$ by the argument used in 11.2.

\proclaim 12.~Hyperbolic reflection groups.

In this section we give a sufficient condition for a Lorentzian
lattice to have a reflection group of finite index in its automorphism
group, or in other words for the reflection group to be arithmetic. We
give some examples to show that this seems to account for most (and
possibly all) of the known examples of Lorentzian lattices with this
property.

\proclaim Theorem 12.1. Suppose $M$ is a Lorentzian lattice of dimension
$1+b^-$. Suppose that $F$ is a modular form of type $\rho_M$
and weight $(1/2-b^-/2,0)$ which is holomorphic on $H$ and meromorphic
at cusps and all of whose Fourier coefficients $c_\lambda(m)$
are real for $m<0$.
Finally suppose that if $c_\lambda(\lambda^2/2)\ne 0$ and
$\lambda^2<0$ then reflection in $\lambda^\perp$ is in $\Aut(M,F,C)$. Then
$\Aut(M,F,C)$ is the semidirect product of a reflection subgroup
and a subgroup fixing the Weyl vector $\rho(M,W,F)$ of a Weyl chamber $W$.
In particular if the Weyl vector has positive norm then
the reflection group of $M$ has finite index in the automorphism  group
and has only a finite number of simple roots.  If the Weyl
vector has zero norm but is nonzero then the quotient of
the automorphism  group of $M$ by the reflection subgroup
has a free abelian subgroup of finite index.
(Warning: If we want the Weyl vector to be  a Weyl vector
in the usual sense of
reflection groups, so that $(\rho(M,W,F),r)=r^2/2$ for all simple roots $r$,
we may need to multiply the simple roots $r$ by nontrivial factors,
so that they are not necessarily primitive vectors of $M$.)

Proof. By assumption, reflection in any wall of a Weyl chamber is in
$\Aut(M,F,C)$ and takes any Weyl chamber to the Weyl chamber on the other side
of the wall.  Therefore the group generated by these reflections is a
hyperbolic reflection group acting transitively on the Weyl chambers of
$F$.  Then the group $\Aut(M,F,C)$ is the
semidirect product of this reflection subgroup by the subgroup fixing a
Weyl chamber (which also fixes the corresponding Weyl vector).

If the Weyl vector has positive norm, then the subgroup fixing it
has finite order. As $\Aut(M,F,C)$ has finite index in $\Aut(M)$,
it then follows that the reflection  group of $M$ has finite index in its
automorphism  group.
Similarly if the Weyl vector has zero norm but is nonzero then the
group fixing it has a free abelian subgroup of finite index.
This proves theorem 12.1.

{\bf Example 12.2.} Take $M$ to be the lattice $II_{1,9}$, $II_{1,17}$, or
$II_{1,25}$ and take $F(\tau)$ to be $E_4(\tau)^2/\Delta(\tau)$,
$E_4(\tau)\Delta(\tau)$, or $1/\Delta(\tau)$. Then the conditions of
theorem 12.1 are satisfied, and the Weyl vector has norm 1240, 620, or
0 in the 3 cases. Hence the first 2 lattices have reflection groups of
finite index in their automorphism groups, and $II_{1,25}$ has a norm
0 Weyl vector for its reflection group. This was first proved by
Conway [C] using the fact that the Leech lattice has covering radius
$\sqrt 2$.  Conway's proof can be run in reverse to deduce that the
Leech lattice has covering radius $\sqrt 2$ from the fact that the
reflection group of $II_{1,25}$ has a norm 0 Weyl vector.

{\bf Example 12.3.} Take $M$ to be the even sublattice of $I_{1,25-n}$.
We take $F$ to be $\Theta_{D_n}(\tau)/\Delta(\tau)$. Then $F$ satisfies
the condition of theorem 12.1 if $n=4$ or $n\ge 6$. In particular we recover
Vinberg and Kaplinskaja's result [V-K] that the reflection
group of $I_{1,25-n}$
has finite index for $n\ge 6$. (For $n=4$ this does not show that
the reflection group of $I_{1,21}$ has finite index in its automorphism
group
as the automorphism  group of $I_{1,21}$ has index 3 in that of
its even sublattice and in particular not all reflections of the
even sublattice of $I_{1,21}$ are reflections of $I_{1,21}$.)

{\bf Example 12.4.} Take $M$ to be $BW\oplus II_{1,1}$ where $BW$ is the
16 dimensional Barnes-Wall lattice, and take $F$ to be
$\Theta_{E_8(2)}(\tau)/\Delta(\tau)$. Then the Weyl vector has norm 0,
so this case is similar to that of $II_{1,25}$: the reflection  group
has a norm 0 Weyl vector.

{\bf Example 12.5.} Take $M$ to be the even Lorentzian lattice of
dimension 20 and determinant 3, and take $F$ to be
$\Theta_{E_6}/\Delta$. Then we see that the reflection group of $M$
has finite index in the automorphism group.

{\bf Example 12.6.} More generally suppose $M$ is any primitive sublattice
of $II_{1,25}$ with the property that any negative norm vector of $M'$
 that is the projection of a norm $-2$ vector of $II_{1,25}$
is a root. Then if we take $F$ to be the theta function
of $M^\perp\cap II_{1,25}$ divided by $\Delta(\tau)$ it
satisfies the conditions of theorem 12.1, so the reflection
group of $M$ has a Weyl vector, and is arithmetic
if this Weyl vector has positive norm.

\proclaim 13.~Holomorphic infinite products.

Suppose that $b^+=2$, $p=1$ (so $m^+=m^-=0$), and $yF(\tau)$ is holomorphic,
so that $c_\gamma(n,k)=0$ for $k\ne 0$. Then we see from theorem 6.2
that the function $\Phi_M(v,1,F)$ has logarithmic singularities.
This suggests that we should try to exponentiate $\Phi_M$. It turns out
that $\Phi_M=\log|\Psi_M|$ for a certain multi-valued automorphic form,
which is holomorphic if the numbers $c_\gamma(n)$ satisfy certain
positivity and integrality conditions. Moreover we can write $\Phi_M$ as
an explicit infinite product for each cusp corresponding to a norm 0
vector of $M$, and can explicitly describe all the zeros of
$\Psi_M$. The special case when $M$ is unimodular is the main
result of [B95].

We recall some facts about hermitian symmetric spaces and set up some
notation.  We let $M$ be any even lattice of signature $(2,b^-)$. We
choose a continuously varying orientation on the 2 dimensional
positive definite subspaces of $M\otimes \R$; there are 2 ways to do
this. We put a complex structure on the Grassmannian of $M$ as
follows.  If $X_M$ and $Y_M$ are an oriented orthogonal base of some
element $v$ of $G(M)$ then we map $v$ to the point
$Z_M=X_M+iY_M\in M\otimes \C$
representing a point of the
complex projective space $\P(M\otimes \C)$. This identifies $G(M)$
with an open subset of $\P(M\otimes \C)$ in a canonical way, and gives
$G(M)$ a complex structure invariant under the subgroup $O_M(\R)^+$ of
index 2 of $O_M(\R)$ of elements preserving the orientation on the 2
dimensional positive definite subspaces, or equivalently of elements
whose spinor norm has the same sign as the determinant. There is a
principal $\C^*$ bundle $P$ over this hermitian symmetric space,
consisting of the norm 0 points $Z_M=X_M+iY_M\in M\otimes \C$ such that $X_M$
and $Y_M$ form an oriented base of an element of $G(M)$. The fact that
$Z_M=X_M+iY_M$ has norm 0 is equivalent to saying that $X_M$ and $Y_M$ are
orthogonal and have the same norm.  We define an automorphic form of
weight $k$ on $G(M)$ to be a function $\Psi_M$ on $P$ which is homogeneous of
degree $-k$ and invariant under some subgroup $\Gamma$ of finite index of
$\Aut(M)^+$.  More generally if $\chi$ is a one dimensional representation
of $\Gamma$ then we say $\Psi $ is  an automorphic form
of character $\chi$ if
$\Psi_M(\sigma(Z_M))= \chi(\sigma)\Psi(Z_M)$ for $\sigma\in \Gamma$.
(Automorphic forms of weight $k$ can also be thought of
as sections of the line bundle corresponding to the principal $\C^*$
bundle $P$ and the representation $z\rightarrow z^{-k}$ of the
structure group $\C^*$.)

Now suppose that we have selected a norm 0 vector $z\in M$ and a
vector $z'\in M'$ with $(z,z')=1$. We let $K$ be the lattice $(M\cap
z^\perp)/\Z z$, and we identify $K\otimes \R$ with the subspace
$M\otimes \R \cap z^\perp\cap z'^\perp$ of $M\otimes \R$ (which we can
do as each element of this subspace represents an element of $K\otimes
\R$). This identifies $K$ with a subgroup of $M\otimes \R$, but this
is not in general a subgroup of $M$. The lattice $K$ is Lorentzian so
$K\otimes \R$ has 2 components of positive norm vectors.  We let $C$
be the open positive cone in $K\otimes \R$, determined as follows: if
$X_M$, $Y_M$ is an oriented basis of some element of $G(M)$ with
$(Y_M,z)=0$, $(X_M,z)>0$, then $Y_M$ represents a positive norm vector of
$K\otimes \R$ whose component only depends on $z$ and the choice or
orientation. This open cone $C$ is called the positive cone.

If $\lambda\in K\otimes \C$, $m,n\in \C$, we write $(\lambda,m,n)$ for
the point $\lambda+mz'+nz\in M\otimes \C$, so that $z=(0,0,1)$,
$z'=(0,1,0)$, and $(\lambda,m,n)^2=\lambda^2+2mn+m^2z'^2$.  We can
embed the set of points $Z=X+iY$ of $K\otimes \C$ with imaginary part
$Y$ in $C$ into $P$ by mapping $Z$ to the unique norm 0 point
$Z_M=(Z,1,-Z^2/2-z'^2/2)$ having inner product 1 with $z$ and
projection $Z$ in $K$. If we compose this map with the natural
projection from $P$ to $G(M)$ we get an isomorphism from $K\otimes
\R+iC$ to $G(M)$, and in particular any automorphic form $\Psi_M$ is
determined by its restriction $\Psi_z$ to the image of $K\otimes
\R+iC$ in $P$.  More explicitly, if $\Psi_M$ is an automorphic form of
weight $k$ then we define $\Psi_z(Z)$ for $Z\in K\otimes\R+iC$ by
$$\Psi_z(Z)=\Psi_M(Z_M)=\Psi_M((Z,1,-Z^2/2-z'^2/2)).$$
If $\Psi_z$ is a function given by the restriction
of an automorphic form of weight $k$ as above then we will also call $\Psi_z$
an automorphic form of weight $k$.
(Warning: the action of $\Aut(M,F)$ on $K\otimes\R+iC$ is not the restriction
of the action on $P$.)

If $Z=X+iY\in K\otimes \R+iC$ represents a point $v$ of $G(M)$ as above, then
(with the notation of section 5)
$$\eqalign{
X_M&=(X,1,Y^2/2-X^2/2-z'^2/2)\cr
Y_M&=(Y,0,-(X,Y))\cr
z_{v^+}&= (z_{v^+},X_M)X_M/X_M^2+(z_{v^+},Y_M)Y_M/Y_M^2
={(X,1,Y^2/2-X^2/2-z'^2/2)\over Y^2}\cr
z_{v^+}^2&= 1/Y^2\cr
\mu &= X\in K\otimes \R \cr
w^+& \hbox{is spanned by } Y\in K\otimes \R\cr
\lambda_{w^+}&= (\lambda,Y)Y/Y^2 \quad\hbox{(for $\lambda\in K'$)}\cr
|\lambda_{w^+}|
&= |(\lambda,Y)||z_{v^+}|.\cr
}
$$

\proclaim Lemma 13.1. Suppose that $\Psi_z$ is a holomorphic function
on $K\otimes \R+iC$ such that $A(\log|\Psi_z(Z)| +k \log|Y|+B)$ is
the restriction of a function $\Phi_M$ on $P$ that is homogeneous of
degree 0 and invariant
under a subgroup $\Gamma$ of finite index in $\Aut(M)^+$
for some integer $k$ and real numbers $A\ne 0$ and $B$.
 Then
$\Psi_z(Z)=\Psi_M((Z,1,-Z^2/2-z'^2/2))$ for an
automorphic form $\Psi_M$
which is holomorphic on $P$ and of weight $k$ for some
one dimensional unitary representation $\chi$ of $\Gamma$.

Proof. Extend $\Psi_z$ to a holomorphic function $\Psi_M$ on $P$ that
is homogeneous of degree $-k$. Then the function
$|\Psi(X_M+iY_M)||Y_M|^k$ is homogeneous of degree 0 and so is
invariant under the action of $\Gamma$ on $P$ by the assumption in the
lemma. But then the function
$|\Psi_M(X_M+iY_M)/\Psi_M(\sigma(X_M+iY_M))| =
(|Y_M|/|\sigma(Y_M)|)^{-k}=1$ is constant for any fixed $\sigma\in
\Gamma$, so $\Psi_M(\sigma(Z_M))=\chi(\sigma)\Psi_M(Z_M)$ for some
constant $\chi(\sigma)$ of absolute value 1.  It is obvious that
$\chi(\sigma_1\sigma_2)=\chi(\sigma_1)\chi(\sigma_2)$, so $\chi$ is a
one dimensional unitary representation of $\Gamma$.  Therefore
$\Psi_z$ is the restriction of a holomorphic function $\Psi_M$ on $P$
of degree $-k$ which transforms under $\Gamma$ according to $\chi$.
This proves lemma 13.1.

\proclaim Lemma 13.2. The constant term at $s=0$ of
$$\Gamma(s+1/2){\pi}^{-s-1/2} (2z_{v^+}^2)^{s}
\sum_{n>0}{\e(n\delta/N)\over n^{2s+1}}
$$
for $N\in \Z$, $N>0$, $\delta\in \Z/N\Z$ is
$-\log(1-\e(\delta/N))$ if $\delta\ne 0$ and
$\log|z_{v^+}| -\Gamma'(1)/2-\log(\sqrt{2\pi})$ if $\delta=0$.

Proof. If $\delta\ne 0$ the function is holomorphic at $s=0$
and the series is just a series for a logarithm. If $\delta=0$
then we can work out the constant term at
$s=0$ by multiplying together the following series:
$$\eqalign{
\pi^{-1/2}\Gamma(s+1/2)&= 1+s(\Gamma'(1)-2\log(2))+O(s^2)\cr
\pi^{-s}&= 1 -s\log(\pi)+O(s^2)\cr
(2z_{v^+}^2)^s &= 1+s\log(2z_{v^+}^2)+O(s^2)\cr
\sum_{n>0}{1\over n^{1+2s}}&=\zeta(1+2s)= {1\over 2s} -\Gamma'(1)+O(s)\cr
}
$$
where $\Gamma'(1)=-.57721=\lim_{n\rightarrow\infty}
(\log(n)-1/1-1/2-\cdots-1/n)$ is Euler's constant.
This proves lemma 13.2.

\proclaim Theorem 13.3. Suppose $M$ is an even lattice of signature
$(2,b^-)$ and $F$ is a modular form of weight $1-b^-/2$ and
representation $\rho_M$ which is holomorphic on $H$ and meromorphic at cusps
and whose
coefficients $c_\lambda(m)$ are integers for $m\le 0$.
Then there is a meromorphic  function $\Psi_M(Z_M,F)$
for $Z\in P$  with the following properties.
\item{1.} $\Psi_M(Z_M,F)$ is an automorphic form of weight $c_0(0)/2$
for the group $\Aut(M,F)$ with respect to some unitary character $\chi$
of $\Aut(M,F)$.
\item{2.} The only zeros or poles  of $\Psi_M$ lie on the
rational quadratic divisors $\lambda^\perp$
for $\lambda\in M$, $\lambda^2<0$ and are zeros of order
$$\sum_{0<x\in \R\atop x\lambda\in M'}c_{x\lambda}(x^2\lambda^2/2).
$$
(or poles if this number is negative).
\item{3.}
$$\log|\Psi_M(Z_M,F)|={-\Phi_M(Z_M, 1,F)\over4}
-{c_0(0)\over 2}(\log|Y_M|+\Gamma'(1)/2+\log\sqrt{2\pi})$$
\item{4.} $\Psi_M$ is a holomorphic function
if the orders of all zeros in item 2 above are nonnegative.
If in addition $M$ has
dimension at least 5, or if $M$ has dimension 4 and contains no 2 dimensional
isotropic sublattice, then $\Psi_M$ is a holomorphic automorphic form.
If in addition $c_0(0)=b^--2$ then the only nonzero Fourier
coefficients of $\Psi_M$ correspond to vectors of $K$ of norm 0.
\item{5.} For each primitive norm 0 vector $z$ of $M$
and for each Weyl chamber $W$ of $K$ the restriction
$\Psi_z(Z,F)$ has an infinite product expansion converging
when $Z$ is in a neighborhood of the cusp of $z$
and $Y\in W$ which is some
constant of absolute value
$$\prod_{\delta\in \Z/N\Z\atop \delta\ne 0}
(1-\e(\delta/N))^{c_{\delta z/N}(0)/2}
$$
times
$$\e((Z,\rho(K,W,F_K)))
\prod_{\lambda\in K'\atop (\lambda,W)>0}
\prod_{\delta\in M'/M\atop \delta|L=\lambda}
(1-\e((\lambda,Z)+(\delta,z')))^{c_{\delta}(\lambda^2/2)}.
$$
(The vector $\rho(K,W,F_K)$ is the Weyl vector, which can usually
be evaluated explicitly using the theorems in section 10.)

Proof.
We first assume that $M$ has a primitive norm 0 vector $z$.
We pick a vector $z'\in M'$ with
$(z',z)=1$.
We use theorem 7.1  to see that $\Phi_M(v,1,F)$
is the constant term at $s=0$ of the following expression
(note that $m^+=m^-=h^+=h^-=h=j=k=0$).
$$\eqalign{
&\Phi_M(v,1,F)\cr
=&{1\over \sqrt 2 |z_{v^+}|}\Phi_K(w,1,F_K)+\cr
&+{  \sqrt 2\over|z_{v^+}|}\sum_{n> 0}
\sum_{\lambda\in K'}
\e((n\lambda,\mu))
\sum_{\delta\in M'/M\atop \delta|L=\lambda}
\e(n(\delta,z'))
\times\cr&\times
\int_{y>0}c_{\delta}(\lambda^2/2)
\exp(-\pi n^2/ 2yz_{v^+}^2-2\pi y\lambda_{w^+}^2)  y^{1-s-5/2} dy
\cr
}
$$
For $v$ given by a norm 0 vector $(X+iY,1,-(X+iY)^2/2-z'^2/2)\in M\otimes\C$
we use lemmas 7.2 and 7.3 to see that
this is the constant term at $s=0$ of
$$\eqalign{
&{8\pi\over  |z_{v^+}|}(Y/|Y|,\rho(K,W,F_K))+\cr
+&{\sqrt 2\over  |z_{v^+}|}
\sum_{n>0}\sum_{\delta\in \Z/N\Z}\e(n\delta/N)
\left({\pi n^2\over 2z_{v^+}^2}\right)^{-s-1/2}
c_{\delta z/N}(0)
\Gamma(s+1/2)
+\cr+&
{\sqrt 2\over  |z_{v^+}|}
\sum_{n>0}
\sum_{\lambda\in K'\atop \lambda\ne 0}
\e((n\lambda,\mu))
\sum_{\delta\in M'/M\atop \delta|L=\lambda}
\e(n(\delta,z'))
{\sqrt 2|z_{v^+}| \over n }
\times\cr&\times
\exp(-2\pi n|\lambda_{w^+}|/|z_{v^+}|)
c_{\delta}(\lambda^2/2)
\cr
=&8\pi (Y,\rho(K,W,F_K))+\cr
+&
2\sum_{\delta\in \Z/N\Z}c_{\delta z/N}(0)
\Gamma(s+1/2){\pi}^{-s-1/2} (2z_{v^+}^2)^{s}
\sum_{n>0}{\e(n\delta/N)\over n^{2s+1}}
+\cr+&
2\sum_{\lambda\in K'\atop \lambda\ne0}
\sum_{n>0}
\e((n\lambda,X))
\sum_{\delta\in M'/M\atop \delta|L=\lambda}
\e(n(\delta,z'))\exp(-2\pi n|(\lambda,Y)|)
{1\over  n }
c_{\delta}(\lambda^2/2)
.\cr
}
$$
If we apply lemma 13.2 we see that this is equal to
$$
\eqalign{
&8\pi (Y,\rho(K,W,F_K))+\cr
+&
c_{0}(0)(\log(z_{v^+}^2)-\Gamma'(1)-\log(2\pi))+
2\sum_{\delta\in \Z/N\Z\atop \delta\ne 0}c_{\delta z/N}(0)
(-\log(1-\e(\delta/N)))+
\cr
+&
4\sum_{\lambda\in K'\atop (\lambda,W)>0}
\sum_{\delta\in M'/M\atop \delta|L=\lambda}
-c_{\delta}(\lambda^2/2)
\log\left|1-\e\left((\lambda,X)+(\delta,z')+{i|(\lambda,Y)|}\right)\right|
.\cr
}
$$
Comparing this with $\Psi_z$ defined by the infinite product for  in 13.3
we see that
$$
\Phi_M(v,1,F)= -4\log|\Psi_z(Z,F)|-4{c_0(0)\over 2}
(\log|Y|+\Gamma'(1)/2+\log\sqrt{2\pi}).
$$
This means that $\Psi_z$ satisfies the conditions of lemma 13.1,
which shows that the function $\Psi_M$ defined in 13.1 and 13.3
is an automorphic form of weight $k=c_0(0)/2$, and also shows
that its restriction $\Psi_z(Z)$ has the infinite product expansion in 13.3.

This proves
theorem 13.3 in the case when $M$ has a primitive norm 0 vector $z$.
If $M$ has no primitive norm 0 vector (which can only happen if
$M$ has dimension
at most 4) we have to
show the existence of the holomorphic function
$\Psi_M$. This follows from the embedding trick (section 8)
by using theorem 13.3 applied to the larger lattice we embed $M$ in.
(There is no need to prove anything about Fourier expansions
because these do not exist when $M$ has no primitive norm 0 vectors!)

We can work out the zeros and poles of $\Psi_M$ using the
singularity theorem 6.2, because these are singularities  of $\Phi_M$.
More precisely we see that the singularities of
$-4\log|\Psi_M(Z_M,F)|$
are of the form
$$\sum_{\lambda\in M'\cap (Z,1,*)^\perp\atop \lambda\ne 0}
-c_\lambda(\lambda^2/2)\log(\lambda_{v^+}^2)
$$
and as
$$\log(\lambda_{v^+}^2) = \log((\lambda,X_M/|X_M|)^2+(\lambda,Y_M/|Y_M|)^2)
= 2\log|(\lambda,Z_M)| -2\log|Y_M|
$$
we see that $\Psi_M$ has zeros and poles as stated in theorem 13.3.

The statements in part 4 of 13.3 follow from the Koecher boundedness
principle and the theory of singular weights.
This proves theorem 13.3.

Remark. The character $\chi$ in 13.3 is often nontrivial. It can often
be worked out explicitly as follows. If $\sigma$ is the reflection
of some negative norm root of $M$ in $\Aut(M,F)^+$ then $\chi(\sigma)$
is 1 or $-1$ if $\Psi_M$ has a zero of even or odd order along
the divisor of points fixed by $\sigma$, and the order of this
zero can be worked out using theorem 13.3. Hence if the abelianization
of $\Aut(M,F)^+$ is generated by such reflections this completely
determines $\chi$. In this case $\chi$ has order 1 or 2, but
in general it can sometimes
have higher order; for example, when $F$ is the theta
function of a one dimensional lattice and $\Psi_z$ is the
square of the Dedekind $\eta$ function it has order 12.

\proclaim Corollary 13.4. Suppose $K$ is an even Lorentzian lattice
of signature $(1,b^--1)$
such that either $K$ has dimension at least 3, or $K$ is anisotropic
of dimension 2.
Suppose that  $F$ is a modular form of weight
$(b^-/2-1/2,-1/2)$ and type $\rho_K$
such  that all the coefficients $c_\delta(m)$ of $F$
are nonnegative integers for $m<0$.
Then the Weyl vector $\rho(K,W,F)$ of any
Weyl chamber lies in the closure of the positive cone of $K\otimes \R$.

Proof. We let $M$ be the lattice $K\oplus II_{1,1}$,
and consider the automorphic form $\Psi_z(Z,F)$. By theorem
13.3 this is a holomorphic  function transforming like an automorphic form.
By the assumption  on $K$ the Koecher boundedness principle applies
to $\Psi_z$, so $\Psi_z$ is holomorphic and therefore all its
nonzero Fourier coefficients correspond to vectors in the closure
of the positive cone of $K\otimes \R$. But the coefficient
of $\rho(K,W,F)$ is nonzero, so $\rho(K,W,F)$
lies in the closure of the positive
cone. This proves corollary 13.4

{\bf Example 13.5.} Take $K=II_{1,1}$ and let
$F(\tau)=j(\tau)-744=q^{-1}+196884q+\cdots$ be the elliptic modular
function minus 744. Then the Weyl vector does not lie in the closure
of the positive cone. If we try to apply the proof above, we find that
$\Psi_z$ is essentially $j(\sigma)-j(\tau)$ which is holomorphic at
finite points but has singularities at infinity, so the Koecher
boundedness principle does not hold for this function.  This shows
that the result of corollary 13.4 is not true for 2 dimensional
lattices which are not anisotropic.

{\bf Example 13.6.} The two functions in example 13.7  have Weyl vectors
which are 0 or are nonzero norm 0 vectors in the closure of the
positive cone, and there are examples with Weyl vectors in the interior
of the positive cone.

We now give some examples of the infinite
products constructed in theorem 13.3. For other examples
which follow from 13.3 see [B92], [B95], and [G-N].

{\bf Example 13.7.} We will recover the denominator
formula of the fake monster
superalgebra used in [B96] to show that the moduli space of Enriques
surfaces is quasiaffine. As a bonus we get a second denominator
formula for another rank 10 superalgebra with no real roots.

We take the lattice $M$ to be the maximal even sublattice of
$I_{2,10}$, so that $M'/M$ is isomorphic to $\Z/2\Z\times \Z/2\Z$.
We denote the elements of this group by $00$, $01$, $10$, and $11$,
with the element $11$ having norm $1/2\bmod \Z$ and the others having norm 0.
We define a modular form $F=\sum_\gamma e_\gamma f_{\gamma}$
of weight $(-4,0)$ and representation
$\rho_M$ by setting
$$\eqalign{
f_{00}(\tau)&= 8\eta(2\tau)^8/\eta(\tau)^{16} = 8+128q+1152q^2+\cdots\cr
f_{10}(\tau)=f_{01}(\tau) &=-8\eta(2\tau)^8/\eta(\tau)^{16}
= -8-128q-1152q^2-\cdots\cr
f_{11}(\tau) =&8\eta(2\tau)^8/\eta(\tau)^{16} +\eta(\tau/2)^8/\eta(\tau)^{16}
= q^{-1/2} +36q^{1/2}+402q^{3/2}+\cdots\cr
}
$$
(This behaves correctly under $S$ because
$\eta(-1/\tau)=\sqrt{\tau/i}\eta(\tau)$, and the only nontrivial thing
to check for the transformations under $T$ is that all integral powers
of $q$ of $f_{11}(\tau)$ vanish.)

By theorem 13.3 the function $\Psi_z(Z)=\Psi_z(Z,F)$ is a holomorphic
automorphic form of weight $c_{00}(0)/2=4$, whose only zeros are zeros
of order 1 at the divisors of norm -1 vectors of $M'$ (coming from the
coefficient of $q^{-1/2}$ in $f_{11}$).

We will also work out the infinite product of $\Psi_z$ at the cusps
corresponding to primitive norm 0 vectors of $M$. The lattice $M$ has
2 orbits of primitive norm 0 vectors under $\Aut(M,F)=\Aut(M)$, one of
level $N=1$ and one of level $N=2$. The function $\Psi_z$ has singular
weight $8/2=4$ so it is a singular automorphic form and therefore the
only nonzero Fourier coefficients correspond to norm 0 vectors. As the
Fourier coefficients of norm 0 vectors of the infinite products are
easy to work out this means we can also work out the Fourier
coefficients of $\Psi_z$ explicitly in the examples below.

Level 1 cusp: In this case we write $M=K\oplus II_{1,1}$, and take the
level 1 norm 0 vector $z$ to be a primitive norm 0 vector of the
$II_{1,1}$. The lattice $K$ is then the lattice of even norm vectors
of $I_{1,9}$ of determinant 4.  The Weyl vector of $K$ is a primitive
norm 0 vector of $K'$ which is half a characteristic vector of
$I_{1,9}$.  Define $c(n)$ by
$$\sum_nc(n)q^n=f_{00}(\tau)+f_{11}(\tau) = q^{-1/2}+8+36q^{1/2}+O(q).$$
The infinite product in theorem 13.3 is
$$\eqalign{
&\e((\rho,Z))\prod_{\lambda\in K'\atop(\lambda,W)>0}
(1-\e((Z,\lambda)))^{\pm c(\lambda^2/2)}\cr
=&
\sum_{w\in G}\det(w)\e((w(\rho),Z))\prod_n(1-\e(n(w(\rho),Z)))^{(-1)^n8}\cr
}
$$
where the sign in the exponent is 1 if $\lambda\in K$ or if $\lambda$
has odd norm, and $-1$ if $\lambda$ has even norm but is not in $K$.
The group $G$ is the reflection group generated by reflections
of the norm $-1$ vectors of $K$.
This is the denominator formula (see [R]) for the superalgebra of
superstrings on a 10 dimensional torus used in [B96]
(and  in [B92] page 415).
In [B96] it is shown that $\Psi_z$ can be considered as
an automorphic form on the period space of Enriques surfaces,
vanishing exactly on the singular Enriques surfaces.

Level 2 cusp: In this case we  decompose $M$ as $M=K\oplus II_{1,1}(2)$
with $K=II_{1,9}$ of determinant 1
(coming from the decomposition $I_{2,10}=II_{1,9}\oplus I_{1,1}$)
and take $z$ to be a primitive norm 0 vector of  $II_{1,1}(2)$.
The Weyl vector of $K$ is now 0.
The denominator formula in this case is
$$\prod_{\lambda\in K\atop (\lambda,W)>0} (1-\e((Z,\lambda)))^{c(\lambda^2/2)}
(1+\e((Z,\lambda)))^{-c(\lambda^2/2)}
=1+ \sum_{\lambda}a(\lambda)\e((Z,\lambda))
$$
where $a(\lambda)$ is the coefficient
of $q^n$ of
$$\eta(\tau)^{16}/\eta(2\tau)^8= 1 - 16q + 112q^2 - 448q^3 + O(q^4)$$
if $\lambda$ is $n$ times a primitive norm 0 vector in the closure of the
positive cone $C$, and 0 otherwise.

This is the denominator formula of another superalgebra of
superstrings on a 10 dimensional torus. (More precisely
it is the twisted denominator formula corresponding to the automorphism
which is 1 on the ordinary elements and $-1$ on the super elements;
the untwisted denominator formula is just $0=0$.) This algebra is
a generalized Kac-Moody superalgebra whose simple roots are
exactly the norm 0 vectors in the closure of the positive cone,
each of which has multiplicity 8 as both an ordinary root
and as a super root. The multiplicities of both the ordinary
and the super root spaces of any
other vector $\lambda\in K$ are both $c(\lambda^2/2)$.
This superalgebra has no real roots, or in other words no tachyons
(as expected for some superstrings), and for each vector $\lambda$
the ordinary and super root spaces have the same dimension.

The two superalgebras above look quite different at first sight:
they have different root lattices, different Weyl vectors,
one has trivial Weyl group and no real roots while the other has
an infinite number of real simple roots. However the denominator functions
of these two algebras are really the same function expanded about
2 different cusps.

{\bf Example 13.8.} We will temporarily abandon our convention
that modular forms are level 1 and vector valued.  Recall that in
[B92] there are several infinite products in 2 variables whose
exponents are coefficients of Hauptmoduls which turn out to be modular
functions of 2 variables. By lemma 2.6 we can construct
a vector valued modular form from any level $N$ modular form,
so by inserting this into theorem 13.3 we can find similar
infinite products corresponding to  arbitrary
modular functions.
We will write out the case when $f$ is a modular form for
$\Gamma_0(N)$ explicitly.
Suppose $f(\tau)=\sum_{n\in \Z}c(n)q^n$ is a complex valued modular function
for the group $\Gamma_0(N)$
with zero constant term $c(0)=0$.
We let $K$ be the lattice $II_{1,1}$, and let $M$ be the
lattice generated by the norm 0 vectors $z$ and $z'$
which are orthogonal to $K$ and have inner product $N$.
Using theorem 13.3 and lemma 2.6 we find that the infinite product
$$\e(\rho_\sigma\sigma+\rho_\tau\tau)\prod_{m>0\atop n\in \Z}
\prod_{d\in (\Z/N\Z)^*}(1-\e(d/N)\e(\sigma md)\e(\tau nd))^{c(mn)}$$
is a modular function of 2 variables
for some $\rho_\sigma$, $\rho_\tau$.

{\bf Example 13.9.} Suppose that $f(\tau)=\sum_nc(n)q^n$ is a modular function
for the normalizer  $\Gamma_0(2)+$ of $\Gamma_0(2)$
with vanishing constant term $c(0)=0$.
We let $K$ be the lattice generated by 2 norm 0 vectors having
inner product 2, and we define a modular function of type
$\rho_K$ by
$$\eqalign{
f_{00}(\tau)&= f(\tau)+f(\tau/2)/2+f((\tau+1)/2)/2\cr
f_{10}(\tau)=f_{01}(\tau) &=f(\tau/2)/2+f((\tau+1)/2)/2\cr
f_{11}(\tau) =&f(\tau/2)/2-f((\tau+1)/2)/2\cr
}
$$
Applying theorem 13.3 (with $M=K\oplus II_{1,1}$) we see that
$$
\e(\rho_\sigma\sigma+\rho_\tau\tau)
\prod_{m>0,n\in \Z}(1-\e(m\sigma)\e(n\tau))^{c(mn)}
(1-\e(2m\sigma)\e(2n\tau))^{c(2mn)}
$$
is a modular function of 2 variables (for some numbers $\rho_\sigma$
and $\rho_\tau$). In particular if $f$ is the Hauptmodul for
$\Gamma_0(2)+$ we recover the denominator function for the baby
monster Lie algebra ([B92, section 10]). There are similar results if
2 is replaced by any prime, or more generally by a square free
integer.

\proclaim 14.~The Shimura-Doi-Naganuma-Maass-Gritsenko-...~correspondence.

Shimura's correspondence [Sh] takes modular forms of half integral
weight $k+1/2$ to modular forms of integral weight $2k$, which should
be thought of as modular forms of weight $k$ for the group
$O_{2,1}(\R)$. Kohnen [Ko] modified Shimura's correspondence to go from a
certain ``plus space'' to modular forms. Doi and Naganuma found a
correspondence from modular forms to Hilbert modular forms, which can
be thought of as automorphic forms on the group $O_{2,2}(\R)$.  Maass
used a map from modular forms to automorphic forms on the group
$Sp_4(\R)$, or equivalently on $O_{2,3}(\R)$, in his work on the
Saito-Kurokawa correspondence [E-Z] (which is essentially the inverse
of Shimura's correspondence followed by Maass's
correspondence). Gritsenko [Gr] generalized the Maass correspondence to go
to automorphic forms on $O_{2,n}(\R)$.  Oda [O] and Rallis and
Schiffmann [R-S] had earlier used the Howe correspondence to construct
a map from automorphic forms on $SL_2$ to automorphic forms on
orthogonal groups which includes Gritsenko's correspondence.

In [B95] there is a generalization (in the level 1 case) of
Gritsenko's correspondence to the case when the modular form is
allowed to have singularities at cusps, so the automorphic form has
poles on rational quadratic divisors. In this section we will find a
similar extension of all of the correspondences above, which works for
forms which are allowed to have poles at cusps (and for all lattices
$M$ of signature $(2,b^-)$).

\proclaim Lemma 14.1. If  $A$ is an  integer then
$$
\sum_j (-1)^j{C\choose j}\times{A-2j+C-B-1\choose A-2j}
=\sum_j (-1)^j{C\choose A-j}{B\choose {j}}
$$
and in particular the sum on the left vanishes if
$B$ and $C$ are nonnegative integers and $B+C<A$.

Proof.
$$\eqalign{
&\sum_j (-1)^j{C\choose A-j}{B\choose {j}}\cr
=& \hbox{coefficient of $x^{A}$ in~} (1+x)^C(1-x)^B\cr
=& \hbox{coefficient of $x^{A}$ in~} (1-x^2)^C(1-x)^{B-C}\cr
=& \sum_j\hbox{coefficient of $x^{2j}$ in~} (1-x^2)^C
\times\hbox{coefficient of $x^{A-2j}$ in~} (1-x)^{B-C}\cr
=&\sum_j (-1)^j{C\choose j}\times(-1)^{A-2j}{B-C\choose A-2j}\cr
=&\sum_j (-1)^j{C\choose j}\times{A-2j+C-B-1\choose A-2j}\cr
}
$$
This proves lemma 14.1.

\proclaim Corollary 14.2. If $C$ and $m^+-h^+$ are integers such that
$0<C<m^+-h^+$
then
$$\sum_{j,m\atop j+m=C}{(-1)^j(m^+-h^++m-j-1)!\over j! (m^+-h^+-2j)!m!}=0.
$$

Proof. We put $B=0$ and $A=m^+-h^+$ in lemma 14.1 and find
$$
\sum_{j,m\atop j+m=C}{(-1)^j(m^+-h^++m-j-1)!C!\over j!
(m^+-h^+-2j)!m!(C-1)!}={C\choose m^+-h^+}.
$$
Corollary 14.2 follows immediately from this.

In the case of holomorphic forms the following theorem
is essentially contained in the results of [O], [R-S], and [Gr].

\proclaim Theorem 14.3. We use notation as in section 13.
Suppose $M$ is an even lattice of signature
$(2,b^-)$ and $F$ is an automorphic form of type $\rho_M$ of weight
$1+m^+-b^-/2$ as in section 6 which is holomorphic on $H$
and meromorphic at cusps.
Assume that $m^+ \ge 1$.
Then there is a
meromorphic function  $\Psi_M(Z_M,F)$
 with the following properties.
\item{1.} $\Psi_M$ is a meromorphic automorphic form
of weight $m^+$.
\item{2.} The only singularities of $\Psi_M$
are poles of order $m^+$ along divisors of the form $\lambda^\perp$ for
$\lambda\in M'$.
\item{3.} $\Psi_M(Z_M,F)=
{1\over 2}(-|Y_M|)^{-m^+}\Phi_M(Z_M/|Y_M|,p,F)$, where
$p(x^++x^-)=(ix_1^+-x_2^+)^{m^+}$.
\item{4.} If all the coefficients $c_\gamma(m)$ of $F$ vanish
for $m<0$ then $\Psi_M$ is a holomorphic function. If in addition $M$ has
dimension at least 5, or if $M$ has dimension 4 and contains no 2 dimensional
isotropic sublattice, then $\Psi_M$ is a holomorphic automorphic form.
\item{5.} If $z$ is a primitive norm 0 vector of $M$
and we choose notation as in sections 6 and 13 then for
sufficiently large $|Y|$ and $m^+>1$
the Fourier expansion
of $\Psi_z(Z,F)$ is given
by
$$\eqalign{
-&\sum_{\delta\in \Z/N\Z}
c_{\delta z}(0)
\sum_{0<\epsilon\le N}
N^{m^+-1}
\e(\delta\epsilon/N)
B_{m^+}(\epsilon/N)/2m^++\cr
+&
\sum_{n> 0}
\sum_{\lambda\in K'\atop(\lambda,W)> 0}
\e((n\lambda,Z))
n^{m^+-1}
\sum_{\delta\in M'/M\atop \delta|L=\lambda}
\e(n(\delta,z'))
c_\delta(\lambda^2/2)\cr
}
$$
and for $m^+=1$ it is given by the expression above plus the
constant function
$-\Phi_K(Y/|Y|,x_2^+,F_K)/ 2\sqrt 2 $.

Proof. We first suppose that $M$ has a primitive norm 0 vector $z$ and
work out the Fourier expansion of $\Phi_M$ at the cusp of $z$.  The
Fourier expansion of theorem 7.1 simplifies in the following ways.
We need  the formula
$$
K_{n+1/2}(z)=\sqrt{\pi/2z}\exp(-z)\sum_{0\le m\le n}
(2z)^{-m}{(n+m)!\over m!(n-m)!}
$$
valid for $n$ a nonnegative integer [E vol 2 section 7.2.6 formula (40)].
We would  like to extend it to be valid for $n=-1$, which
we can arrange by taking the sum over all $m\ge 0$ with the convention
that $(-1)!/(-1)!=1$.
Using 7.2 and substituting in this formula
we see that if $\lambda_{w^+}\ne 0$ then the integral over $y$ in 7.1
for $s=0$ is equal to
$$\eqalign{
&2c_\delta(\lambda^2/2)
\left({n\over 2|z_{v^+}||\lambda_{w^+}|}\right)^{-h^+-j+m^+-{1/2}}
K_{-h^+-j+m^+-1/2}(2\pi n|\lambda_{w^+}|/|z_{v^+}|)\cr
}
$$
(as $k=h=h^-=0$, $b^+=2$) which is equal to
$$\eqalign{
&2c_\delta(\lambda^2/2)
\left({n\over 2|z_{v^+}||\lambda_{w^+}|}\right)^{-h^+-j+m^+-{1/2}}
\times\cr&\times
\sqrt{\pi|z_{v^+}|/4\pi n|\lambda_{w^+}|}
\exp(-2\pi n|\lambda_{w^+}|/|z_{v^+}|)
\times\cr&\times
\sum_{0\le m}
(4\pi n|\lambda_{w^+}|/|z_{v^+}|)^{-m}
{(-h^+-j+m^+-1+m)!\over m!(-h^+-j+m^+-1-m)!}
.\cr
}
$$

We are given that if $\lambda\in K\otimes\R$ then
$$\eqalign{
p(v(\lambda))
&=i^{m^+}(\lambda,X/|Y|+iY/|Y|)^{m^+}\cr
=& \sum_{h^+} (\lambda,X)^{h^+}
{m^+\choose h^+}i^{2m^+-h^+}
|Y|^{-m^+}(\lambda,Y)^{m^+-h^+}\cr
=& \sum_{h^+} (\lambda,z_{v^+})^{h^+}
{m^+\choose h^+}i^{2m^+-h^+}
|Y|^{2h^+-m^+}(\lambda,Y)^{m^+-h^+}\cr
}
$$
so that
$$
p_{w,h^+,0}(w(\lambda))
={m^+\choose h^+}i^{2m^+-h^+}
|Y|^{h^+}(\lambda,Y/|Y|)^{m^+-h^+}.
$$

As a consequence we see that
$$\eqalign{
&(-\Delta)^j(\bar p_{w,h^+,0})(w(\lambda))\cr
=& (-1)^ji^{h^+-2m^+}{m^+\choose h^+}
{(m^+-h^+)!\over (m^+-h^+-2j)!}|Y|^{h^+}(\lambda,Y/|Y|)^{m^+-h^+-2j}.
\cr
}
$$

If we substitute these into the Fourier expansion of
$\Phi_M$ given in theorem 7.1 we find that $\Psi_M(Z_M,F)$
as defined in part 3 of 14.3 is given by
$$(-1)^{m^+}(P_1+P_2+P_3)/2$$
where
$$P_1= {|Y|^{-m^+}\over \sqrt 2|z_{v^+}|}\Phi_K(Y/|Y|,(ix_2^+)^{m^+},F_K)$$
is the term involving $\Phi_K$, $P_2$
is the sum of the terms with $\lambda=0$,
and
$$\eqalign{
&P_3=
{\sqrt2\over |z_{v^+}||Y|^{m^+}}
\sum_{n> 0}
\sum_{\lambda\in K'\atop\lambda\ne 0}
\e((n\lambda,\mu))
n^{h^+}
\sum_{\delta\in M'/M\atop \delta|L=\lambda}
\e(n(\delta,z'))
\sum_{h^+}i^{h^+-2m^+}{m^+\choose h^+}(2i)^{-h^+}
\times\cr&\times
\sum_j
{(-1)^j\over j!(8\pi)^j}{(m^+-h^+)!\over (m^+-h^+-2j)!}
|Y|^{h^+}(\lambda,Y/|Y|)^{m^+-h^+-2j}
\times\cr&\times
2c_\delta(\lambda^2/2)
\left({n\over 2|z_{v^+}||\lambda_{w^+}|}\right)^{-h^+-j+m^+-{1/2}}
\times\cr&\times
\sqrt{\pi|z_{v^+}|/4\pi n|\lambda_{w^+}|}
\exp(-2\pi n|\lambda_{w^+}|/|z_{v^+}|)
\times\cr&\times
\sum_{0\le m}
(4\pi n|\lambda_{w^+}|/|z_{v^+}|)^{-m}
{(-h^+-j+m^+-1+m)!\over m!(-h^+-j+m^+-1-m)!}
\cr
}
$$
is the sum of the terms with $\lambda\ne 0$.

Fortunately most of the terms in $P_3$ cancel. If we fix
$h^+$ and a value of of $C=j+m$ and  compare the expression
for $P_3$ with corollary 14.2 and use the fact that
$$
{1\over (8\pi)^j}
(4\pi n|\lambda_{w^+}|/|z_{v^+}|)^{-m}
(\lambda,Y/|Y|)^{-2j}
\left({n\over 2|z_{v^+}||\lambda_{w^+}|}\right)^{-j}
=
\left({|z_{v^+}|\over 4\pi n |\lambda_{w^+}|}\right)^{m+j}
$$
depends on $m$ and $j$ only through $m+j$
we see that all the terms in $P_3$
with $C>0$, or equivalently with $m\ne 0$ or $j\ne 0$, cancel out.
We then find that the only factors involving $h^+$
are of the form
$$\eqalign{
&\sum_{h^+}(2i)^{-h^+}{m^+\choose h^+}
{n^{h^+}|Y|^{h^+}(\lambda,Y/|Y|)^{-h^+}i^{h^+}
\over n^{h^+}(|2z_{v^+}||\lambda_{w^+}|)^{-h^+}}\cr
=& \sum_{h^+}{m^+\choose h^+}
\left({(Y,\lambda)\over |(Y,\lambda)|}\right)^{-h^+}
\cr
}
$$
which is equal to 0
if $(\lambda,Y)<0$ (because $m^+>0$) and to
$2^{m^+}$ if $(\lambda,Y)>0$.
If we put these simplifications into the expression for $P_3$
we find
$$\eqalign{
P_3&=
{\sqrt2 \over |z_{v^+}||Y|^{m^+}}
\sum_{n> 0}
\sum_{\lambda\in K'\atop(\lambda,W)>0}
\e((n\lambda,\mu))
\sum_{\delta\in M'/M\atop \delta|L=\lambda}
\e(n(\delta,z'))(-1)^{m^+}
\times\cr&\times
(\lambda,Y/|Y|)^{m^+}
c_\delta(\lambda^2/2)
\times\cr&\times
2\left({n\over 2|z_{v^+}||\lambda_{w^+}|}\right)^{m^+-{1/2}}
\times\cr&\times
\sqrt{\pi|z_{v^+}|/4\pi n|\lambda_{w^+}|}
\exp(-2\pi n|\lambda_{w^+}|/|z_{v^+}|)
\cr
=&
2 (-1)^{m^+}
\sum_{n> 0}
\sum_{\lambda\in K'\atop(\lambda,W)> 0}
\e((n\lambda,\mu))
n^{m^+-1}
\sum_{\delta\in M'/M\atop \delta|L=\lambda}
\e(n(\delta,z'))
\times\cr&\times
c_\delta(\lambda^2/2)
\exp(-2\pi n|(\lambda,Y)|)\cr
=&
2 (-1)^{m^+}
\sum_{n> 0}
\sum_{\lambda\in K'\atop(\lambda,W)> 0}
\e(n(\lambda,X+iY))
n^{m^+-1}
\sum_{\delta\in M'/M\atop \delta|L=\lambda}
\e(n(\delta,z'))
c_\delta(\lambda^2/2).
}
$$

We work out the constant term $P_2$ as follows.
Note that for $\Re(s)$ large the expression for the terms with
$\lambda=0$ is the limit of the expression for some nonzero $\lambda$ as
$\lambda$ tends to 0. Hence we can work out $P_2$ by
taking the expression for $P_3$, changing $n^{m^+}-1$ to $n^{m^+-1-2s}$,
and taking the constant term at $s=0$ of its analytic continuation.
If we do this we find that
$$\eqalign{
P_2&=
2(-1)^{m^+}
\sum_{\delta\in M'/M\atop \delta|L=0}
c_\delta(0)
\sum_{n> 0}
n^{m^+-1-2s}
\e(n(\delta,z'))\cr
&=
2(-1)^{m^+}
\sum_{\delta\in \Z/N\Z}
c_{\delta z}(0)
\sum_{0<\epsilon\le N}
\sum_{n\ge 0}
(Nn+\epsilon)^{m^+-1-2s}
\e(\delta\epsilon/N)\cr
&=
2(-1)^{m^+}
\sum_{\delta\in \Z/N\Z}
c_{\delta z}(0)
\sum_{0<\epsilon\le N}
N^{m^+-1}
\e(\delta\epsilon/N)
\sum_{n\ge 0}
(n+\epsilon/N)^{m^+-1-2s}\cr
&=
2(-1)^{m^+}
\sum_{\delta\in \Z/N\Z}
c_{\delta z}(0)
\sum_{0<\epsilon\le N}
N^{m^+-1}
\e(\delta\epsilon/N)
(-B_{m^+}(\epsilon/N)/m^+)\cr
}
$$
by [E, 1.10, formulas (1) and (11)].

We work out the term $P_1$.  We know that $$P_1= (|Y|^{-m^+}/ \sqrt
2|z_{v^+}|)\Phi_K(Y/|Y|,(x_2^+)^{m^+},F_K).$$ By theorem 10.3 the
function $\Phi_K(*,(x_2^+)^{m^+},F)$ is a polynomial of degree
$m^--m^++2k_{max}+1=1-m^+$ and so $P_1$ is zero if $m^+>1$.  If
$m^+=1$ then we see by theorem 10.3 that $\Phi_K$ is a constant, and
as $|Y||z_{v^+}|=1$ we see that $P_1$ is ${\Phi_K(Y/|Y|,x_2^+,F_K)/
\sqrt 2 }$.

If we add together the expressions we have found for
$P_1$, $P_2$, and $P_3$ we obtain the Fourier expansion
of theorem 14.3.

In the case when $M$ has no primitive norm 0 vector $z$
we can show that $\Psi$ is holomorphic by using
the embedding trick in section 8 as in the proof of 13.3.

We can work out the singularities of $\Psi_M$ directly from theorem 6.2.
The function $(ix_1-x_2)^{m^+}$ is harmonic, so by 6.2 we find that
the function $\Psi_M$ has a singularity of type
$$\eqalign{
&{(-1)^{m^+}\over 2}\sum_{\lambda\in M'\cap v_0^\perp\atop \lambda\ne 0}
c_\lambda(\lambda^2/2)(i\lambda,\bar Z/|Y|)^{m^+}
(2\pi |(\lambda,Z/|Y|)|^2)^{-m^+}\Gamma(m^+)\cr
=&
\sum_{\lambda\in M'\cap v_0^\perp\atop (\lambda,W)>0}
{c_\lambda(\lambda^2/2)(m^+-1)!
\over
(2\pi i (\lambda,Z/|Y|))^{m^+}
}
\cr
}
$$
so in particular we see that the singularities are all poles of order
exactly $m^+$
along rational quadratic divisors.

This proves theorem 14.3.

{\bf Example 14.4.} We will work out a case of the singular Shimura
correspondence explicitly.  We restrict to the case when $F$ has
weight $1/2+m^+$ for $m^+$ even and assume that $F$ has type $\rho_K$
where $K$ is a one dimensional lattice generated by a vector of norm
2. Such modular forms are equivalent to modular forms of level 4
satisfying Kohnen's plus space condition; see [E-Z] chapter 5 or [K].
The Shimura correspondence takes forms of weight $m^++1/2$ to forms of
weight $2m^+$, but the correspondence above in the case $b^-=1$ takes
forms of weight $m^++1/2$ to forms of weight $m^+$. The reason for
this factor of 2 between the two weights is that we are constructing
forms on $O_{2,1}(\R)$, and the map from $SL_2(\R)$ to the identity
component of $O_{2,1}(\R)$ is a double cover. Hence we pick up a
factor of 2 in the weights (which are essentially representations of
the maximal torus) when we go from $O_{2,1}(\R)$ to $SL_2(\R)$.
Theorem 14.3 then implies that if $f(\tau)=\sum c(n)q^n$ is a modular
form for $\Gamma_0(4)$ of weight $m^++1/2$ such that $c(n)$ vanishes
unless $n\equiv 0,1\bmod 4$ then
$$\Psi_M(\tau)= {-c(0)B_{m^+}\over2m^+}+\sum_n \sum_m q^{mn}n^{m^+-1} c(m^2)
$$
is a modular form of weight $2m^+$. When $f$ is holomorphic
this is a special case of theorem 1 of [K]. Theorem 14.3 says
that it is still correct even if $f$ has poles at the cusps,
although the function $\Psi_M$ will then have poles of order $m^+$
at some quadratic irrationals.

For our explicit example of this we
want $F(\tau)$ to be a weight $5/2$ modular form of level 4 of the
form $q^{-3}+O(q)$ satisfying Kohnen's plus space condition,
so we define
$F(\tau)$ by
$$\eqalign{
E(\tau)&=\sum_{n>0,n \,{\rm odd}} \sigma_1(n)q^n = q+4q^3+6q^5\cdots\cr
\theta(\tau)&=\sum_{n\in Z} q^{n^2} = 1+2q+2q^4+\cdots\cr
F(4\tau)&=E(\tau)\theta(\tau)(\theta(\tau)^4-2E(\tau))
(\theta(\tau)^4-16E(\tau))E_8(4\tau)/\Delta(4\tau)+6720\sum_nH(2,n)q^n\cr
    &=q^{-3}+64q-32384q^4+131535q^5-4257024q^8+11535936q^9+O(q^{12})\cr
&=\sum_n c(n)q^n\cr
}$$
(where $H(2,n)$ is Cohen's function [Co])
so that $F$ has weight $5/2$. Applying theorem 14.3
we see that $\Psi_M(\tau,F)$ is a modular form
of weight $2(5/2-1/2)=4$, with a singularity of type
$(2\pi)^{-2}3^{-1}(\tau-\omega)^{-2}$ at a cube root of unity $\omega$
and a zero at the cusp $i\infty$.
Hence we find that $\Psi_M(\tau)$ must be
$$\eqalign{
64\Delta(\tau)/E_4(\tau)^2
&= 64(q - 504q^2 + 180252q^3 - 56364992q^4 +O(q^5))\cr
&=\sum_nb(n)q^n.\cr
}$$
The Fourier expansion in theorem 14.3 states that
$b(n)=\sum_{d|n}dc(n^2/d^2)$, which can be checked explicitly
in the example above.

The function $\Delta(\tau)/E_4(\tau)^2$ can also be written as an infinite
product whose exponents are given by coefficients of
a modular form of weight $1/2$ with a pole at the cusp; see
theorem 14.1 and the examples following it in [B95].

More generally, the classical Shimura correspondence works well
for forms of weight $m^++1/2$ at least $5/2$, but behaves strangely
for weight $3/2$ (the images of cusp forms need not be cusp forms)
and very badly for weight $1/2$. We see that this odd behavior
in low weights is caused by the term involving $\Phi_K$, which
is a piecewise polynomial of degree at most $1-m^+$, so it vanishes
for weights at least $5/2$. For weight $3/2$ it adds an extra constant
(so the image of a cusp form need not be a cusp form)
and in weight $1/2$ it adds a linear term, which is essentially
the Weyl vector in theorem 13.3. In the case considered by
Maass and Gritsenko with $b^->2$ and holomorphic functions $F$,
we see that the term involving $\Phi_K$ always vanishes
except when $b^-=3$ and $F$ has weight $1/2$.

\proclaim 15.~Examples related to mirror symmetry and Donaldson polynomials.

{\bf Example 15.1.} Take $M$ to be the lattice $II_{3,19}$ and take $F$ to be
$E_4(\tau)/\Delta(\tau)= q^{-1} + 264 + 8244q + 139520q^2 + O(q^3)$.
Then the function $\Phi_M(v,1,F)$ is a function on the Grassmannian
$G(M)$ invariant under $\Aut(II_{3,19})$ whose only singularities are
on the subspaces of the form $r^\perp$ for $r^2=-2$. Recall that the
period space of Ricci flat metrics of volume 1 on a marked K3 surface
is exactly the set of points where this function $\Phi_M$ is
nonsingular.  (See [B-P-V, chapter 8, sections 11-14].) Hence the
function $\Phi_M$ can be thought of as a function on the moduli space of
K3 surfaces with a Ricci flat metric. Todorov and Jorgenson [J-T] have
announced the construction of  similar functions by
taking a regularized determinant of a Laplacian operator on a K3
surface; it seems natural to conjecture that the
function  $\Phi_M(v,1,F)$ can be constructed in the same way.
See also [HM96].

{\bf Example 15.2.} Similarly we can take $M$ to be the
lattice $II_{4,20}$ and again
take $F$ to
be $E_4(\tau)/\Delta(\tau)$.
Then we get a function $\Phi_M$ on the moduli space of ``K3 surfaces
with a B-field modulo mirror symmetry'', which is (more or less) the
quotient by $\Aut(II_{4,20})$ of the subset of the Grassmannian
$G(M)$ of points not orthogonal to a norm $-2$ vector of $M$.  See [A-M]
for more details.

{\bf Example 15.3.} There seems to be some connection between the
automorphic forms with singularities on hyperbolic space (which are
piecewise polynomials by 10.3) and Donaldson polynomials for
4-manifolds with $b^+=1$. If we can find an automorphic form with
singularities with the same wall crossing formula as a Donaldson
piecewise polynomial invariant then we can subtract them to obtain a
polynomial invariant of the 4-manifold not depending on the choice of
chamber. We will give an example where it is possible to do this.

In [D] Donaldson
defines an invariant for 4-manifolds with $b^+=1$
which is essentially
a piecewise linear function on the space $H^2\otimes \R$ (where
$H^2$ is the second homology group of the manifold).
Define the Weyl chambers to be the components of positive norm vectors
which are not orthogonal to any norm $-1$ vector.
Donaldson shows that his invariant is given by the inner product
by a fixed vector $\rho(W)$ on the interior
of each Weyl chamber $W$, and satisfies the wall crossing formula
$\rho(W_1)=\rho(W_2)+2r$ if $r$ is a norm $-1$ vector such that
$r^\perp$ is the wall between $W_1$ and $W_2$ and which has positive inner
product with $W_1$.

Suppose that the second homology group is isometric to $I_{1,b^-}$,
and suppose that $b^-\le 9$. Then there is a Weyl vector for the
reflection group generated by norm $-1$ vectors (given by the Weyl
vector of example 13.7 when $b^-=9$), and these Weyl vectors satisfy
the same wall crossing formulas as the Weyl vectors defining
Donaldson's invariant (up to a factor of 2). Hence if we subtract the
singular automorphic form corresponding to these Weyl vectors from
Donaldson's invariant we get a function whose wall crossing formulas
are all 0 and is therefore a polynomial. So if $b^-\le 9$ we can find
an invariant of the manifold given by a linear function which does not
depend on the choice of Weyl chamber.  This is in some sense a sort of
average of the Donaldson polynomials of the different Weyl chambers.

Donaldson worked out his invariant in the cases when the manifold is
either $\P^2(\C)$ blown up at 9 points, or a Dolgachev surface (when
$b^-=9$). Donaldson's results imply that in the first case the
invariant polynomial is 0 and in the second case the invariant
polynomial is nonzero. (Donaldson used this to show that the two
manifolds are not diffeomorphic even though they are homeomorphic.)

When $b^-\ge 10$ Donaldson's piecewise linear function does not always
seem to be the same as one of the automorphic forms with singularities
in this paper.

\proclaim 16.~Open problems.

{\bf Problem 16.1.}
In theorem 12.1 we give a sufficient condition for
a Lorentzian lattice to have a reflection group of finite or virtually
free abelian index in its automorphism  group. Is this also a necessary
condition? Can it be used to classify the Lorentzian lattices with this
property? (Nikulin showed that the number of such lattices
is essentially finite.)

{\bf Problem 16.2.}
A closely related question is that of finding all ``interesting''
generalized Kac-Moody algebras. It is not quite clear what
``interesting'' should mean, but it should certainly include cases
when the denominator function is an automorphic form of singular
weight, and possibly all cases when the denominator function is an
automorphic form. These appear to correspond roughly to cases when the
Lorentzian lattice $M$ has a reflection subgroup with a norm 0 Weyl
vector. (However there are also many cases when the generalized
Kac-Moody algebra has no real roots so does not obviously correspond to some
reflection group acting on $M$.) Gritsenko and Nikulin [G-N] have recently
written several preprints giving many examples of automorphic forms
related to generalized Kac-Moody superalgebras.

{\bf Problem 16.3.}  If we take the lattice to be of signature (2,1)
or (2,2) then we get lots of examples of meromorphic sections of line
bundles over modular curves or Hilbert modular surfaces with known
zeros and poles from theorem 13.3. More generally we
can get meromorphic functions on some higher dimensional varieties
in the same way. Can these be used to give
interesting relations between elements of the Picard or N\'eron-Severi
groups represented by the divisors of zeros of these sections?
In particular is it possible to prove the Gross--Zagier
theorem along these lines?

{\bf Problem 16.4.}
We have worked throughout with quadratic forms over the integers.
It seems natural to ask if everything can be extended to quadratic forms over
the rings of integers of
algebraic number fields and function fields.
There is one obvious major problem in carrying out any extension:
because of the Koecher boundedness principle, holomorphic automorphic forms
with poles at cusps are rare in higher dimensions. For example,
they do not exist for $SL_2$ of any totally real number field other than
the rational numbers, so it is unclear how to extend the results to Hilbert
modular varieties. However it may be possible to do something
with $SL_2$ of the integers of a real quadratic field
or a quaternion algebra over the rational numbers, when the corresponding
symmetric space is hyperbolic space of dimension 3 or 5.

{\bf Problem 16.5.}
Describe how the correspondence in this paper behaves under the action
of Hecke operators. (This is another place where it would probably be
easier to use the Weil representation over the adeles.) In the
holomorphic cases correspondences such as the Shimura correspondence [Sh]
take eigenforms to eigenforms, but for forms with singularities this
statement is usually vacuous as eigenforms do not exist in general. A
possible replacement for this might be that there is a homomorphism
from the Hecke algebra of $O_{b^+,b^-}$ to that of $SL_2$ which is
compatible with the correspondences in theorems 13.3 and 14.3.
(In the case of theorem 13.3 the Hecke operators for
$O_{2,b^-}$ should act multiplicatively rather than additively
on the meromorphic infinite products.)

{\bf Problem 16.6.}
What local properties do the functions $\Phi_M$ have?
Are they eigenfunctions of the Laplacian at their nonsingular points?
More generally, are their restrictions to nonsingular values killed
by an ideal of finite index in the center of the universal enveloping
algebra of $O_M(\R)$, or in other words do they satisfy the same local
conditions as
automorphic forms except at their singular points?
This may follow from the explicit Fourier series expansion,
most of whose terms look like eigenfunctions of the Laplacian.
Kontsevich recently  told me that he has calculated the behavior
of $\Phi_M$ under differential operators; see [Kon, section 3.3].

{\bf Problem 16.7.}
When $b^+=1$ the wall crossing formula and the fact that $\Phi_M$ is a
piecewise polynomial are remarkably similar to statements about
Donaldson invariants of 4 manifolds with $H^2$ equal to $I_{1,b^-}$.
When do Donaldson polynomials for
some 4 manifolds have  the same wall crossing formulas as some functions
$\Phi_M$, as in example 15.3? When they do  the piecewise
polynomial Donaldson invariants split as the sum of a polynomial
invariant and an automorphic form $\Phi_M$ with singularities, which
gives polynomial invariants for 4-manifolds with $b^+=1$ not
depending on a choice of Weyl chamber. Example
15.3 shows that this happens for piecewise linear Donaldson invariants
when $b^-\le 9$. The wall crossing formulas for more general
4-manifolds with $b^+=1$ are given by G\"ottsche in [G, theorem 3.3] and are
similar to the wall crossing formulas in this paper; for example,
both wall crossing formulas are polynomials in the quadratic form and
a linear function vanishing on the wall.  When $b^+>1$ it is harder to
see a possible connection, because the Donaldson polynomials are
polynomials on $M\otimes \R$, while the functions $\Phi_M$ are neither
polynomials nor defined on $M\otimes \R$.

{\bf Problem 16.8.}  Investigate the functions on other Hermitian
symmetric spaces.  We only get infinite products which are holomorphic
automorphic forms on Hermitian symmetric spaces with hermitian
symmetric subspaces of complex dimension 1 less, in other words the
symmetric spaces of $O_{2,b^-}(\R)$ and $U(1,n)$, and the symmetric
spaces of $U(1,n)$ (which are the unit balls in $\C^{n}$) can be
embedded in the symmetric spaces of $O_{2,2n}$. On other symmetric
spaces, such as Siegel upper half planes of genus $g$ greater than 2,
we do not get holomorphic automorphic forms as infinite products, but
we do get real analytic automorphic forms with singularities along
Siegel upper half planes of genus $g-1$ by embedding the Siegel upper
half plane in the Grassmannian $G(\R^{2g,2g})$ and restricting an
automorphic form with singularities on this Grassmannian. What can be
done with these? Is is possible to find holomorphic sections of vector
bundles on Siegel upper half planes with known zeros?

{\bf Problem 16.9.}
What congruence conditions (especially at the primes $2$ and $3$
dividing $|M'/M|$) does the Weyl vector satisfy when $F$
has integral coefficients? What are the ``best possible'' congruences
satisfied by lattices, or in other words what is the lattice
generated by the theta functions of lattices of some fixed genus?

{\bf Problem 16.10.} Can one reverse the correspondence from
modular forms to automorphic forms with singularities,  and
reconstruct modular forms from automorphic forms with singularities?
In particular when are there isomorphisms between spaces?

{\bf Problem 16.11.} In [H] Hejhal constructs some ``pseudo cusp forms''
which are functions on the upper half plane with logarithmic
singularities at imaginary quadratic irrationals, which are almost
eigenfunctions of the Laplacian with eigenvalues given
by the imaginary part of zeros of the Riemann zeta function.
Can these pseudo cusp forms
be constructed using  the singular Howe correspondence in this paper?
Some good candidates for constructing these pseudo cusp forms
using the Howe correspondence might be the functions mentioned in problem
16.13.

{\bf Problem 16.12.} Generalize the correspondence of theorem 14.3 to vector
valued forms by using polynomials with $m^->0$.

{\bf Problem 16.13.}
What happens if theorem 7.1 is applied to functions which
are not almost holomorphic?
For example Freitag asked if there are analogues of Maass wave forms
with singularities at cusps, in which case these could be inserted into
theorem 7.1. There are many examples of such functions
in [H83]; for example, the functions $F_n$ on page 658.

\proclaim References.

\item{[A-M]} P. Aspinwall, D. Morrison,  String theory on K3 surfaces,
preprint hep-th/9404151
\item{[B-P-V]} W. Barth, C. Peters, A. Van de Ven,
``Compact complex surfaces'', Springer Verlag, 1984.
\item{[B92]}{R. E. Borcherds,
     Monstrous moonshine and monstrous Lie superalgebras,
     Invent. Math. 109, 405-444 (1992).}
\item{[B95]} R. E. Borcherds,  Automorphic forms on
     $O_{s+2,2}(R)$ and infinite products.
     Invent. Math. 120, p. 161-213 (1995)
\item{[B96]}R. E. Borcherds,  The moduli space of Enriques surfaces and the
     fake monster Lie superalgebra. Topology vol. 35 no. 3, 699-710, 1996.
\item{[B-C]} A. Borel, W. Casselman, Automorphic forms, representations, and
$L$-functions, Proc. Symp. Pure Math. Vol. XXXIII, A.M.S.1979.
\item{[Co]} H. Cohen, Sums involving the values at negative integers of
   $L$-functions of quadratic characters. Math. Ann. 217, 271-285 (1975).
\item{[C]}{J. H. Conway, The automorphism group of the
   26 dimensional even Lorentzian lattice. J. Algebra 80 (1983) 159-163.
   This paper is reprinted as chapter 27 of [C-S].
        }
\item{[C-S]}{J. H. Conway, N. J. A. Sloane.  Sphere packings, lattices
   and groups. Springer-Verlag New York 1988, Grundlehren der mathematischen
   Wissenschaften 290.}
\item{[D]} S. Donaldson, Irrationality and the $h$-cobordism conjecture,
J. Differential Geometry 26 (1987) 141-168.
\item{[E-Z]}{M. Eichler, D. Zagier, ``The theory of Jacobi forms'', Progress
in mathematics vol. 55,
Birkh\"auser, Boston, Basel, Stuttgart 1985.}
\item{[E]}{A. Erdelyi, W. Magnus, F. Oberhettinger, F. G. Tricomi,
``Higher transcendental functions'', vols. 1-3, and ``Tables
of integral transforms'', vols 1-2,
McGraw-Hill Book Company Inc, New York, Toronto, London 1953.}
\item{[F-F]} I. Frenkel, A. Feingold, A hyperbolic Kac-Moody algebra
and the theory of Siegel modular forms of genus 2, Math. Ann. 263
(1983) no. 1, 87-144.
\item{[G]} L. G\"ottsche, Modular forms and Donaldson invariants for
   4-manifolds with $b_+ =1$, Journal of the A. M.S., Vol. 9, No. 3,
   July 1966, 827-843.  Also see the preprint alg-geom/9612020 by
   L. G\"ottsche and D. Zagier, ``Jacobi forms and the structure of
   Donaldson invariants for 4-manifolds with $b^+=1$''.
\item{[Gr]} V. A. Gritsenko, Jacobi functions of $n$ variables.
 J. Soviet Math. 53 (1991), 243-252.
(The original Russian is in Zap. Nauch. Seminars of LOMI v. 168 (1988),
p. 32-45.)
\item{[G-N]} V. A. Gritsenko, V. V. Nikulin, Siegel automorphic form
corrections of some Lorentzian
Kac-Moody Lie algebras. Amer. J. Math. 119 (1997), no. 1, 181--224.
Automorphic correction of a Lorentzian Kac-Moody
algebra. C. R. Acad. Sci. Paris S\'er. I Math. 321 (1995), no. 9, 1151--1156.
Also see the preprints by these authors on the alg-geom preprint server
http://xxx.lanl.gov/archive/alg-geom/
\item{[H-M]} J. Harvey, G. Moore,
Algebras, BPS states, and strings. Nuclear Phys. B 463 (1996),
no. 2-3, 315--368.
preprint hep-th/9510182.
\item{[HM96]} J. Harvey, G. Moore, Exact Gravitational Threshold Correction
in the FHSV Model, preprint hep-th/9611176.
\item{[H]} D. Hejhal, Some observations concerning eigenvalues of the
 Laplacian and Dirichlet $L$-series, in H. Halberstam and C. Hooley (Eds.),
``Recent Progress in Analytic Number Theory''. Academic, London, 1981, 95-110.
\item{[H83]} D. Hejhal, The Selberg trace formula for $PSL(2,R)$ Vol. 2.
Springer Lecture notes in Mathematics vol. 1001, 1983.
\item{[J-T]} J. Jorgenson, A. Todorov, Enriques surfaces
and analytic discriminants, preprint 1995. An analytic discriminant
for polarized K3 surfaces, 1994 preprint.
\item{[Ko]} W. Kohnen, Modular forms of half integral weight on $\Gamma_0(4)$,
Math. Ann. 248, 249-266 (1980).
\item{[Kon]} M. Kontsevich, Product formulas for modular forms on $O(2,n)$.
S\'eminaire Nicolas Bourbaki 821, November 1996.
\item{[L-S-W]} W. Lerche, A. N. Schellekens, N. P. Warner, ``Lattices
and strings'', Phys. Representation. 177 (1989) 1.
\item{[Ni]} S. Niwa, Modular forms of half integral weight and the integral of
certain theta functions, Nagoya Math J. 56 (1975), 147-163.
\item{[O]} T. Oda,
On modular forms associated with indefinite quadratic forms of
signature $(2,n-2)$, Math. Ann. {\bf 231} (1977), 97--144.
\item{[R-S]} S. Rallis and G. Schiffmann,
On a relation between $SL_2$ cusp forms
and cusp forms on tube domains associated to orthogonal groups,
Trans. AMS {\bf 263} (1981), 1--58.
\item{[R]} U. Ray, A character formula for generalized Kac-Moody
superalgebras. J. Algebra 177 (1995), no. 1, 154--163.
\item{[Sh]} G. Shimura,
On modular forms of half integral weight. Ann. of Math. (2) 97
(1973), 440--481.
\item{[S]} T. Shintani, On construction of holomorphic cusp forms
of half integral weight,
Nagoya Math. J. {\bf 58} (1975), 83--126.
\item{[V-K]} Vinberg, Kaplinskaja, On the groups $O_{18,1}(\Z)$ and
$O_{19,1}(\Z)$, Soviet Math. 19, No 1 (1978) 194-197.
\item{[W]} A. Weil, Sur certains groupes d'op\'erateurs unitaires,
Acta Math. 111, pp. 143-211.
\item{[Z]} D. Zagier, Nombres de classes et formes modulaires de poids $3/2$,
 C. R. Acad. Sci. Paris S\'er. A-B 281 (1975), no. 21, Ai, A883--A886.
\bye